\documentclass[12pt]{article}
\usepackage{color}
\usepackage{amsthm}
\usepackage{amsmath}
\usepackage{amssymb}
\usepackage{enumerate}
\usepackage{pgf,tikz}\usetikzlibrary{matrix, calc, arrows,decorations.pathreplacing}
\usepackage{cite}
\usepackage{pdfcolfoot}

\usepackage{accents}
\DeclareRobustCommand
{\mathringbig}[1]{\accentset{\smash{\raisebox{-0.1ex}{$\scriptstyle\circ$}}}{#1}\rule{0pt}{2.3ex}}
\usepackage{color}
\usepackage{xcolor}  

\newcommand{\mypropbox}[1]{\vspace{1mm}\par\noindent\fbox{\parbox{0.98\textwidth}{\vspace{-2mm}#1\vspace{-2mm}}}\par\vspace{1mm}}

\usepackage[bookmarks=true, bookmarksopen=true, bookmarksopenlevel=5]{hyperref} 

\definecolor{myblue}{rgb}{0,0,0.6}     
\hypersetup{pdftitle={Density results for Sobolev, Besov and Triebel-Lizorkin spaces on rough sets},  
            pdfauthor={A. M. Caetano, D. P. Hewett, A. Moiola},
     colorlinks=true, linkcolor=myblue,  citecolor=myblue, filecolor=myblue,   urlcolor=myblue,  }
\allowdisplaybreaks[4]

\newcommand{\rf}[1]{(\ref{#1})}
\newcommand{\mmbox}[1]{\fbox{\ensuremath{\displaystyle{ #1 }}}}
\newcommand{\cS}{\mathcal{S}}
\newcommand{\cD}{\mathcal{D}}
\newcommand{\cH}{\mathcal{H}}
\newtheorem{thm}{Theorem}[section]
\newtheorem{lem}[thm]{Lemma}
\newtheorem{defn}[thm]{Definition}
\newtheorem{prop}[thm]{Proposition}
\newtheorem{cor}[thm]{Corollary}
\newtheorem{rem}[thm]{Remark}
\newcommand{\tH}{\widetilde{H}}
\newcommand{\dimH}{{\rm dim_H}}
\newcommand{\OO}{{(\Omega)}}
\newcommand{\baro}{{\overline{\Omega}}} 
\newcommand{\Tr}{{\mathrm{Tr}}}
\newcommand{\IH}{{\mathbb H}}
\newcommand{\GG}{{(\Gamma)}}
\DeclareMathOperator{\supp}{supp}
\DeclareMathOperator{\diam}{diam}
\DeclareMathOperator{\dist}{dist}
\newcommand{\subdense}{{\,\overset{\mathrm{ds}}{\subset}\,}}
\newcommand{\subcl}{{\,\overset{\mathrm{cl}}{\subset}\,}}
\newcommand{\supcl}{{\,\overset{\mathrm{cl}}{\supset}\,}}

\begin{document}
 \global\long\def\N{\mathbb{N}}
 \global\long\def\No{\mathbb{N}_{0}}
 \global\long\def\R{\mathbb{R}}
 \global\long\def\Z{\mathbb{Z}}
 \global\long\def\C{\mathbb{C}}
 \global\long\def\Nn{\mathbb{N}^{n}}
 \global\long\def\Non{\mathbb{N}_{0}^{n}}
 \global\long\def\Rn{\mathbb{R}^{n}}
 \global\long\def\Zn{\mathbb{Z}^{n}}
 \global\long\def\Cn{\mathbb{C}^{n}}
 \global\long\def\K{\mathbb{K}}
 \global\long\def\Bspq{B_{p,q}^{s}(\mathbb{R}^{n})}
 \global\long\def\Fspq{F_{p,q}^{s}(\mathbb{R}^{n})}
 \global\long\def\Aspq{A_{p,q}^{s}(\mathbb{R}^{n})}
 \global\long\def\Aspqp{A_{p',q'}^{-s}(\mathbb{R}^{n})}
 \global\long\def\tildeA{\widetilde{A}_{p,q}^{s}(\Omega)}
 \global\long\def\tildeB{\widetilde{B}_{p,q}^{s}(\Omega)}
 \global\long\def\tildeF{\widetilde{F}_{p,q}^{s}(\Omega)}
 \global\long\def\subA{A_{p,q,\overline{\Omega}}^{s}}
 \global\long\def\subB{B_{p,q,\overline{\Omega}}^{s}}
 \global\long\def\subF{F_{p,q,\overline{\Omega}}^{s}}
 \global\long\def\subAg{A_{p,q,\Gamma}^{s}}
 \global\long\def\tildeAg{\widetilde{A}_{p,q}^{s}(\Gamma^{c})}
 \global\long\def\subAgp{A_{p',q',\Gamma}^{-s}}
 \global\long\def\tildeAgp{\widetilde{A}_{p',q'}^{-s}(\Gamma^{c})}

\newcommand{\TAO}{R\!A^s_{p,q,\baro}}
\newcommand{\TAOO}{T\!\widetilde A^s_{p,q}(\overline\Omega)}
\newcommand{\SqS}{{\Gamma}}
\newcommand{\SNangle}{\beta}
\newcommand{\SN}{\Gamma}
\newcommand{\LLambda}{\Gamma}
 
\title{Density results for Sobolev, Besov and \\ Triebel--Lizorkin spaces on rough sets}

\author{A. M. Caetano$^{\text{a}}$\footnote{A.C.\ acknowledges the support of CIDMA (Center for Research and
Development in Mathematics and Applications) and FCT (Foundation for Science and Technology) within projects UID/MAT/04106/2019 and UIDB/04106/2020.},
D. P. Hewett$^{\text{b}}$\footnote{D.H. acknowledges support from EPSRC grant EP/S01375X/1.} 
and A. Moiola$^{\text{c}}$\footnote{A.M.\ acknowledges support from EPSRC grant EP/N019407/1, from GNCS-INDAM
and MIUR through the ``Dipartimenti di Eccellenza'' Programme (2018--2022) -- Dept. of Mathematics, University of Pavia.}
\\[2pt]
$^{\text{a}}${\footnotesize
Center for R\&D in Mathematics and Applications},\\ {\footnotesize Departamento de Matem\'atica, Universidade de Aveiro, Aveiro, Portugal}\\[4pt]
$^{\text{b}}${\footnotesize Department of Mathematics, University College London, London, United Kingdom}\\[4pt]
$^{\text{c}}${\footnotesize Dipartimento di Matematica ``F. Casorati'', Universit\`a degli studi di Pavia, Pavia, Italy}
}

\maketitle

\begin{abstract}
We investigate two density questions for Sobolev, Besov and Triebel--Lizorkin spaces on rough sets. Our main results, stated in the simplest Sobolev space setting, are that: (i) for an open set $\Omega\subset\Rn$, $\mathcal{D}(\Omega)$ is dense in $\{u\in H^s(\Rn):\supp{u}\subset \overline{\Omega}\}$ whenever $\partial\Omega$ has zero Lebesgue measure and $\Omega$ is ``thick'' (in the sense of Triebel); and (ii) 
for a $d$-set $\Gamma\subset\Rn$ ($0<d<n$), $\{u\in H^{s_1}(\Rn):\supp{u}\subset \Gamma\}$ is dense in $\{u\in H^{s_2}(\Rn):\supp{u}\subset \Gamma\}$ whenever $-\frac{n-d}{2}-m-1<s_{2}\leq s_{1}<-\frac{n-d}{2}-m$ for some $m\in\N_0$. 
For (ii), we provide concrete examples, for any $m\in\N_0$, where density fails when $s_1$ and $s_2$ are on opposite sides of $-\frac{n-d}{2}-m$.
The results (i) and (ii) are related in a number of ways, including via their connection to the 
question of whether $\{u\in H^s(\Rn):\supp{u}\subset \LLambda\}=\{0\}$ for a given closed set $\LLambda\subset\Rn$ and $s\in\R$. They also both arise naturally in the study of boundary integral equation formulations of acoustic wave scattering by fractal screens. 
We additionally provide analogous results in the more general setting of Besov and Triebel--Lizorkin spaces.\\[2mm]
\indent \emph{Keywords:} density, Sobolev space, Besov space, Triebel-Lizorkin space, rough set, fractal, thick domain, d-set, screen, acoustic scattering.\\[1mm]
\indent \emph{MSC:} 46E35, 28A80.\\[2mm]
\copyright 2021. Licensed under the CC BY-NC-ND 4.0 license http://creativecommons.org/licenses/by-nc-nd/4.0/. Formal publication: https://doi.org/10.1016/j.jfa.2021.109019.

\end{abstract}

\section{Introduction \label{sec:Intro}}

Consider the following two density questions for the classical Hilbert Sobolev spaces $H^s(\Rn)$:
\begin{itemize}
\item[Q1:] When does $\widetilde{H}^s(\Omega)$ equal $H^s_{\overline\Omega}$ for $\Omega\subset \Rn$ a proper domain? 

(Equivalently: when is $\cD\OO=C^\infty_0(\Omega)$ dense in $H^s_{\overline\Omega}$?)

\item[Q2:] When is $H^{s_1}_{\Gamma}$ dense in $H^{s_2}_{\Gamma}$ for $s_1>s_2$ and $\Gamma\subset \Rn$ a closed set with empty interior?
\end{itemize}
Here, following the notational conventions in \cite{McLean}, for an open set $\Omega\subset\Rn$ and a closed set $\LLambda\subset\Rn$ the spaces $\widetilde{H}^{s}(\Omega)$ and $H_{\LLambda}^{s}$
are the closed subspaces of $H^{s}(\Rn)$, $s\in\R$, defined in the following way:
\begin{align*}
\widetilde{H}^{s}(\Omega)&:= \text{ the closure of }\mathcal{D}(\Omega)\text{ in }H^{s}(\Rn);\\
H_{\LLambda}^{s}&:=\{f\in H^{s}(\Rn):\,{\rm supp}\, f\subset \LLambda\}.
\end{align*}
One can also consider the analogous questions in the much more general setting of Besov and Triebel--Lizorkin spaces, which we shall do in the main body of the paper. But to make our initial discussions as accessible as possible, we focus in this introductory section on the special case of $H^s(\Rn)$. Our particular interest in this case stems from the second two authors' recent investigations into wave scattering by fractal screens \cite{chandler2017sobolev,ScreenPaper}, in which questions Q1 and Q2 arise quite naturally. We shall say more about the connection with this motivating application in \S\ref{Motivation}.

The answer to questions Q1 and Q2 obviously depends on both the regularity parameter $s$ and the type of domain considered. One classical result relating to Q1, appearing for example in McLean's book \cite[Thm.~3.29]{McLean}, is that $\widetilde H^s(\Omega)=H^s_{\overline\Omega}$ for all $s\in\R$ whenever $\Omega$ is $C^0$, in the sense that for every point $x\in\partial\Omega$ there exists a neighbourhood $U$ of $x$ and a Cartesian coordinate system in which 
$\Omega\cap U$ coincides with the hypograph of some continuous function from $\R^{n-1}$ to $\R$.
This result was extended by Chandler-Wilde, Hewett and Moiola in \cite[Thm.~3.24]{chandler2017sobolev} to domains that are $C^0$ except at a countable set of points $P\subset\partial\Omega$, such that $P$ has a finite number of limit points in each bounded subset of $\partial\Omega$, albeit for a limited range of $s$, namely $|s|\leq 1$ for $n\geq 2$ and $|s|\leq 1/2$ for $n=1$. This includes domains formed as unions of polygons/polyhedra touching at vertices, the ``double brick'' domain, curved cusp domains, spiral domains, and Fraenkel's ``rooms and passages'' domain --- for illustrations see \cite[Fig.~4]{chandler2017sobolev}.

Another general result one can state is that if $\widetilde{H}^s(\overline{\Omega}^\circ)=H^s_{\overline\Omega}$ then $\widetilde{H}^s(\Omega)=H^s_{\overline\Omega}$ if and only if 
$H^{-s}_{\LLambda}=\{0\}$ for every closed $\LLambda\subset \overline{\Omega}^\circ\setminus\Omega$ \cite[Lem.~3.17(v)]{chandler2017sobolev}. This result extends previous work of Maz'ya \cite[Thm.~13.2.1]{Maz'ya} and Triebel \cite{Tri:08}, which concerned the case where $\overline{\Omega}=\Rn$. It also facilitates the construction of counterexamples for which the answer to Q1 is negative \cite[\S3.5]{chandler2017sobolev}. 
Indeed, let $\Omega_0$ be a proper domain such that $\widetilde H^s(\Omega_0)= H^s_{\overline\Omega_0}$ and $\Gamma \subset\Omega_0$ a compact set with empty interior such that $H^{-s}_\Gamma \neq\{0\}$. Then $\Omega=\Omega_0\setminus \Gamma $ satisfies $\widetilde H^s(\Omega)\neq H^s_{\overline\Omega}$. As a somewhat extreme example, one can take $\Gamma $ to be the ``Swiss cheese'' set defined by Polking in \cite{Po:72a}, for which $H^{s}_\Gamma \neq\{0\}$ for all $s\leq n/2$, and $\Omega_0$ to be any open ball containing $\Gamma $. Then $\Omega=\Omega_0\setminus \Gamma $ satisfies $\widetilde H^s(\Omega)\neq H^s_{\overline\Omega}$ for all $s\geq -n/2$.
See also Lemma \ref{lem:TildeSubImpliesNull} below for a related result.

The main contribution of the current paper to the study of Q1, presented in \S\S\ref{sec:TildeSubscript}--\ref{sec:ExamplesThick}, is to extend the classical $C^0$ result in a different direction, to the case of \textit{thick domains}, in the sense of Triebel \cite[\S3]{Tri08} --- see Definition~\ref{def:thicks} below. 
One of our main results is Corollary \ref{cor:thick}, which implies that
$$
\text{if }\Omega\text{ is thick and }|\partial\Omega|=0
\;\text{ then }\;\widetilde{H}^s(\Omega)=H^s_{\overline\Omega}\;\text{ for all }s\in\R.
$$ 
This includes in particular the classical Koch snowflake domain and some of its generalisations (see \S\ref{sec:ExamplesThick}), which fail to be $C^0$ at any of its boundary points. 
Our proof uses duality arguments and the identification of $H^s_{\overline\Omega}$ with a certain space of distributions on $\Omega$ (see Lemma \ref{lem:restriction}), for which a wavelet decomposition is available (see Theorem \ref{thm:TAWavelet}, which follows from \cite[Thm.~3.13]{Tri08}).

Regarding Q2, in certain special cases it is possible to give a complete answer using known results. For instance, if $\Gamma=\R^{n-1}\times\{0\}$ 
is a $(n-1)$-dimensional hyperplane then the standard decomposition $H^s_\Gamma=\sum_{0\leq j<-s-1/2}H^{s+j+1/2}(\R^{n-1})\otimes \delta^{(j)}$, $\delta^{(j)}$ being the $j$th derivative of the one-dimensional delta function in the variable perpendicular to $\Gamma$, (see e.g.\ \cite[Lem.~3.39]{McLean}) implies that $H^{s_1}_\Gamma$ is non-trivial and dense in $H^{s_2}_\Gamma$ if and only if $-3/2 - m \leq s_2 \le s_1 < -1/2 - m$ for some $m\in\N_0$
(for the ``only if'' part, see a detailed proof of a related result in
Remark~\ref{rem:CounterexDensity}). 
By standard arguments involving coordinate charts, analogous results hold for
smooth $(n-1)$-dimensional submanifolds of $\Rn$. 
On the other hand, the only existing results we know of applicable to completely general $\Gamma$ are negative, coming from the fact that \cite[Prop.~2.4]{HewMoi:15} for every closed $\Gamma\subset\R^n$ with empty interior there exists $s_\Gamma\in[-n/2,n/2]$ (termed the ``nullity threshold'' in \cite{HewMoi:15}) such that $H^s_\Gamma\neq\{0\}$ for $s<s_\Gamma$ and $H^s_\Gamma=\{0\}$ for $s>s_\Gamma$. Hence if $s<s_\Gamma<t$ then $H^t_\Gamma=\{0\}$ cannot be dense in $H^s_\Gamma\neq \{0\}$. 

Our main contribution in this paper to the study of Q2, presented in \S\ref{sec:Density}, is to generalise, except for the limit case $s_2=-3/2-m$,
the ``if'' part of the hyperplane result mentioned above to the case where $\Gamma$ is a $d$-set (intuitively, a closed set with the same Hausdorff dimension $d$ in a neighbourhood of each of its points, see Definition \ref{def:dSet} below) for some $0<d<n$.
In particular, Theorem \ref{thm:maindensity2} implies that 
\begin{gather*}
\label{}
\text{if $\Gamma\subset\R^n$ is a $d$-set for some $0\!<\!d\!<\!n$}\\
\text{ and }-\frac{n-d}2-m-1<s_2 \leq  s_1<-\frac{n-d}2-m, \text{ for some }m\in\N_0,\\
\text{ then } H^{s_1}_\Gamma \text{ is dense in }H^{s_2}_\Gamma.
\end{gather*}
A key tool used to prove this fact is Proposition~\ref{prop:netrusov}, a consequence of a result due to Netrusov, which states in particular that if $\frac{n-d}2+m<s<\frac{n-d}2+m+1$, $m\in\N_0$, then $\tH^s(\Gamma^c)$ is the kernel of a trace operator $\Tr_{\Gamma,m}$ (defined on $H^s(\R^n)$) involving partial derivatives of order at most $m$.
Theorem~\ref{thm:gammasubgamma} shows that, under the same conditions on $s$ and $m$, the adjoint of the trace operator $\Tr_{\Gamma,m}$ provides a natural identification of the space of distributions $H^{-s}_\Gamma$ defined on $\R^n$ and supported in the $d$-set $\Gamma$ with the dual of the trace space $\Tr_{\Gamma,m}(H^s(\R^n))$. 
We also provide counterexamples showing that the assumptions made on the indices (e.g.\ on $s_1$ and $s_2$ above) are close to optimal; see Proposition~\ref{prop:d0} and Remarks~\ref{rem:CounterexDensity}--\ref{rem:CounterexDensityRough}.

As already mentioned, the results in the following sections will be presented in the wider generality of quasi-Banach Besov and Triebel--Lizorkin spaces.

\section{Motivation: scattering by fractal screens\label{Motivation}}
As mentioned above, our study is motivated by recent work by two of the authors into boundary integral equation (BIE) formulations of wave scattering by fractal screens \cite{chandler2017sobolev, ScreenPaper}, where the questions Q1 and Q2 arise naturally in the study of well-posedness and BIE solution regularity. To give context to the current study we now briefly explain this connection. 

Consider the problem of time-harmonic acoustic scattering (governed by the Helmholtz equation $\Delta u + k^2u =0$, $k>0$) in $\R^{n+1}$, $n=1,2$, by a planar screen, a bounded set $S\subset\R^n$ embedded in the hyperplane $\{x\in\R^{n+1}:\,x_{n+1}=0\}$. When $S=\Omega$ for $\Omega\subset\R^n$ a bounded open set, it was shown in \cite{ScreenPaper} that the classical Dirichlet and Neumann scattering problems (as stated in \cite{stephan87}, and see \cite[Defs.~3.10 and 3.11]{ScreenPaper}) are well-posed (and equivalent to the weak formulations in \cite{wilcox:75,NeittaanmakiRoach:87}, which view the screen as the closed set $S=\overline\Omega$) 
if and only if $\widetilde H^s(\Omega)=H^s_{\overline\Omega}$ and $H^{-s}_{\partial\Omega}=\{0\}$, with $s=-1/2$ for the Dirichlet case and $s=+1/2$ for the Neumann case. 
The unknown Cauchy data $\phi$ satisfies an associated BIE $L\phi=f$, where the data $f$ depends on the incident (source) wave field and $L$ is a bounded linear integral operator mapping bijectively between the space $\widetilde H^s(\Omega)=H^s_{\overline\Omega}$ and the space $H^{-s}(\Omega)\cong (H^{-s}_{\Omega^c})^\perp$ (orthogonal complement in $H^{-s}(\Rn)$).
One corollary of our results in the current paper is that the classical Dirichlet screen problem is well-posed whenever $\Omega\subset\R^n$ is a thick domain with $|\partial\Omega|=0$. In particular this holds for the Koch snowflake screen, for which well-posedness was raised as an open question in \cite[Examp.~8.7]{ScreenPaper}. 
On the other hand, the classical Neumann problem is \textit{not} well-posed for the Koch snowflake since $H^{-1/2}_{\partial\Omega}\neq \{0\}$ \cite[Examp.~8.7]{ScreenPaper}.

When $S=\Gamma$ for a compact set $\Gamma\subset\Rn$ with empty interior, it is also possible to formulate well-posed scattering problems, with the associated BIE posed in the space $H^s_{\Gamma}$, with data in $(\widetilde{H}^{-s}(\Gamma^c))^\perp$, again with $s=-1/2$ for the Dirichlet case and $s=+1/2$ for the Neumann case \cite{ScreenPaper,FractalBEM}.
Accordingly, the BIE solution (and hence the corresponding scattered wavefield) is non-zero (for non-zero incident data) if and only if the space $H^s_{\Gamma}$ is non-trivial. Furthermore, when $H^s_{\Gamma}$ is non-trivial and the BIE solution $\phi$ is non-zero, it is important to know whether $\phi$ possesses any extra smoothness 
(beyond membership of $H^s_{\Gamma}$) 
that can be exploited, for instance, to prove approximation error estimates for numerical discretizations. 
A natural question is whether $\phi$ lies in $H^t_\Gamma$ for some $t>s$. To our knowledge this question is almost completely open, with the only results we know of being negative, namely that if $s<s_\Gamma$ (where $s_\Gamma$ is the nullity threshold defined at the end of \S\ref{sec:Intro}) then a non-zero BIE solution $0\neq \phi\in H^s_\Gamma$ cannot lie in $H^t_\Gamma$ for any $t>s_\Gamma$ because $H^t_\Gamma=\{0\}$. A satisfactory answer to the question of solution regularity will necessitate a study of the relevant boundary integral operators, which we do not want to go into here. 
The density question Q2, however, is a weaker condition that can be investigated purely using function space theory. It represents a necessary condition for increased solution regularity, in the sense that if the BIE solution were known to lie in $H^t_\Gamma$ for all data $f$ in some dense subspace of the range of $L$ (for example, plane incident waves, see \cite{Ch:13}), then the boundedness of $L^{-1}$ would imply that $H^t_\Gamma$ is dense in $H^s_\Gamma$. Question Q2 also provides a pathway to proving convergence of numerical discretizations: if approximation error estimates can be proved for elements of $H^t_\Gamma$ for some $t>s$, and $H^t_\Gamma$ is dense in $H^s_\Gamma$, then one can prove convergence of the numerical discretization, by first approximating $\phi\in H^s_\Gamma$ by some $\tilde{\phi}\in H^t_\Gamma$ and then applying the numerical approximation theory to $\tilde{\phi}\in H^t_\Gamma$ --- for details see \cite{FractalBEM}.

\section{Preliminaries}
\label{sec:Prelim}
In this paper we are concerned with finding sufficient conditions under which the answers to Q1 and Q2 are affirmative. 
While Q1 and Q2 were posed in the context of the Sobolev spaces $H^s(\Rn)$, the approach to be used relies on results available in the more general framework of Triebel--Lizorkin spaces $\Fspq$ and Besov spaces $\Bspq$, where $s\in\R$ and $0<p,q<\infty$. Hence, whenever it does not complicate the argument we work in this more
general setting. Furthermore, we adopt the convention of using the letter $A$ instead of $F$ or $B$ in our notation when we want to mention both cases, so that statements can be read either by replacing $A$ by $F$ all over or by replacing $A$ by $B$ all over. 
With this convention we define the spaces
\begin{align}
\label{eqn:TildeDef}
\tildeA &:= \text{ the closure of }\mathcal{D}(\Omega)\text{ in }\Aspq;\\
\label{eqn:SubscriptDef}
A_{p,q,\LLambda}^{s}&:=\{f\in \Aspq:\,{\rm supp}\, f\subset \LLambda\}.
\end{align}
for open $\Omega\subset\Rn$ and closed $\LLambda\subset\Rn$. 
We note that $\subA$ is denoted $\widetilde{A}^s_{p,q}(\overline{\Omega})$ by Triebel in \cite[Def.~2.1(ii)]{Tri08}; our notation is an extension of that used in \cite{HewMoi:15,chandler2017sobolev,ScreenPaper}.

As for the definition of the Triebel--Lizorkin and Besov spaces themselves, they are quite standard and can be found in several reference works of
Triebel, e.g.\ in \cite[\S2.3.1]{Tri83} or in the more recent book \cite[Def.~1.1]{Tri08}
which we are going to refer to extensively. In \cite[p.~37]{Tri83}
the reader can also recall the definition of the Bessel-potential
spaces $H_{p}^{s}(\Rn)$, $s\in\R$, $1<p<\infty$, and both in \cite[\S2.3.5]{Tri83}
and in \cite[Rmk.~1.2]{Tri08} one can find the relation
\begin{equation}
H_{p}^{s}(\Rn)=F_{p,2}^{s}(\mathbb{R}^{n}),\quad s\in\R,\quad1<p<\infty,\label{eq:H=00003DF}
\end{equation}
between the Bessel-potential Sobolev spaces and the Triebel--Lizorkin spaces. The reader who is not familiar with such spaces might
also want to consult \cite[\S2.3.2, \S2.3.3]{Tri83}, where some of their basic properties are presented, which we may use without further warning. We note that the spaces $H^{s}(\Rn)$ considered above are, by definition, the same as $H_{2}^{s}(\Rn)$. 
We emphasize that the equality relation in (\ref{eq:H=00003DF})
indicates equality as sets but in general only equivalence of norms. In other words, it says that the identity operator is a linear and topological isomorphism between the two spaces.

We will make frequent use of the following standard duality result\footnote{In this paper, dual spaces consist of bounded \textit{linear} functionals. In previous work by the second two authors (e.g.\ \cite{HewMoi:15,chandler2017sobolev}), they are assumed to consist of bounded \emph{anti-linear} functionals, for reasons of notational convenience. Complex conjugation provides an isometric anti-linear bijection between the two types of dual space: if $l$ is a bounded linear (resp.\ antilinear) functional then $\overline{l}$ defined by $\overline{l}(u):=\overline{l(u)}$ is a bounded antilinear (resp.\ linear) functional with the same norm as $l$.}.
Here and henceforth the numbers $p'$ and $q'$ stand for the conjugate exponents of $p$ and $q$ respectively. 
We denote by $\cS(\Rn)$ the Schwartz space and by $\cS'(\Rn)$ its dual, the space of tempered distributions.
\begin{prop}[{\!\!\cite[Thm.~2.11.2]{Tri83}}]
\label{prop:dualA}Given any $s\in\R$ and $1<p,q<\infty$, the operator
\[
I_{p',q'}^{-s,A}:\, A_{p',q'}^{-s}(\Rn)\longrightarrow\big(\Aspq\big)'
\]
defined by 
\[
(I_{p',q'}^{-s,A}f)(g)=\lim_{k\to\infty}\left\langle f,g_{k}\right\rangle, \qquad \forall f\in A_{p',q'}^{-s}(\Rn),\,\forall g\in\Aspq,\;
\]
where $\left\langle \cdot,\cdot\right\rangle $ is the dual pairing on $\mathcal{S}'(\Rn)\times \mathcal{S}(\Rn)$ and $(g_{k})_{k\in\N}\subset\mathcal{S}(\Rn)$ is any sequence converging to $g$ in $\Aspq$,  
is a linear and topological isomorphism.
\end{prop}
\begin{rem}
\label{rem:SPrime}
It follows from the proof of \cite[Thm.~2.11.2]{Tri83}
that there exists $c>0$ such that for each $f\in A_{p',q'}^{-s}(\Rn)$ and $\varphi\in\mathcal{S}(\Rn)$,
\begin{equation}
\left|\left\langle f,\varphi\right\rangle \right|\leq c\,\|f|A_{p',q'}^{-s}(\Rn)\|\|\varphi|\Aspq\|.\label{eq:duality}
\end{equation}
This, together with the density of the embedding $\mathcal{S}(\Rn)\hookrightarrow\Aspq$, guarantees that the construction of $I_{p',q'}^{-s,A}f$ above makes sense (in particular, does not depend on the choice of approximating sequence $(g_{k})_{k\in\N}$), and defines an element of $(\Aspq)'$.
That $I_{p',q'}^{-s,A}$ is a linear and topological isomorphism is then precisely the content of \cite[Thm.~2.11.2]{Tri83}.\end{rem}
\begin{cor}
\label{cor:comute}Given any $s\in\R$ and $1<p,q<\infty$, any \textup{$f\in A_{p',q'}^{-s}(\Rn)$}
and any\textup{ $g\in\Aspq$}, 
\begin{equation}
(I_{p',q'}^{-s,A}f)(g)=(I_{p,q}^{s,A}g)(f).\label{eq:comute}
\end{equation}
\end{cor}
\begin{proof}
Consider $(f_{k})_{k\in\N}\subset\mathcal{S}(\Rn)$ converging to
$f$ in $A_{p',q'}^{-s}(\Rn)$ and $(g_{k})_{k\in\N}\subset\mathcal{S}(\Rn)$
converging to $g$ in $\Aspq$, and write
\begin{align*}
&{|(I_{p',q'}^{-s,A}f)(g)-(I_{p,q}^{s,A}g)(f)|}\\
 & \leq  |(I_{p',q'}^{-s,A}f)(g)-\left\langle f,g_{k}\right\rangle |+|\left\langle f,g_{k}\right\rangle -\left\langle f_{k},g_{k}\right\rangle |+|\left\langle g_{k},f_{k}\right\rangle -\left\langle g,f_{k}\right\rangle |+|\left\langle g,f_{k}\right\rangle -(I_{p,q}^{s,A}g)(f)|.
\end{align*}
The first and last terms on the right-hand side clearly tend
to $0$ when $k$ goes to $\infty$, by definition of the operators
$I_{p',q'}^{-s,A}$ and $I_{p,q}^{s,A}$. That the same happens to the
middle terms follows from (\ref{eq:duality}) and the hypotheses considered
here. Of course, we are using the facts $-(-s)=s$,
$(p')'=p$ and $(q')'=q$.\end{proof}
\begin{rem}
\label{rem:extension}
The operator $I_{p',q'}^{-s,A}$ is by construction an extension of the dual pairing $\langle \cdot, \cdot \rangle$, in the sense that if $f\in A_{p',q'}^{-s}(\Rn)$ and $g\in\mathcal{S}(\Rn)$ then $(I_{p',q'}^{-s,A}f)(g)=\left\langle f,g\right\rangle$. 
Therefore it is common to continue writing $\left\langle f,g\right\rangle $
instead of $(I_{p',q'}^{-s,A}f)(g)$ even when $g\not\in\mathcal{S}(\Rn)$.
In particular, with this convention the identity (\ref{eq:comute}) can be written as $\left\langle f,g\right\rangle =\left\langle g,f\right\rangle $.
\end{rem}

The following proposition provides an important connection between the ``tilde'' and ``subscript'' spaces introduced in \rf{eqn:TildeDef} and \rf{eqn:SubscriptDef}. 
Here, and henceforth, the superscript ``$a$'' stands for annihilator. 
We note that this result was proved for the special case of $H^s(\Rn)$ in \cite[Lem.~3.2]{chandler2017sobolev}.
\begin{prop}
\label{prop:annihilators}Given a closed set $\Gamma\subset\Rn$, $s\in\R$ and $1<p,q<\infty$,
\[
\big(\tildeAg\big)^{a}=I_{p',q'}^{-s,A}(\subAgp)\qquad\mbox{and}\qquad(\subAg)^{a}=I_{p',q'}^{-s,A}\big(\tildeAgp\big).
\]
\end{prop}
\begin{proof}
For the first identity, by the continuity of $\langle\cdot, \cdot\rangle$ and the density of $\cD(\Gamma^c)$ in $\widetilde A^s_{p,q}(\Gamma^c)$ we have
\begin{align*}
\nonumber
\big(I_{p',q'}^{-s,A}\big)^{-1}\Big(\big(\widetilde A^s_{p,q}(\Gamma^c)\big)^{a}\Big)
&=\Big\{\psi\in A^{-s}_{p',q'}(\R^n), \; \langle \psi,u\rangle =0 \; \forall u\in \widetilde A^s_{p,q}(\Gamma^c)\Big\}\\
&=\Big\{\psi\in A^{-s}_{p',q'}(\R^n), \; \langle \psi,u\rangle =0 \; \forall u\in \cD(\Gamma^c)\Big\}
= A^{-s}_{p',q',\Gamma}.
\end{align*}
For the second identity, replacing $(s,p,q)$ by $(-s,p',q')$ in the first identity gives $(\tildeAgp)^{a}=I_{p,q}^{s,A}(\subAg)$, and then applying a left annihilator to both sides we have by \cite[Prop.~1.10.15(c)]{megginson1998introduction} (noting that $\tildeAgp$ is closed) and Corollary \ref{cor:comute} that 
\begin{align*}
\label{}
\tildeAgp={}^a\big((\tildeAgp)^{a}\big)&={}^a\big(I_{p,q}^{s,A}(\subAg)\big)
= \Big\{ \psi\in A^{-s}_{p',q'}(\R^n), \; \langle u,\psi\rangle =0 \; \forall u\in \subAg \Big\}\\
&= \Big\{ \psi\in A^{-s}_{p',q'}(\R^n), \; \langle \psi,u \rangle =0 \; \forall u\in \subAg \Big\}
= \big(I_{p',q'}^{-s,A}\big)^{-1}\big((\subAg)^{a}\big).
\end{align*}
\end{proof}

A key concept arising in our study of both Q1 and Q2 is that of ``$A^{s}_{p,q}$-nullity''. 
\begin{defn}
A closed set $\LLambda\subset\Rn$ is said to be ``$A^{s}_{p,q}$\emph{-null}'' if
$A^{s}_{p,q,\LLambda}=\{0\}$.
\end{defn}
Conditions for $H^s_p$-nullity were studied in detail in \cite{HewMoi:15} using classical potential theoretic results on capacities \cite{AdHe}. Combining the results in \cite[Thm.~2.12]{HewMoi:15} with the standard embeddings in 
\cite[Prop.~2.3.2.2, Thm.~2.7.1]{Tri83} and some knowledge about delta functions leads to the following general statements concerning sets $\LLambda$ with zero Lebesgue measure.

\begin{prop}\label{prop:Nullity}
Let $\LLambda\subset\R^n$ be non-empty and closed with $|\LLambda|=0$, and define $d:=\dimH(\LLambda)\in[0,n]$, where $\dimH(\cdot)$ denotes the Hausdorff dimension. 
Then, for any $0<q<\infty$,
\begin{enumerate}[(i)]
\item $\LLambda$ is $A^s_{p,q}$-null (i.e.\ $A^s_{p,q,\LLambda} = \{0\}$) if either 
\[1<p<\infty \text{ and }s>
\dfrac{d-n}{p'} = (n-d)\left(\frac1p-1\right) \qquad \text{or} \qquad
0<p\leq 1 \text{ and }s>n\left(\frac1p-1\right)\geq 0;
\]
\item $\LLambda$ is not $A^s_{p,q}$-null (i.e.\ $A^s_{p,q,\LLambda} \neq \{0\}$) if either 
\[1<p<\infty \text{ and }s<\dfrac{d-n}{p'}= (n-d)\left(\frac1p-1\right)\le0 \qquad \text{or} \qquad
0<p\leq 1 \text{ and } s<n\left(\frac1p-1\right).
\]
\end{enumerate}
\end{prop}
\begin{proof}

(i) Let $1<p<\infty$ and $s>(d-n)/p'$, and choose $\tilde{s}$ satisfying $s>\tilde{s}>(d-n)/p'$. Then by \cite[Thm.~2.12]{HewMoi:15} we have that $F^{\tilde{s}}_{p,2,\LLambda}=
H^{\tilde{s}}_{p,\LLambda}=\{0\}$, and by \cite[Prop.~2.3.2.2(iii)]{Tri83} it follows that also $B^{\tilde{s}}_{p,\min\{p,2\},\LLambda}=\{0\}$.
Since $s>\tilde{s}$ we can use \cite[Prop.~2.3.2.2(ii)]{Tri83} to deduce that $B^{s}_{p,q,\LLambda}=\{0\}$ and $F^{s}_{p,q,\LLambda}=\{0\}$ for all $0<q<\infty$, as claimed. 

Now let $0<p\leq 1$ and $s>n(1/p-1)$, and choose $1<\tilde{p}<\infty$ satisfying
\[s>\frac{d}{\tilde{p}'} + n \left(\frac{1}{p}-1\right)\geq n \left(\frac{1}{p}-1\right). \]
(If $d=0$ this is trivially true for all $1<\tilde{p}<\infty$.) Then $\tilde{s}:=s + n\left(\frac{1}{\tilde{p}}-\frac{1}{p}\right)$ satisfies 
\[ s>\tilde{s}> \frac{d}{\tilde{p}'}+ n \left(\frac{1}{p}-1\right) + n\left(\frac{1}{\tilde{p}}-\frac{1}{p}\right) = \frac{d-n}{\tilde{p}'},
\]
and, arguing as above, we have that
$B^{\tilde{s}}_{\tilde{p},q,\LLambda}=\{0\}$ and $F^{\tilde{s}}_{\tilde{p},q,\LLambda}=\{0\}$ for $0<q<\infty$. Furthermore, using \cite[Thm.~2.7.1]{Tri83} we deduce that
$B^{s}_{p,q,\LLambda}=\{0\}$ and $F^{s}_{p,q,\LLambda}=\{0\}$ for $0<q<\infty$, as claimed. 

(ii) Let $1<p<\infty$ and $s<(d-n)/p'$, and choose $\tilde{s}$ satisfying $s<\tilde{s}<(d-n)/p'$. Then by \cite[Thm.~2.12]{HewMoi:15} we have that $F^{\tilde{s}}_{p,2,\LLambda}=
H^{\tilde{s}}_{p,\LLambda}\neq \{0\}$, and by \cite[Prop.~2.3.2.2(iii)]{Tri83} it follows that $B^{\tilde{s}}_{p,\max\{p,2\},\LLambda}\neq \{0\}$.
Since $s<\tilde{s}$ we can use \cite[Prop.~2.3.2.2(ii)]{Tri83} to deduce that $B^{s}_{p,q,\LLambda}\neq\{0\}$ and $F^{s}_{p,q,\LLambda}\neq\{0\}$ for all $0<q<\infty$, as claimed. 

Now let $0<p\leq 1$ and $s<n(1/p-1)$. Since there is at least one point $a\in\Gamma$, from \cite[Rmk.~2.2.4.3]{RuS96} and \cite[Prop.~2.3.8]{Tri83} it follows that $\delta_a \in B^{s}_{p,q,\{a\}} \subset B^{s}_{p,q,\LLambda}$ for any $0<q<\infty$. Using \cite[Prop.~2.3.2.2(iii)]{Tri83} we get that $A^{s}_{p,q,\LLambda} \neq \{0\}$ for $0<q<\infty$, as claimed.
\end{proof}

\begin{rem}
\label{rem:nullity}
The excluded cases ($1<p<\infty$ with $s=(d-n)/p'$ and $0< p\leq 1$ with $s=n(1/p-1)$) are delicate and are not discussed here. 
For $1<p<\infty$ and $H^{(d-n)/p'}_{p}$-nullity, all possible behaviours are detailed and exemplified in \cite[Cor.~2.15 and Thm.~4.5]{HewMoi:15}. In particular we note that if $0<d<n$ and $\LLambda$ is a compact $d$-set (see Definition \ref{def:dSet} below) or a $d$-dimensional hyperplane (with $d\in\N$) then $H^{(d-n)/p'}_{p,\LLambda}=\{0\}$ for all $1<p<\infty$ \cite[Thm.~2.17]{HewMoi:15}.
\end{rem}

The concept of a ``$d$-set'', already mentioned in the previous remark, will play an important role in our later considerations. We give a definition here.
\begin{defn}\label{def:dSet}
Let $\Gamma$ be a non-empty closed subset of $\Rn$ and $0\leq d\leq n$.
$\Gamma$ is said to be a {\em$d$-set} if there exist $c_{1},c_{2}>0$
such that
\[
c_{1}r^{d}\leq\mathcal{H}^{d}\big(B(\gamma,r)\cap\Gamma\big)\leq c_{2}r^{d},\qquad\gamma\in\Gamma,\quad0<r\leq1,
\]
where $B(\gamma,r)$ is the closed ball of radius $r$ with centre at $\gamma$ and $\mathcal{H}^{d}$ stands for the $d$-dimensional Hausdorff measure on $\Rn$.
\end{defn}

As we shall show in Propositions \ref{prop:class-snow-d-set} and \ref{prop:square-snow-d-set}, the boundaries of the snowflake domains considered in \S\ref{sec:ExamplesThick} are all examples of (compact) $d$-sets in $\R^{2}$ with $0<d<2$.
For less exotic but nonetheless important examples, given $d\in \{1,2,\ldots,n\}$, every $d$-dimensional closed Lipschitz manifold is a $d$-set in $\Rn$.
For more information about $d$-sets, see, e.g., \cite[II.1]{JW84} and \cite[I.3]{Tri97}. 
In particular (see \cite[Cor. 3.6]{Tri97}), for a $d$-set $\Gamma$ with $0<d<n$ one has that $|\Gamma|=0$ and $\dimH(\Gamma) =d$.

\section{Equality between \texorpdfstring{$\tildeA$ and $\subA$}{tilde spaces subscript spaces}}\label{sec:TildeSubscript}
Our aim in this section is to determine conditions under which
\begin{align}
\label{eqn:TildeSubscript}
\tildeA=\subA,
\end{align}
where $\Omega\subset\Rn$ is a domain (non-empty open set) and $\tildeA$ and $\subA$ are defined as in \eqref{eqn:TildeDef}--\eqref{eqn:SubscriptDef}. 
Since \eqref{eqn:TildeSubscript} holds trivially when $\Omega=\Rn$, our interest is in the case where $\Omega$ is a \emph{proper domain}, i.e.\ $\Omega\neq \Rn$.
We start by remarking that the inclusion
\[
\tildeA\subset\subA
\]
is clear, since $\mathcal{D}(\Omega)\subset\subA$ and the latter
is a closed subspace of $\Aspq$. Therefore, to prove \eqref{eqn:TildeSubscript} we shall be merely concerned with proving that $\subA\subset\tildeA$.

We deal first with the simplest case where $A=F$, $s=0$, $1<p<\infty$
and $q=2$. By (\ref{eq:H=00003DF}) this means the setting of $H_{p}^{0}(\Rn)$, $1<p<\infty$, or, to put it
simpler, $L_{p}(\Rn)$, $1<p<\infty$. Actually,
since the proof works also when $p=1$, we include this case in the following proposition. 
Here $\widetilde{L}_{p}(\Omega)$ and $L_{p,\overline{\Omega}}$ are defined in the obvious way, and
\begin{equation}\label{circL}
\mathringbig{L}_p(\Omega):=\{f\in L_p(\Rn):f=0 \text{ a.e.\ in }\Omega^c\}.
\end{equation}
\begin{prop}
\label{prop:Lp}
Let $\Omega$ be a domain in $\Rn$ and let $1\leq p<\infty$. 
Then $\widetilde{L}_{p}(\Omega)=\mathringbig{L}_p(\Omega)$. 
If $|\partial\Omega|=0$ then also $\widetilde{L}_{p}(\Omega)=L_{p,\overline{\Omega}}$.
\end{prop}
\begin{proof}
Since $\mathcal{D}(\Omega)\subset \mathringbig{L}_p(\Omega)$ and the latter is a closed subspace of $L_p(\R^n)$, then $\widetilde{L}_{p}(\Omega)\subset\mathringbig{L}_p(\Omega)$. Let now $f\in \mathringbig{L}_p(\Omega)$. Then $\|f|_{\Omega}|L_{p}(\Omega)\|=\|f|L_{p}(\Rn)\|$, which, together with the fact that $\cD\OO|_\Omega$ is dense in $L_{p}(\Omega)$ (e.g., \cite[Cor.~3.5]{McLean}), proves that also $f \in \widetilde{L}_{p}(\Omega)$. 
If $|\partial\Omega|=0$ then obviously $L_{p,\overline{\Omega}}=\mathringbig{L}_p(\Omega)$, so that also $L_{p,\overline{\Omega}}=\widetilde{L}_{p}(\Omega)$ by the first part.
\end{proof}

As mentioned above, we shall make frequent reference to some results of Triebel in \cite{Tri08}. However, there is an an unfortunate clash between the notation in \cite{Tri08} and some of the notation introduced above, which follows the conventions adopted in the second two authors' previous papers \cite{HewMoi:15,chandler2017sobolev,ScreenPaper}. 
We already pointed out immediately after \rf{eqn:SubscriptDef} that the space we call $\subA$ is denoted $\widetilde{A}^s_{p,q}(\overline{\Omega})$ in \cite[Def.~2.1(ii)]{Tri08}. In \cite[Def.~2.1(ii)]{Tri08} Triebel introduces another space that will be important for our purposes, defined in Definition \ref{def:TA} below. 
Triebel denotes this space $\widetilde{A}^s_{p,q}(\Omega)$, but since we are already using the notation $\widetilde{A}^s_{p,q}(\Omega)$ (see \rf{eqn:TildeDef}), we 
instead denote this new space $\TAO$, with the ``R'' highlighting the fact that $\TAO$ is a space of restrictions of distributions in $\subA$.
\begin{defn}[{\!\!\cite[Def.~2.1(ii)]{Tri08}}]
\label{def:TA}
Let $\Omega$ be a domain in $\Rn$.
Let $s\in\R$ and $0<p,q<\infty$.
\[
\TAO:=\{f\in\mathcal{D}'(\Omega):\, f=g|_{\Omega}\mbox{ for some }g\in \subA\};
\]
\[
\|f|\TAO\|:=\inf\,\|g|\Aspq\|,
\]
where the infimum is taken over all $g\in \subA$ with $g|_{\Omega}=f$. 
\end{defn}
\begin{rem}
\label{rem:thesame}
The norm on $\TAO$ defined above is in general stronger than that inherited from the usual restriction space $A_{p,q}^{s}(\Omega):=\{f\in\mathcal{D}'(\Omega):\, f=g|_{\Omega}\mbox{ for some }g\in A_{p,q}^{s}(\Rn)\}$, where the norm involves an infimum over all $g\in A_{p,q}^{s}(\Rn)$ such that $g|_{\Omega}=f$.
\end{rem}

It is mentioned in \cite[Rmk.~2.2]{Tri08} that there is a one-to-one correspondence between $\subA$ and $\TAO$ if, and only if, 
\begin{equation}
A_{p,q,\partial\Omega}^{s}=\{0\},
\label{eq:trivialonboundary}
\end{equation}
i.e.\ $\partial\Omega$ is $A^s_{p,q}$-null.
Indeed, using standard arguments from the theory of distributions
we can be more precise and state the following:
\begin{lem}
\label{lem:restriction}Let $\Omega$ be a domain in $\Rn$.
Let $s\in\R$ and $0<p,q<\infty$. If $\partial\Omega$ is $A^s_{p,q}$-null (i.e., (\ref{eq:trivialonboundary}) holds), then the restriction operator
\[
|_\Omega:\, \subA\longrightarrow \TAO
\]
is an isometric isomorphism; in particular,
\begin{equation}
\|f|_\Omega|\TAO\|=\|f|\Aspq\|\quad\mbox{ for all }f\in \subA.\label{eq:restriction}
\end{equation}

\end{lem}
The importance of this result is the following: there are some
results for $\TAO$ in \cite{Tri08} that
we would like to transfer to $\subA$; this is possible whenever the above lemma applies. In particular, by Proposition \ref{prop:Nullity} this holds for $\Omega\neq\Rn$ whenever $|\partial\Omega|=0$ and either 
\[1<p<\infty \text{ and }s>\dfrac{\dimH{\partial\Omega}-n}{p'} \qquad \text{or} \qquad
0<p\leq 1 \text{ and }s>n(1/p-1)\geq 0.\]

In order to state the main results of this section later on, we shall need the following notions from \cite[Def.~3.1(ii)--(iv), Rmk.~3.2]{Tri08}. 
Here, and henceforth, for a set $S\subset\R^n$, we denote by $\mathcal Q(S)$ the set of the (open) cubes contained in $S$ and with the edges parallel to the Cartesian axes; and for any $Q\in\mathcal Q(S)$ we denote by $l(Q)$ the length of its edges.
\begin{defn}
\label{def:thicks}Let $\Omega$ be a proper domain in $\Rn$.\begin{enumerate}[(i)]
\item $\Omega$ is said to be {\em$E$-thick (exterior thick)} if for any choice
of $c_{1},c_{2},c_{3},c_{4}>0$ and $j_{0}\in\N$ there are $c_{5},c_{6},c_{7},c_{8}>0$
such that for any $j\in\N$, $j\geq j_{0}$, and any \emph{interior
cube} $Q^{i}\in \mathcal{Q}(\Omega)$
with
\[
c_{1}2^{-j}\leq l(Q^{i})\leq c_{2}2^{-j}\quad\mbox{ and }\quad{\rm c_{3}2^{-j}\leq dist}(Q^{i},\partial\Omega)\leq c_{4}2^{-j},
\]
there exists an \emph{exterior cube} $Q^{e}\in\mathcal{Q}(\Omega^{c})$ with
\[
c_{5}2^{-j}\leq l(Q^{e})\leq c_{6}2^{-j}\quad\mbox{ and }\quad c_{7}2^{-j}\leq{\rm dist}(Q^{e},\partial\Omega)\leq{\rm dist}(Q^{i},Q^{e})\leq c_{8}2^{-j}.
\]

\item $\Omega$ is said to be {\em$I$-thick (interior thick)} if for any choice
of $c_{1},c_{2},c_{3},c_{4}>0$ and $j_{0}\in\N$ there are $c_{5},c_{6},c_{7},c_{8}>0$
such that for any $j\in\N$, $j\geq j_{0}$, and any exterior cube
$Q^{e}\in\mathcal{Q}(\Omega^{c})$ with
\[
c_{1}2^{-j}\leq l(Q^{e})\leq c_{2}2^{-j}\quad\mbox{ and }\quad c_{3}2^{-j}\leq{\rm dist}(Q^{e},\partial\Omega)\leq c_{4}2^{-j},
\]
there exists an interior cube $Q^{i}\in\mathcal{Q}(\Omega)$ with
\[
c_{5}2^{-j}\leq l(Q^{i})\leq c_{6}2^{-j}\quad\mbox{ and }\quad c_{7}2^{-j}\leq{\rm dist}(Q^{i},\partial\Omega)\leq{\rm dist}(Q^{e},Q^{i})\leq c_{8}2^{-j}.
\]
\item $\Omega$ is said to be {\em thick} if it is both $E$-thick and $I$-thick.
\end{enumerate}
\end{defn}

\begin{rem}
\label{rem:ThickTest}
It is easily seen that: 
\begin{enumerate}
\item Once the definition of $E$-thickness (or $I$-thickness) has been checked for some $j_0\in \N$ then it automatically holds for all $j_0$;
\item The definitions of $E$-thickness and $I$-thickness can be equivalently stated with $2^{-j}$ replaced throughout by $\xi^j$ for any $0<\xi<1$.
\end{enumerate}
\end{rem}

In \cite[Prop.~3.6(i)--(iv), Prop.~3.8(i),(iii)]{Tri08}, some relations with well-known concepts are presented:
\begin{prop}
\label{prop:snowflake}$\phantom{a}$
\begin{enumerate}[(i)]
\item Any $(\varepsilon,\delta)$-domain \cite{jones1981quasiconformal} $\Omega$
in $\Rn$ is $I$-thick with $|\partial\Omega|=0$.
\item Any bounded Lipschitz domain \cite[Def.~3.4(iii)]{Tri08} in $\Rn$,
$n\geq2$, is thick.
\item The classical Koch snowflake domain as per \cite[Fig.~3.5]{Tri08} in $\R^{2}$
is a thick $(\varepsilon,\delta)$-domain.
\item Let $\Omega$ be a domain in $\Rn$. Then
$\Rn=\Omega\cup\partial\Omega\cup\overline{\Omega}^{c}$ and $\partial(\overline{\Omega}^{c})\subset\partial\Omega$.
Furthermore,
$\partial\Omega=\partial(\overline{\Omega}^{c})$ if, and only if, $(\overline{\Omega})^{\circ}=\Omega$.
\item If $\Omega$ is an $E$-thick domain in $\Rn$, then $(\overline{\Omega})^{\circ}=\Omega$
and $\overline{\Omega}^{c}$ is $I$-thick.
\item If $\Omega$ is an $I$-thick domain in $\Rn$ and $\overline{\Omega}\not=\Rn$,
then $\overline{\Omega}^{c}$ is $E$-thick.
\end{enumerate}
\end{prop}
We shall also need to use wavelet representations of some spaces,
needing in particular to consider so-called orthonormal $u$-wavelet
basis in $L_{2}(\Omega)$. However, we don't want to go into details,
so we shall keep things at the bare minimum.
\begin{defn}[{\!\!\cite[Def.~2.31]{Tri08}}]
 Let $\Omega$ be a proper domain in $\Rn$.
Let $u\in\N$. A collection (of real functions)
\[
\{\Phi_{r}^{j}:\: j\in\No,\: r=1,\ldots,N_{j}\}\quad\mbox{ with }N_{j}\in\N\cup\{\infty\}
\]
is called an \emph{orthonormal $u$-wavelet basis} in $L_{2}(\Omega)$
if it is both an $u$-wavelet system according to \cite[Def.~2.4]{Tri08}
and an orthonormal basis in $L_{2}(\Omega)$.\end{defn}
\begin{rem}
\label{rem:wavelets} We do not go deeper into the long definition
of what an \emph{$u$-wavelet system} is because, in addition to what
we are going to write down below, we shall only need the following
two properties: for all $j$ and $r$ as above,
\begin{enumerate}[(i)]
\item $\Phi_{r}^{j}$ belongs to $C^{u}(\Rn)$;
\item ${\rm supp}\,\Phi_{r}^{j}\subset\Omega$ .
\end{enumerate}
From these two properties it follows that $\Phi_{r}^{j}\in\tildeA$
when $u>s$: that $\Phi_{r}^{j}\in\Aspq$ follows from the fact that
then $\Phi_{r}^{j}$ is essentially an atom in $\Aspq$ --- see, e.g.,
\cite[Cor.~4.11]{AC16b}, read in the constant exponents case; that
it can be approximated in $\Aspq$ by functions in $\mathcal{D}(\Omega)$
follows from the density of $\mathcal{S}(\Rn)$ in $\Aspq$ with the
help of a suitable cut-off function and pointwise multiplier properties
(see, e.g., \cite[\S2.8.2]{Tri83}), since the second property
above guarantees that there is some room between ${\rm supp}\,\Phi_{r}^{j}$
and $\Omega^{c}$.\end{rem}
\begin{thm}[{\!\!\cite[Thm.~2.33]{Tri08}}]
Let $\Omega$ be a proper domain in $\Rn$.
For any $u\in\N$ there exist orthonormal $u$-wavelet bases in $L_{2}(\Omega)$.
\end{thm}
The next result, which will be crucial for our intentions, follows
from \cite[Def.~3.11, Thm.~3.13]{Tri08}:
\begin{thm}
\label{thm:TAWavelet}
Let $\Omega$ be an $E$-thick domain in $\Rn$. Let $0<p,q<\infty$
and 
\[
s>n\left(\frac{1}{\min\{1,p,q\}}-1\right)\quad\mbox{ if }A=F,\qquad s>n\left(\frac{1}{\min\{1,p\}}-1\right)\quad\mbox{ if }A=B.
\]
Let $u>s$ be a natural number and 
$\{\Phi_{r}^{j}:\: j\in\No,\: r=1,\ldots,N_{j}\}$ with $N_{j}\in \N\cup\{\infty\}$
be an orthonormal $u$-wavelet basis in $L_{2}(\Omega).$ 
Then $\{\Phi_{r}^{j}\}$ is an unconditional basis in $\TAO$.
\end{thm}
We can now prove one of the main results in this section: 
\begin{thm}
\label{thm:rightequality}Let $\Omega$ be an $E$-thick domain in
$\Rn$ with $|\partial\Omega|=0$. Let $p$, $q$ and $s$ be as in
Theorem~\ref{thm:TAWavelet}. 
Then $\tildeA=\subA$.\end{thm}
\begin{proof}
The hypotheses on $p$, $q$ and $s$, and the fact that $|\partial\Omega|=0$, together imply by Proposition \ref{prop:Nullity}(i) that $\partial\Omega$ is $A^s_{p,q}$-null, so that Lemma \ref{lem:restriction} applies. 
Given any $f\in\subA$, we have
that $f|_\Omega\in \TAO$, which, by the preceding theorem, is the
limit in $\TAO$, when the natural $N$ tends to $\infty$, of
finite linear combinations $f_{N}$ of functions $\Phi_{r}^{j}$.
From Remark \ref{rem:wavelets} it follows that $\tilde f_{N}\in\tildeA\subset\subA$, where $\tilde f_{N}$ is the extension of $f_N$ to $\R^n$ by zero.
Hence from (\ref{eq:restriction}) we get that $\|f-\tilde f_{N}|\Aspq\|=\|f|_\Omega-f_{N}|\TAO\|$, which tends to
$0$ when $N$ goes to $\infty$, so we conclude that $f$ is in the closure of $\tildeA$ in $\Aspq$, that is, $f\in\tildeA$ too.
\end{proof}
\begin{cor}
\label{cor:rightequalityH}
Let $\Omega$ be an $E$-thick domain in
$\Rn$ with $|\partial\Omega|=0$. 
Then $\widetilde{H}_{p}^{s}(\Omega)=H_{p,\overline{\Omega}}^{s}$ whenever $s\geq0$ and $1<p<\infty$.\end{cor}
\begin{proof}
The case $s>0$ follows from (\ref{eq:H=00003DF}) and the above theorem.
The case $s=0$ follows from Proposition \ref{prop:Lp} (and in this
case we don't even need the domain $\Omega$ to be $E$-thick).
\end{proof}
For completeness we remark that for some parameters $p$ and $q$ and with some extra conditions on $\partial\Omega$ it is possible
to get the conclusion of Theorem \ref{thm:rightequality} for some
negative values of $s$. This follows by conjugating \cite[Def.~3.11, Thm.~3.13, Prop.~3.19 and Rmk.~3.20]{Tri08}.
And once we have the conclusion of Theorem \ref{thm:rightequality}
both for some positive and for some negative values of $s$, under
even more stringent conditions interpolation techniques can be applied
to get the conclusion for some parameters $p$ and $q$ when $s=0$
--- see \cite[Def.~3.11, Prop.~3.19, Rmk.~3.20 and Prop.~3.21]{Tri08}.
However, we shall not pursue the above avenue of research here. For the example of $H^s(\Rn)$ motivating our studies, the case $s=0$ is already covered by Corollary \ref{cor:rightequalityH}. And we shall reach negative values of $s$ under somewhat different conditions using duality.
\begin{lem}\label{cor:omegaandomegac}
Let $\Omega$ be a domain in $\Rn$. Let $s\in\R$ and $1<p,q<\infty$. Then
\[
\tildeA=\subA\qquad\mbox{if and only if}\qquad\widetilde{A}_{p',q'}^{-s}(\overline{\Omega}^{c})=A_{p',q',\Omega^{c}}^{-s}.
\]
\end{lem}
\begin{proof}
$\tildeA=\subA$ iff $(\tildeA)^{a}=(\subA)^{a}$, since $\tildeA$
and $\subA$ are closed subspaces of $\Aspq$. On
the other hand, the latter identity is equivalent to $A_{p',q',\Omega^{c}}^{-s}=\widetilde{A}_{p',q'}^{-s}(\overline{\Omega}^{c})$,
as follows by applying Proposition \ref{prop:annihilators} with $\Gamma=\Omega^{c}$
in its first identity and $\Gamma=\overline{\Omega}$ in its second
identity. 
\end{proof}
If we specialize the above result to the case $A=F$ and
$p=q=2$ we recover \cite[Lem.~3.26]{chandler2017sobolev}.

One immediate corollary of Lemma \ref{cor:omegaandomegac} is the second part of the following lemma, which provides another connection between density results and $A^s_{p,q}$-nullity. 
\begin{lem}
\label{lem:TildeSubImpliesNull}
Let $\Omega\subset\R^n$ be a proper domain, and let $s>0$ and $1<p<\infty$. 
\begin{enumerate}[(i)]
\item If $\tildeA=\subA$ and $0<q<\infty$ then $A^s_{p,q,\partial\Omega}=\{0\}$.
\item If $\widetilde{A}^{-s}_{p',q'}(\Omega)=A^{-s}_{p',q',\overline\Omega}$ and $1<q<\infty$ then  $A^s_{p,q,\partial\Omega}=\{0\}$.
\end{enumerate}
\end{lem}
\begin{proof}
(i) First note that the assumptions on $s$, $p$ and $q$ imply that $\Aspq\hookrightarrow L_p(\Rn)$. 
Suppose that $\tildeA=\subA$. Then by Proposition \ref{prop:Lp} (recall \eqref{circL} for the definition of $\mathringbig{L}_p(\Omega)$) we have
\[A^s_{p,q,\partial\Omega} = A^s_{p,q,\partial\Omega}\cap A^s_{p,q,\overline\Omega} = A^s_{p,q,\partial\Omega}\cap \tildeA \hookrightarrow L_{p,\partial\Omega}\cap \widetilde{L}_p(\Omega) = L_{p,\partial\Omega}\cap \mathringbig{L}_p(\Omega)=\{0\}. \] 

(ii) Since $\overline{\overline{\Omega}^c}\subset \Omega^c$, we have the inclusions $A^{s}_{p,q,\Omega^c}\supset A^{s}_{p,q,\overline{\overline{\Omega}^c}}\supset \widetilde{A}^{s}_{p,q}(\overline\Omega^c)$. If $\widetilde{A}^{-s}_{p',q'}(\Omega)=A^{-s}_{p',q',\overline\Omega}$ and $1<q<\infty$, then Lemma \ref{cor:omegaandomegac} implies that $A^{s}_{p,q,\Omega^c}=\widetilde{A}^{s}_{p,q}(\overline\Omega^c)$, so in fact the previously mentioned inclusions are all equalities, i.e.\ $A^{s}_{p,q,\Omega^c}=A^{s}_{p,q,\overline{\overline{\Omega}^c}}= \widetilde{A}^{s}_{p,q}(\overline\Omega^c)$. But then
\[A^s_{p,q,\partial\Omega} = 
A^s_{p,q,\partial\Omega} \cap A^s_{p,q,\Omega^c}
=A^s_{p,q,\partial\Omega} \cap \widetilde{A}^{s}_{p,q}(\overline\Omega^c)
\hookrightarrow L_{p,\partial\Omega}\cap \widetilde{L}_p(\overline\Omega^c) = L_{p,\partial\Omega}\cap \mathringbig{L}_p(\overline\Omega^c)=\{0\}.
\]
\end{proof}
Note that the statement of Lemma \ref{lem:TildeSubImpliesNull} does not extend to $s<0$; counterexamples for $n=2$ include the thick domains considered in \S\ref{subsec:classsnow}, for which $\tildeA=\subA$ for all $s\in\R \setminus \{ 0 \}$, $1<p,q<\infty$; indeed, given any $s<0$, following Remark \ref{rem:Hdim} we can pick a domain of that class whose boundary has Lebesgue measure zero and Hausdorff dimension $d$ such that $s<\frac{d-2}{p'}=\frac{d-n}{p'}$, in which case $A^s_{p,q,\partial\Omega}\neq\{0\}$ by Proposition \ref{prop:Nullity}.

We now proceed to another of the main results in this section.
\begin{thm}
\label{thm:leftequality}Let $\Omega$ be an $I$-thick domain in
$\Rn$ with $(\overline{\Omega})^{\circ}=\Omega$ and $|\partial\Omega|=0$.
Given any $s<0$ and $1<p,q<\infty$, $\tildeA=\subA$.\end{thm}
\begin{proof}
According to Lemma \ref{cor:omegaandomegac},
it is enough to show that
\begin{equation}
\widetilde{A}_{p',q'}^{-s}(\overline{\Omega}^{c})=A_{p',q',\Omega^{c}}^{-s}.\label{eq:almostlefteq}
\end{equation}
Since $\Omega\neq\Rn$, the assumption $(\overline{\Omega})^{\circ}=\Omega$ implies that also $\overline{\Omega}\not=\Rn$. Then since $\Omega$ is $I$-thick, from 
Proposition \ref{prop:snowflake}(vi) it follows that $\overline{\Omega}^{c}$
is $E$-thick. 
Furthermore, since $|\partial\Omega|=0$ it follows by Proposition \ref{prop:snowflake}(iv) that $|\partial(\overline{\Omega}^{c})|=0$. 
We can then apply Theorem \ref{thm:rightequality} with
$(-s,p',q',\overline{\Omega}^{c})$ in place of $(s,p,q,\Omega)$ to obtain $\widetilde{A}_{p',q'}^{-s}(\overline{\Omega}^{c})=A_{p',q',\overline{\overline{\Omega}^{c}}}^{-s}$, 
from which (\ref{eq:almostlefteq}) follows using  $\Omega^{c}\setminus\overline{\overline{\Omega}^{c}}=(\overline{\Omega})^{\circ}\setminus\Omega=\emptyset$,
the last identity being true by hypothesis.
\end{proof}
The next two corollaries follow immediately from Definition \ref{def:thicks},
Proposition \ref{prop:snowflake}, Theorem \ref{thm:leftequality}
and either Theorem \ref{thm:rightequality} or Corollary \ref{cor:rightequalityH}.

\mypropbox{
\begin{cor}\label{cor:equalnot0}
Let $\Omega$ be a thick domain in $\Rn$ with
$|\partial\Omega|=0$. Given any $s\in\R\setminus\{0\}$ and $1<p,q<\infty$,
$\tildeA=\subA$.
\end{cor}
}

\mypropbox{
\begin{cor}\label{cor:thick}
Let $\Omega$ be a thick domain in $\Rn$ with $|\partial\Omega|=0$.
Given any $s\in\R$ and $1<p<\infty$, $\widetilde{H}_{p}^{s}(\Omega)=H_{p,\overline{\Omega}}^{s}$.
\end{cor}
}

Due to Proposition \ref{prop:snowflake}, both corollaries above apply
to the case when $\Omega$ is the classical Koch snowflake domain in $\R^{2}$. We shall consider some further examples in the next section.

\begin{rem}
Thickness is not necessary to ensure $\tH_{p}^{s}(\Omega)=H_{p,\overline{\Omega}}^{s}$.
Indeed, for $n=2$ the domain $\Omega=\{x_2<|x_1|^{1/2}\}$ is not $E$-thick and $\baro^c$ is not $I$-thick, but $\tH^{s}(\Omega)=H_{\overline{\Omega}}^{s}$ and $\tH^{s}(\baro^c)=H_{\Omega^c}^{s}=H_{\overline{\baro^c}}^{s}$ for all $s\in\R$ by \cite[Thm.~3.29]{McLean}.
\end{rem}

\begin{rem}As suggested after Corollary \ref{cor:rightequalityH}, the conclusion
of Corollary \ref{cor:equalnot0} can also be obtained for $s=0$
if we are willing to restrict appropriately the range of domains to
be considered. But we shall not pursue this here. 
\end{rem}

\begin{rem}
Since by definition the $\tildeA$ spaces all share a common dense subspace $\cD\OO$, whenever $A^{s_1}_{p_1,q_1}(\Rn)\subset A^{s_2}_{p_2,q_2}(\Rn)$ it follows that the $\widetilde A^{s_1}_{p_1,q_1}\OO$ is dense in $\widetilde A^{s_2}_{p_2,q_2}\OO$ for all open sets $\Omega$.
Therefore, one consequence of Corollary \ref{cor:thick} is that if $\Omega$ is a thick domain in $\Rn$ with $|\partial\Omega|=0$, then $H^{s_1}_{p,\baro}$ is dense in $H^{s_2}_{p,\baro}$ for all $s_2<s_1\in\R$ and $1<p<\infty$. 
This complements the results in \S\ref{sec:Density}, where we study conditions under which $H^{s_1}_{p,\Gamma}$ is dense in $H^{s_2}_{p,\Gamma}$ for $\Gamma$ a closed set with empty interior. 
\end{rem}

\section{Examples of thick domains}
\label{sec:ExamplesThick}
In the previous section we proved that a sufficient condition for the equality $\widetilde{H}_{p}^{s}(\Omega)=H_{p,\overline{\Omega}}^{s}$ is that $\Omega$ is a thick domain with $|\partial\Omega|=0$. In this section we prove thickness for a general class of domains (possibly with fractal boundaries) formed as the limit of a sequence of smoother (``prefractal'') domains. This includes a family of generalisations of the classical Koch snowflake domain, for which a proof of thickness was sketched in \cite[Prop.~3.8(iii)]{Tri08}, and the ``square snowflake'' domain considered in \cite{sapoval1991vibrations}.

Our general result is Proposition \ref{prop:thick}. Before stating this result we need to describe the framework we have in mind. 
Suppose we have a nested increasing sequence $(\Gamma^-_j)_{j\in\N_0}$ of bounded open sets
\[\Gamma^-_0\subset\Gamma^-_1\subset \ldots \subset \Gamma^-_j \subset \Gamma^-_{j+1}\subset\ldots , \]
and a nested decreasing sequence $(\Gamma^+_j)_{j\in\N_0}$ of compact sets
\[\Gamma^+_0\supset\Gamma^+_1\supset \ldots \supset \Gamma^+_j \supset \Gamma^+_{j+1}\supset\ldots , \]
such that $\Gamma^-_j$ is non-empty for all except a finite number of $j$ and 
\[ \Gamma^-_j\subset \Gamma^+_j \qquad \text{for all }j\in\N_0. \]
Define
\[ \Gamma^-:=\bigcup_{j=0}^\infty \Gamma^-_j \qquad \text{and}\qquad
\Gamma^+:=\bigcap_{j=0}^\infty \Gamma^+_j.\]
Then $\Gamma^-$ is non-empty and open and $\Gamma^+$ is non-empty and compact, with
\[ \overline{\Gamma^-} \subset \Gamma^+. \]
Furthermore, defining the compact set
\[\Delta_j := \Gamma^+_j\setminus \Gamma^-_j, \]
it holds that
\[\partial\Gamma^- \subset \Delta_j \text{ and } \partial\Gamma^+ \subset \Delta_j, \qquad \text{for each }j\in\N_0.\]

We are now ready to state our general result concerning thickness. 
Note that conditions \eqref{Cond1}, \eqref{Cond2} and \eqref{Cond3} in Proposition \ref{prop:thick} are statements about fixed order approximations $\Gamma^\pm_j$, and do not involve the limiting objects $\Gamma^\pm$. 
One can think of \eqref{Cond2} and \eqref{Cond3} as ``$j$-uniform'' thickness estimates on the sequence of approximations $\Gamma^\pm_j$. 
We remark that a necessary condition for \rf{Cond1} to hold is that $d_H(\partial\Gamma^-_j,\partial\Gamma^+_j)\leq c\xi^{j}$, where $d_H$ is the Hausdorff distance.
\begin{prop}
\label{prop:thick}
Let $\Gamma_j^-,\Gamma_j^+,\Gamma^-,\Gamma^+,\Delta_j$, for $j\in\N_0$, be as above.
Suppose that there exists constants $0<\xi<1$, $j_*\in\N_0$ and $c, c_1^\pm,c_2^\pm,c_3^\pm,c_4^\pm>0$ such that, for each $j\ge j_*$, using the $\mathcal{Q}(S)$ and $l(Q)$ notation introduced before Definition \ref{def:thicks},
\begin{align}
\label{Cond1}
\forall x\in \Delta_j, \,\dist(x,\partial\Gamma^-_j)&\leq c\xi^{j} \text{ and } \dist(x,\partial\Gamma^+_j)\leq c\xi^{j};
\\
 \label{Cond2}
\forall x^-\in \partial\Gamma^-_j, \, \exists Q^i \in \mathcal Q(\Gamma^-_j )
\text{ s.t. } &c_1^- \xi^{j}\le l(Q^i)\le c_2^- \xi^{j} \text{ and } \nonumber\\
&c_3^- \xi^{j}\le \dist(Q^i,\partial\Gamma^-_j)\le  \dist(Q^i,x^-) \le c_4^-\xi^{j};
\\\text{and}\qquad 
\label{Cond3}
\forall x^+\in \partial\Gamma^+_j, \, \exists Q^e \in \mathcal Q\big((\Gamma^+_j)^c \big) 
\text{ s.t. } &c_1^+ \xi^{j}\le l(Q^e)\le c_2^+ \xi^{j} \text{ and }\nonumber\\
&c_3^+ \xi^{j}\le \dist(Q^e,\partial\Gamma^+_j) \le \dist(Q^e,x^+) \le c_4^+\xi^{j}.
\end{align} 
Then $\Gamma^-$ is thick, with $\overline{\Gamma^-}=\Gamma^+$ and $\Gamma^-=(\Gamma^+)^\circ$.
Moreover, if  $|\partial \Gamma^-|=0$, then $\widetilde{A}^s_{p,q}(\Gamma^-)=A^s_{p,q,\overline{\Gamma^-}}=A^s_{p,q,\Gamma^+}$ for all $s\in\R\setminus\{0\}$ and $1<p,q<\infty$, and $\widetilde{H}_{p}^{s}(\Gamma^-)=H_{p,\overline{\Gamma^-}}^{s}=H_{p,{\Gamma^+}}^{s}$ for all $s\in\R$ and $1<p<\infty$.
\end{prop}
\begin{proof}
The fact that $\overline{\Gamma^-}=\Gamma^+$ is an obvious consequence of \eqref{Cond1}, since if $x\in\Gamma^+\setminus \Gamma^-$ then $x\in\Delta_j$ for every $j$, and by \eqref{Cond1} there exists a sequence of points $x_j\in\partial\Gamma^-_j\subset \overline{\Gamma^-}$ converging to $x$, so that $x\in \overline{\Gamma^-}$. 
Similarly, it's easy to check that $\partial\Gamma^-=\partial\Gamma^+=\bigcap_{j=0}^\infty \Delta_j$ and hence that $\Gamma^-=(\Gamma^+)^\circ$.

In proving thickness we recall from Remark \ref{rem:ThickTest} that it is enough to verify the conditions of Definition \ref{def:thicks} for a single value of $j_0$ and with $2^{-j}$ replaced by $\xi^j$ throughout.

To prove $I$-thickness for $\Gamma^-$, fix $c_1,c_2,c_3,c_4>0$ and let $j\geq j_0:=j_*$. Let $Q^e\in \mathcal{Q}((\Gamma^-)^c)$ satisfy $c_1 \xi^{j}\le l(Q^e)\le c_2 \xi^{j}$ and $c_3 \xi^{j}\le \dist(Q^e,\partial\Gamma^-)\le c_4\xi^{j}$, and let $x\in\partial\Gamma^-$ be such that $\dist(Q^e,\partial\Gamma^-)=\dist(Q^e,x)$. Since $x\in \partial\Gamma^-\subset\Delta_j$ there exists $x^-\in\partial\Gamma^-_j$ such that $|x-x^-|\leq c\xi^{j}$. 
By condition \eqref{Cond2}, there exists $Q^i \in \mathcal{Q}(\Gamma^-_j)\subset\mathcal{Q}(\Gamma^-)$ such that $c_1^- \xi^{j}\le l(Q^i)\le c_2^- \xi^{j}$ and $c_3^- \xi^{j}\le \dist(Q^i,\partial \Gamma^-_j) \le \dist(Q^i,x^-) \le c_4^-\xi^{j}$. 
Then, since $\dist(Q^i,\partial\Gamma^-)\geq\dist(Q^i,\partial\Gamma^-_j)$ and $\dist(Q^i,Q^e)\leq \dist(Q^i,x^-) + |x-x^-| + \dist(Q^e,x)$, we see that the definition of $I$-thickness for $\Gamma^-$ is satisfied with $c_5=c_1^-$, $c_6=c_2^-$, $c_7=c_3^-$ and $c_8 = c_4 + c + c_4^-$.

To prove $E$-thickness for $\Gamma^-$, fix $c_1,c_2,c_3,c_4>0$ and let $j\geq j_0:=j_*$. Let $Q^i\in \mathcal{Q}(\Gamma^-)$ satisfy $c_1 \xi^{j}\le l(Q^i)\le c_2 \xi^{j}$ and $c_3 \xi^{j}\le \dist(Q^i,\partial\Gamma^-)\le c_4\xi^{j}$, and let $x\in\partial\Gamma^-$ be such that $\dist(Q^i,\partial\Gamma^-)=\dist(Q^i,x)$. Since $x\in \partial\Gamma^-\subset\Delta_j$ there exists $x^+\in\partial\Gamma^+_j$ such that $|x-x^+|\leq c\xi^{j}$. By condition \eqref{Cond3}, there exists $Q^e \in \mathcal{Q}((\Gamma^+_j)^c)\subset  \mathcal{Q}((\Gamma^-)^c)$ such that $c_1^+ \xi^{j}\le l(Q^e)\le c_2^+ \xi^{j}$ and $c_3^+ \xi^{j}\le \dist(Q^e,\partial \Gamma^+_j) \le \dist(Q^e,x^+) \le c_4^+\xi^{j}$.
Then, since $\dist(Q^e,\partial\Gamma^-)\geq\dist(Q^e,\partial\Gamma^+_j)$ and $\dist(Q^i,Q^e)\leq \dist(Q^e,x^+) + |x-x^+| + \dist(Q^i,x)$, we see that the definition of $E$-thickness for $\Gamma^-$ is satisfied with $c_5=c_1^+$, $c_6=c_2^+$, $c_7=c_3^+$ and $c_8 = c_4 + c + c_4^+$.

The final assertions of the proposition then follow easily from Corollaries \ref{cor:equalnot0} and \ref{cor:thick}.
\end{proof}

We now apply Proposition \ref{prop:thick} to prove thickness for some concrete examples.

\subsection{The classical snowflakes}
\label{subsec:classsnow}

\begin{figure}[htb]
\includegraphics[width=\textwidth, clip, trim=130 130 100 115 ]{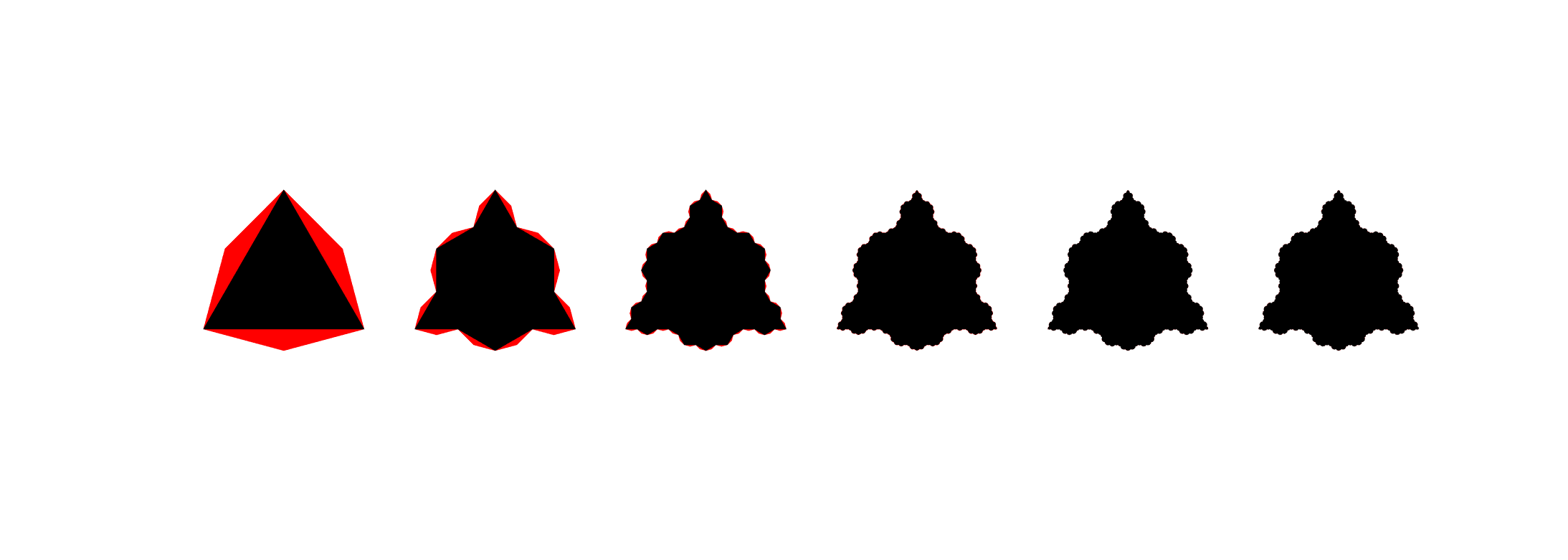}
\includegraphics[width=\textwidth, clip, trim=130 130 100 115 ]{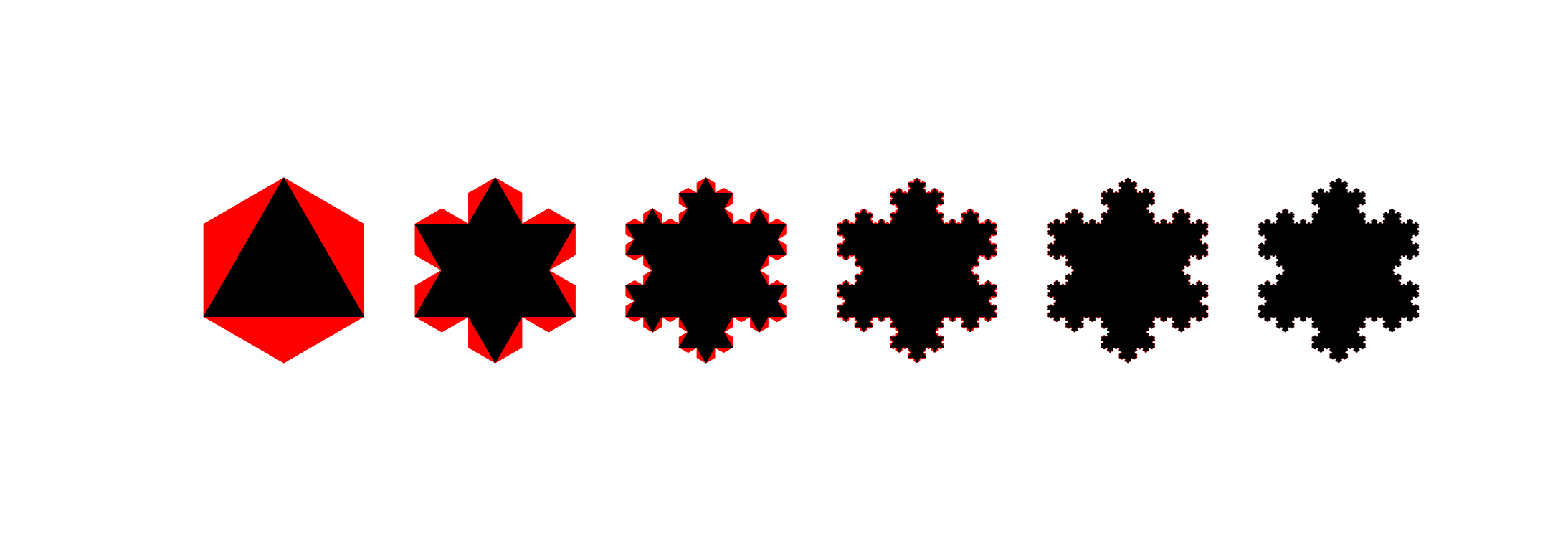}
\includegraphics[width=\textwidth, clip, trim=130 130 100 115 ]{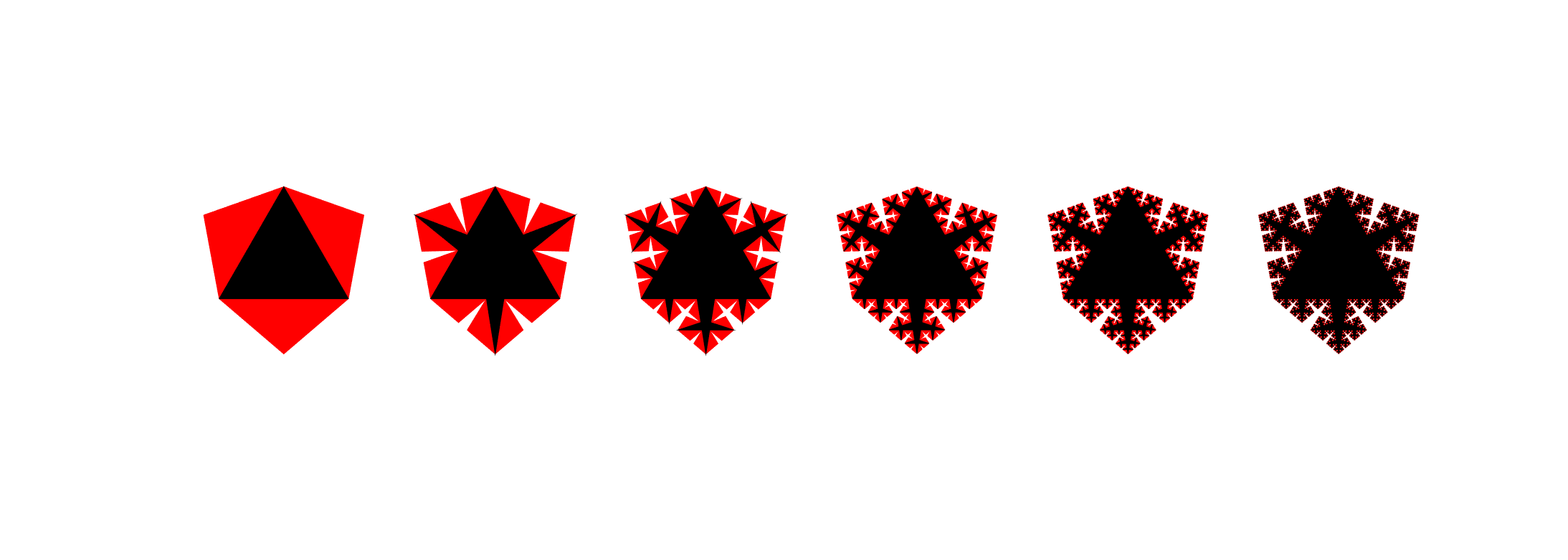}
\caption{The first 6 prefractals $\SN_0^\pm,\ldots,\SN_5^\pm$ of the classical snowflakes for $\beta=\frac\pi3$ (top), $\beta=\frac\pi6$ (centre), $\beta=\frac\pi{20}$ (bottom).
The inner prefractals $\SN_j^-$ are the black shapes, the outer ones $\SN_j^+$ are the union of the red and the black shapes, and the differences $\Delta_j$ are the red parts.
}
\label{fig:SnowflakesShapesInnerOuter}
\end{figure}

We first consider the family of ``classical snowflakes'' studied in \cite{EvansPhD}, which generalise the standard Koch snowflake. 
These snowflakes are open subsets of $\R^2$, defined as limits of nested (increasing) sequences of open polygonal prefractals.  
In order to apply Proposition \ref{prop:thick} and deduce thickness, we need to introduce a sequence of nested (decreasing) closed prefractals, which generalise those considered in \cite{banjai2017poincar} for the standard Koch snowflake. 
The interior and exterior prefractals for three examples (including the Koch snowflake) are shown in Figure \ref{fig:SnowflakesShapesInnerOuter}. 

The snowflakes are parametrised by a number $0<\SNangle<\pi/2$, which represents half the width of each convex angle of the interior prefractals (except possibly the three angles of the first interior prefractal).  Given 
$0<\SNangle<\frac\pi2$ 
we define 
$\xi:=\frac1{2(1+\sin\SNangle)}$, which satisfies $\frac14<\xi<\frac12$ and represents the ratio of the side lengths of two successive prefractals.
The standard Koch snowflake corresponds to the choice $\SNangle=\pi/6$, so that $\xi=1/3$. We note that $\xi$ is denoted $\alpha^{-1}$ in \cite[\S1.1]{EvansPhD}.

We define an increasing sequence of nested polygons $\SN_j^-$, $j\in\N_0$, as follows.
Each $\SN_j^-$ is an open polygon with $M_j^-:=3\cdot4^j$ edges of length $\xi^j$.
Let $\SN_0^-$ be the equilateral triangle with vertices $(0,0)$, $(1,0)$, $(\frac12,\frac12\sqrt3)$.
Then $\SN_j^-$ is the union of $\SN_{j-1}^-$ and $M_{j-1}^-$ identical disjoint isosceles triangles (together with their bases) with basis length $\xi^{j-1}(1-2\xi)$, side length $\xi^j$, height $\xi^{j-1}\sqrt{\xi-\frac14}$, apex angle $2\SNangle$, disjoint from $\SN_{j-1}^-$ and placed in such a way that  the midpoint of the basis of the $k$th such triangle coincides with the midpoint of the $k$th side of $\SN_{j-1}^-$, for $k=1,\ldots,M_{j-1}^-$. 

The external closed polygons $\SN_j^+$, $j\in\N_0$, are defined as follows.
Each $\SN_j^+$ is a closed polygon with $M_j^+:=6\cdot4^j$ edges of length $\xi^{j+\frac12}$.
The first one $\SN_0^+$ is the convex hexagon obtained as union of $\SN_0^-$ and the three isosceles closed triangles with base the three sides of $\SN_0^-$, respectively, and height $\sqrt{\xi-\frac14}$ ($\SN_0^+$ is a regular hexagon only if $\beta=\frac\pi6$).
Then $\SN_j^+$ is the difference of $\SN_{j-1}^+$ and $M_{j-1}^+$ identical disjoint isosceles triangles (together with their bases) with basis length $\xi^{j-\frac12}(1-2\xi)$, side length $\xi^{j+\frac12}$, height $\xi^{j-\frac12}\sqrt{\xi-\frac14}$, apex angle $2\SNangle$, contained in $\SN_{j-1}^+$ and placed in such a way that the midpoint of the basis of the $k$th such triangle coincides with the midpoint of the $k$th side of $\SN_{j-1}^+$, for $k=1,\ldots,M_{j-1}^+$. 

The prefractals satisfy $\SN_j^-\subset\SN_j^+$, $\SN_j^-\subset\SN_{j+1}^-$ and $\SN_{j+1}^+\subset\SN_{j}^+$, as required in the framework for Proposition \ref{prop:thick}.
The limit snowflakes are defined as $\SN^-:=\bigcup_{j\in\N_0}\SN_j^-$ and $\SN^+:=\bigcap_{j\in\N_0}\SN_j^+$ and the boundary approximations are $\Delta_j:=\SN_j^+\setminus\SN_j^-$ (the red parts in Figure~\ref{fig:SnowflakesShapesInnerOuter}).

\begin{prop}\label{prop:SNthick3propr}
For every $0<\SNangle<\frac\pi2$, the classical snowflake domain $\SN^-$ is thick, with $\overline{\SN^-}=\SN^+$, 
and $\widetilde{A}^s_{p,q}(\Gamma^-)=A^s_{p,q,\overline{\Gamma^-}}=A^s_{p,q,\Gamma^+}$ for all $s\in\R\setminus\{0\}$ and $1<p,q<\infty$, and $\widetilde{H}_{p}^{s}(\Gamma^-)=H_{p,\overline{\Gamma^-}}^{s}=H_{p,{\Gamma^+}}^{s}$ for all $s\in\R$ and $1<p<\infty$.
\end{prop}
\begin{proof}
We prove that the sequences $\SN^\pm_j$ satisfy the assumptions of Proposition~\ref{prop:thick}. 
We first note that since $\partial\SN^-\subset\Delta_j$ and $|\Delta_j|=4\xi^2|\Delta_{j-1}|=(2\xi)^{2j}|\Delta_0|=\frac32(2\xi)^{2j}\sqrt{\xi-\frac14}\xrightarrow{j\to\infty}0$ by $\xi<\frac12$, we have $|\partial\SN^-|=0$.
Next we verify the three conditions \eqref{Cond1}, \eqref{Cond2} and \eqref{Cond3}. 
To that end we choose $\xi$ in \eqref{Cond1}--\eqref{Cond3} to be the $\xi$ in the definition of $\SN_j^\pm$, and fix $j_0=1$.

To prove \eqref{Cond1}, take any $x\in\Delta_j$.
By definition of the prefractals, $x\in T$, where $T$ is an isosceles triangle with base length $\xi^j$, height $\xi^j\sqrt{\xi-\frac14}$, base contained in $\partial\SN_j^-$ and legs in $\partial\SN_j^+$.
Thus $\dist(x,\partial\SN_j^-)\le \xi^j\sqrt{\xi-\frac14}$ and $\dist(x,\partial\SN_j^+)\le \frac12\xi^j$ and so \eqref{Cond1} holds with $c=\frac12$.

\begin{figure}[htb]\centering
\begin{tikzpicture}[scale=5]
\fill[black!20!white] (0,0)--(0.216,-0.184)--(0.432,0)--(0.283,-0.242)--(0.5,-0.427)--(0.716,-0.242)--(0.567,-0)--(0.783,-0.185)--
  (1,0)--(1.051,0.28)--(0.783,0.374)--(1.068,0.367)--(1.119,0.646)--(0.851,0.741)--(0.716,0.491)--(0.768,0.771)--(0.5,0.866)--  (0.232,0.771)--(0.283,0.492)--(0.148,0.741)--(-0.112,0.646)--(-0.068,0.37)--(0.216,0.374)--(-0.052,0.279)--(0,0);
\fill[black!50!white] (0,0)--(0.432,0)--(0.5,-0.427)--(0.567,0)--(1,0)--(0.783,0.374)--(1.119,0.646)--(0.716,0.491)--
  (0.5,0.866)--(0.283,0.491)--(-0.12,0.646)--(0.216,0.374)--(0,0);
\draw(.5,.5)node{$\SN_j^-$};\draw(.13,.62)node{$\SN_j^+$};
\draw[very thick] (0.768,0.771)--(0.851,0.741)--(0.716,0.491)--(0.768,0.771);\draw(.84,.8)node{$T^+$};
\draw (.795,.72)circle(.035);
\fill (.795+.012,.72+.012)--(.795-.012,.72+.012)--(.795-.012,.72-.012)--(.795+.012,.72-.012);
\draw[->](.9,.7)--(.795+.012,.72);\draw(.95,.7)node{$Q^e$};
\fill(0.7835,0.616)circle(.012);\draw(.83,.6)node{$x^+$};
\fill(0.2,0)circle(.012);\draw(.2,-.03)node{$x^-$};
\draw[very thick] (0.432,0)--(0.5,-0.427)--(0.567,0)--(0.432,0);\draw(.6,-.2)node{$T^-$};
\draw (.5,-.058)circle(.055);
\fill (.5+.016,-.058+.016)--(.5+.016,-.058-.016)--(.5-.016,-.058-.016)--(.5-.016,-.058+.016);
\draw[->](.41,-.16)--(.5-.016,-.058-.016);\draw(.38,-.15)node{$Q^i$};
\fill(0.432,0)circle(.012);\draw(.4,.04)node{$\tilde x^-$};
\draw[<->](0.96,0.02)--(0.77,0.34);
\draw (0.8,0.15) node{$\xi^j$};
\draw[<->](1.03,0.02)--(1.07,0.25);
\draw (1.18,0.14) node{$\xi^{j+1/2}$};
\end{tikzpicture}
\caption{\label{fig:ClassicalQiQe}
A schematic representation of the proof of conditions \eqref{Cond2} and \eqref{Cond3} for the classical snowflake $\SN^-$ in the proof of Proposition~\ref{prop:SNthick3propr}.
Given $x^-\in\partial\SN_j^-$ we construct $Q^i\in \mathcal Q(\SN_j^-)$; given $x^+\in\partial\SN_j^+$ we construct $Q^e\in \mathcal Q((\SN_j^+)^c)$.
In this example $j=1$, $\beta=\pi/20$ and $\tilde x^+=x^+$.}
\label{fig:SNtria}
\end{figure}
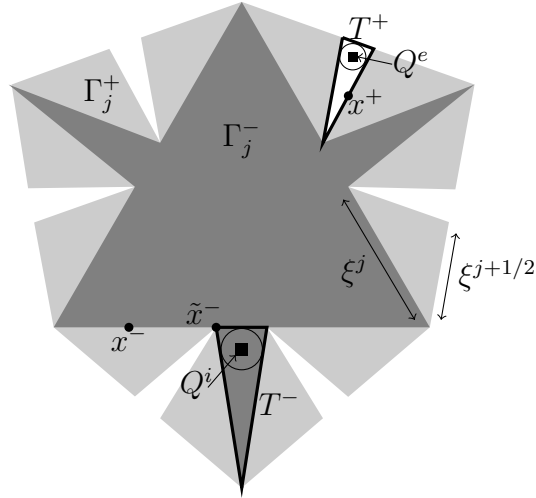

Now take $x^-\in\partial\SN_j^-$. 
Since $j\geq j_0=1$ there exists a connected component $T^-$ of $\SN_j^-\setminus\overline{\SN_{j-1}^-}$ and a point $\tilde{x}^-\in\partial T^-$ such that $|x^--\tilde{x}^-|\le\xi^j$.
By the construction of $\SN_j^-$, $T^-$ is an open isosceles triangle $T^-\subset\SN_j^-$ with leg length $\xi^j$ and apex angle $2\beta$.
(In particular if $x^-\in\partial \SN_j^-\setminus\partial\SN_{j-1}^-$ then we can take $\tilde{x}^-=x^-$.) 
For an illustration see Figure \ref{fig:ClassicalQiQe}. 
The triangle $T^-$ contains an open disc of radius $\rho_\beta\xi^j$ (whose boundary is the inscribed circle) where
$\rho_\beta:=\sin\beta\tan(\frac\pi4-\frac\beta2)>0$ 
depends only on $\beta$, and inside this disc we can construct a square $Q^i \in \mathcal Q (\Gamma_j^-)$ of side length $\frac{\rho_\beta}{\sqrt2}\xi^j$ sharing the same centre as the disc. Then $\frac{\rho_\beta}{2}\xi^j\le \dist(Q^i,\partial T^-) \le \dist(Q^i,\partial\SN_j^-)\le \dist(Q^i,x^-)\le |x^--\tilde{x}^-|+\dist(Q^i,\tilde{x}^-)\le 2\xi^j$, so that \eqref{Cond2} holds with $c_1^-=c_2^-=\frac{\rho_\beta}{\sqrt2}$, $c_3^-=\frac{\rho_\beta}{2}$ and $c_4^-=2$.

Similarly, for any $x^+\in\partial\SN_j^+$ there exists an open isosceles triangle $T^+\subset (\SN_{j-1}^+)^\circ\setminus\SN_j^+\subset(\SN_j^+)^c$ with leg length $\xi^{j+\frac12}$ and apex angle $2\beta$ and a point $\tilde{x}^+\in\partial T^+$ such that $|x^+-\tilde{x}^+|\le\xi^{j+1/2}$. Again, for an illustration see Figure \ref{fig:ClassicalQiQe}.  The same reasoning as above permits the construction of $Q^e\in \mathcal Q((\SN_j^+)^c)$ such that \eqref{Cond3} holds with $c_m^+=\xi^{\frac12}c_m^-$ for $m=1,\ldots,4$.
\end{proof}

The next result will be useful in section \ref{sec:Density}, where we deal with $d$-sets (cf. Definition \ref{def:dSet}).

\begin{prop}\label{prop:class-snow-d-set}
The boundaries of the classical snowflakes introduced above and parametrized by $\beta \in (0,\frac{\pi}{2})$ are $d$-sets for $d=-\frac{\log 4}{\log \xi}$, where $\xi=\frac{1}{2(1+\sin \beta)}$ as before.
\end{prop}

\begin{proof}
Let $\beta$, $\xi$ and $d$ be as in the statement of the proposition.

\emph{Step 1.} Since a finite union of $d$-sets is clearly still a $d$-set, it is enough to prove that the part of the boundary {built over} each one of the three legs of the initial equilateral triangle is a $d$-set. And since the Hausdorff measure is invariant under translations and rotations, we shall do the forthcoming analysis after a rigid motion has been performed in such a way that each leg of the initial triangle coincides with the segment $[(0,0),(1,0)]$ in $\R^2$ and the corresponding part of the boundary lies above it. Our objective is then to prove that this is a $d$-set.

\emph{Step 2.} We shall use the same notation as before, except that we prepend the fraction $\frac{1}{3}$ to it. So, the part of the boundary to be considered is denoted $\frac{1}{3}\partial \Gamma^-$ and equals $\bigcap_{j=0}^\infty \frac{1}{3}\Delta_j$ (as in the first paragraph of the proof of Proposition~\ref{prop:thick}), where $\frac{1}{3}\Delta_j$ stands for the part of the boundary approximation $\Delta_{j}$ built only over the segment $[(0,0),(1,0)]$. From the way \eqref{Cond1} is proved in Proposition \ref{prop:SNthick3propr}, we see that also the following holds:
$$\forall x\in {\textstyle \frac{1}{3}}\Delta_j, \;\; \dist(x,{\textstyle \frac{1}{3}}\partial\Gamma^-_j) \leq {\textstyle \frac{1}{2}}\xi^{j}\, \text{ and }\, \dist(x,{\textstyle \frac{1}{3}}\partial\Gamma^-)\leq {\textstyle \frac{1}{\sqrt{2}}}\xi^{j}.$$
The second inequality is due to the fact that in the isosceles triangle $T$ mentioned in the proof the end points of the base belong to ${\textstyle \frac{1}{3}}\Delta_{j'}$ for every $j'$, and therefore to ${\textstyle \frac{1}{3}}\partial\Gamma^-$. Then, similarly as observed just before Proposition \ref{prop:thick}, $d_H({\textstyle \frac{1}{3}}\partial\Gamma_j^-,{\textstyle \frac{1}{3}}\partial\Gamma^-)\leq {\textstyle \frac{1}{\sqrt{2}}}\xi^j$, therefore
\begin{equation}\label{eq:convclassnow}
{\textstyle \frac{1}{3}}\partial\Gamma_j^- \xrightarrow[j \to \infty]{} {\textstyle \frac{1}{3}}\partial\Gamma^-
\end{equation}
in the Hausdorff metric in the space of non-empty compact subsets of $\R^2$.

\emph{Step 3.} We are going to show now that ${\textstyle \frac{1}{3}}\partial\Gamma^-$ is also the fractal (invariant set) determined by four contractions $\psi_i$, $i=1,\ldots,4$, in $\R^2$ according to \cite[Thm.~4.2]{Tri97} and that these contractions are indeed similarities (similitudes) with contraction ratio equal to $\xi$ and satisfy the open set condition of \cite[Def.~4.5(ii)]{Tri97}. Afterwards, by \cite[Thm.~4.7]{Tri97} we can conclude that ${\textstyle \frac{1}{3}}\partial\Gamma^-$ is a compact $D$-set with $D\geq 0$ determined by $\sum_{i=1}^4 \xi^D =1$, from which it follows that $D=d$, finishing the proof.

The mentioned contractions are defined as follows, where $H(\xi)$ denotes the homothety with centre at the origin and ratio $\xi$, $R(\theta)$ denotes  counterclockwise rotation through angle $\theta$ about the origin, and $T(x,y)$ denotes translation by a vector $(x,y)$: 
\begin{align*}
\psi_1 &= H(\xi), &
\psi_2 &= T(\xi,0) \circ R({\textstyle\frac\pi2}-\beta) \circ H(\xi),\\
\psi_3 &= T(1-\xi-\xi\sin\beta,\xi\cos\beta) \circ R(\beta-{\textstyle\frac\pi2}) \circ H(\xi), &
\psi_4 &= T(1-\xi,0) \circ H(\xi).
\end{align*}	

These are, clearly, similarities of ratio $\xi$ and determine, according to \cite[Thm.~4.2]{Tri97}, the unique non-empty compact set $K$ in $\R^2$ such that
$$K=\psi(K):=\bigcup_{i=1}^4 \psi_i(K).$$
Still according to \cite[Thm.~4.2]{Tri97}, $K$ can be obtained as
\begin{equation}\label{eq:Kaslimit}
K=\lim_{j\to\infty} \psi^j(\Lambda) := \lim_{j\to\infty} (\underbrace{\psi \circ \ldots \circ \psi}_j)(\Lambda)
\end{equation}
for any non-empty compact subset $\Lambda$ of $\R^2$,
the limit being taken in the metric space of all non-empty compact sets in $\R^2$ equipped with the Hausdorff metric.

Since each $\psi_i$ maps an edge of $\frac13\partial\Gamma_{j-1}^-$ to one of $\frac13\partial\Gamma_{j}^-$, choosing $\Lambda=[(0,0),(1,0)]$, it is easy to see that
$$\psi^j(\Lambda)={\textstyle \frac{1}{3}}\partial\Gamma_j^-, \quad j\in \N.$$
Combining this with \eqref{eq:convclassnow} and \eqref{eq:Kaslimit}, we get that
$$K={\textstyle \frac{1}{3}}\partial\Gamma^-.$$

In order to finish the proof, it only remains to exhibit a non-empty open set $O$ in $\R^2$ such that
$$\psi(O)=\bigcup_{i=1}^4 \psi_i(O)\subset O \quad \text{and} \quad \psi_i(O)\cap \psi_k(O)=\emptyset \; \text{ for }\; i\not=k.$$
It is easily seen that we can take for $O$ the interior of $\frac{1}{3}\Delta_0$.
\end{proof}

\begin{rem}\label{rem:Hdim}
Combining the above result with the information given after Definition \ref{def:dSet}, we have that the Hausdorff dimension of $\partial\SN^-$ is $-\log 4/\log\xi$, with the boundary of the standard Koch snowflake ($\SNangle=\frac\pi6$, $\xi=\frac13$) having dimension $-\log 4/\log(1/3)$. Moreover, since $\xi$ ranges over all values in $(\frac{1}{4},\frac{1}{2})$, we have produced a class of domains in $\R^2$ whose boundaries have Hausdorff dimensions ranging over all values in $(1,2)$.   
\end{rem}

\subsection{The square snowflake}\label{subsec:squaresnow}
\begin{figure}[htb]
\includegraphics[width=\textwidth, clip, trim = 25 90 0 90]{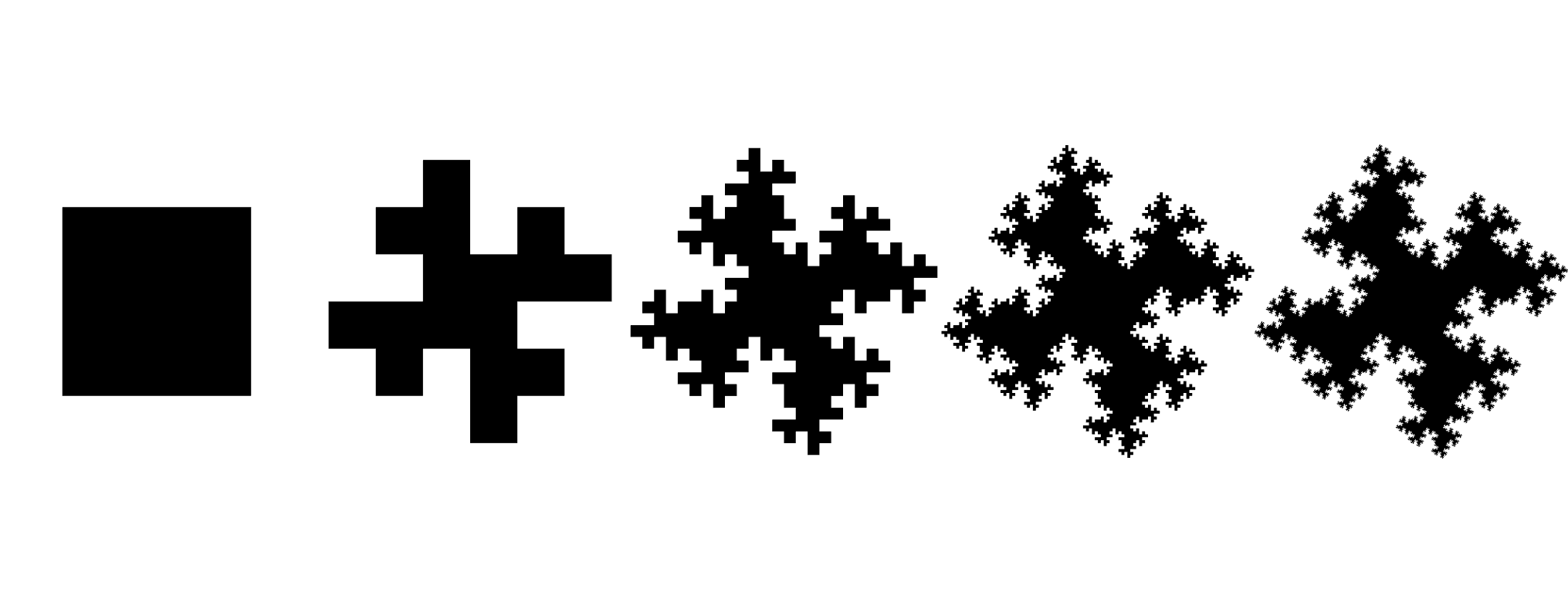}
\caption{The first five prefractals $\SqS_0,\ldots,\SqS_4$ of the square snowflake.}
\label{fig:SquareSnowflakeShapesBlack}
\end{figure}

We now consider the ``square snowflake'' studied in \cite{sapoval1991vibrations} (see also \cite[\S7.6]{grebenkov2013geometrical} and the references therein). 
Like the classical snowflakes studied in the previous section, this is an open set $\SqS^- \subset\R^2$ with fractal boundary. 
The starting point for the definition of $\SqS^-$ is a sequence of \textit{non-nested} polygonal prefractals $\SqS_j$, $j\in\N_0$, the first five of which are shown in Figure \ref{fig:SquareSnowflakeShapesBlack}. 

The sequence of prefractals $\SqS_j\subset\R^2$, $j=0,1,\ldots$ is defined as follows.
Each prefractal $\SqS_j$ is a polygon whose boundary is the union of $N_j:=4\cdot 8^j$ segments of length $\ell_j:=4^{-j}$ aligned to the Cartesian axes.
Let $\SqS_0=(0,1)^2$ be the unit open square.
For $j\in\N$, $\partial\SqS_j$ is constructed by replacing each horizontal edge $[(x,y),(x+\ell_{j-1},y)]$ and each vertical edge $[(x',y'),(x',y'+\ell_{j-1})]$ of $\partial\SqS_{j-1}$ respectively by the following polygonal lines composed of 8 edges each:
\begin{equation}\label{eq:SqSrule}
   \begin{aligned}
\left[\begin{aligned}(x,&y)\\(x+\ell_{j-1},&y)
\end{aligned}\right]\rightsquigarrow
\\[10mm]\mmbox{\begin{tikzpicture}[scale=.3]
\draw[very thick](-5,0)--(-2,0);
\draw[very thick](0,0)--(1,0)--(1,1)--(2,1)--(2,-1)--(3,-1)--(3,0)--(4,0);
\draw(-1,0)node{$\rightsquigarrow$};
\end{tikzpicture}}\hspace{5mm}
   \end{aligned}
\left[\begin{aligned}
(x,&y)\\
(x+\ell_j,&y)\\
(x+\ell_j,&y+\ell_j)\\
(x+2\ell_j,&y+\ell_j)\\
(x+2\ell_j,&y)\\
(x+2\ell_j,&y-\ell_j)\\
(x+3\ell_j,&y-\ell_j)\\
(x+3\ell_j,&y)\\
(x+\ell_{j-1},&y)
\end{aligned}\right],\qquad
   \begin{aligned}
\left[\begin{aligned}(x',&y')\\(x',&y'+\ell_{j-1})\end{aligned}\right]\rightsquigarrow
\\[10mm]\mmbox{\;\begin{tikzpicture}[scale=.3]
\draw[very thick](-3,-2)--(-3,2);
\draw[very thick](1,-2)--(1,-1)--(0,-1)--(0,0)--(2,0)--(2,1)--(1,1)--(1,2);
\draw(-1.5,0)node{$\rightsquigarrow$};
\end{tikzpicture}}\hspace{10mm}
   \end{aligned}
\left[\begin{aligned}
(x',&y')\\
(x',&y'+\ell_j)\\
(x'-\ell_j,&y'+\ell_j)\\
(x'-\ell_j,&y'+2\ell_j)\\
(x',&y'+2\ell_j)\\
(x'+\ell_j,&y'+2\ell_j)\\
(x'+\ell_j,&y'+3\ell_j)\\
(x',&y'+3\ell_j)\\
(x',&y'+\ell_{j-1})
\end{aligned}\right].
\end{equation}
(Note that the fourth and the fifth segments obtained are aligned; in the following however we count them as two different edges of $\SqS_j$.)
Each polygonal path $\partial\SqS_j$ constructed with this procedure is the boundary of a simply connected polygon $\SqS_j$ of unit area, composed of $16^j$ squares of side length $\ell_j$.
We note that the closures of the prefractals tile the plane:
$\bigcup_{(k_1,k_2)\in\Z^2}(\overline\SqS_j+(k_1,k_2))=\R^2$ 
for any $j$ and $(\SqS_j+(k'_1,k'_2))\cap(\SqS_j+(k_1,k_2))=\emptyset$ for all $(k'_1,k'_2),(k_1,k_2)\in\Z^2$ such that $ (k'_1,k'_2)\not=(k_1,k_2) $.

The resulting sequence of prefractals $(\Gamma_j)_{j\in\N_0}$ is not nested: for each $j\in\N$ neither $\SqS_j\subset\SqS_{j-1}$ nor $\SqS_j\supset\SqS_{j-1}$.
Indeed, the two set differences $\SqS_j\setminus\SqS_{j-1}$ and $\SqS_{j-1}\setminus\SqS_j$ are made of $4\cdot 8^{j-1}=2^{3j-1}$ disjoint squares of side length $\ell_j$.
Thus the limit set $\SqS^-$ cannot be defined simply as the infinite union or intersection of the prefractals defined above.  
However, we are going to construct, as before, two nested sequences $(\Gamma^\pm_j)_{j}$ of open and closed prefractals approximating monotonically an open set $\Gamma^-$ and its closure, as in Proposition \ref{prop:thick}, such that the boundary $\partial \Gamma^-$ of $\Gamma^-$ is the limit, in the Hausdorff metric, of the non-nested prefractal boundaries $(\partial\Gamma_j)_{j}$ (cf.\ eq.~\eqref{eq:convsquaresnow} and Proposition \ref{prop:square-snow-d-set}).

We first denote by $E_{j,k}$, $k=1,\ldots,N_j$, the sides of $\partial\SqS_j$, each of which has length $\ell_j$.
Then let $S_{j,k}$, $k=1,\ldots,N_j$, be the closed squares with diagonals $E_{j,k}$, respectively; they have disjoint interiors and are tilted at 45 degrees to the Cartesian axes.
We then define the set $\Delta_j:=\bigcup_{k=1}^{N_j}S_{j,k}$, which is compact with 
Lebesgue measure equal to $4 \cdot 8^j \cdot \frac12 \ell_j^2=2^{1-j}$.
The relevance of this construction is the following: given an edge $E_{j,k}=[(x,y),(x+\ell_j,y)]$ (or $E'_{j,k}=[(x,y),(x,y+\ell_j)]$) of $\partial \SqS_j$, its ``evolution'', i.e.\ all the segments obtained from the successive applications of the rules in \eqref{eq:SqSrule}, are contained in the closed square $S_{j,k}$ with vertices $(x,y),(x+\frac12\ell_j,y-\frac12\ell_j),(x+\ell_j,y),(x+\frac12\ell_j,y+\frac12\ell_j)$ (or $(x,y),(x+\frac12\ell_j,y+\frac12\ell_j),(x,y+\ell_j),(x-\frac12\ell_j,y+\frac12\ell_j)$, respectively), which is one of the squares composing $\Delta_j$.
This implies that these sets are nested and contain the boundaries of the prefractals of higher order: $\partial \SqS_{j'}\subset\Delta_{j'}\subset\Delta_j$ for all $j'\ge j$.

\begin{figure}[htb]
\includegraphics[width=\textwidth,clip,trim=5 80 10 80]{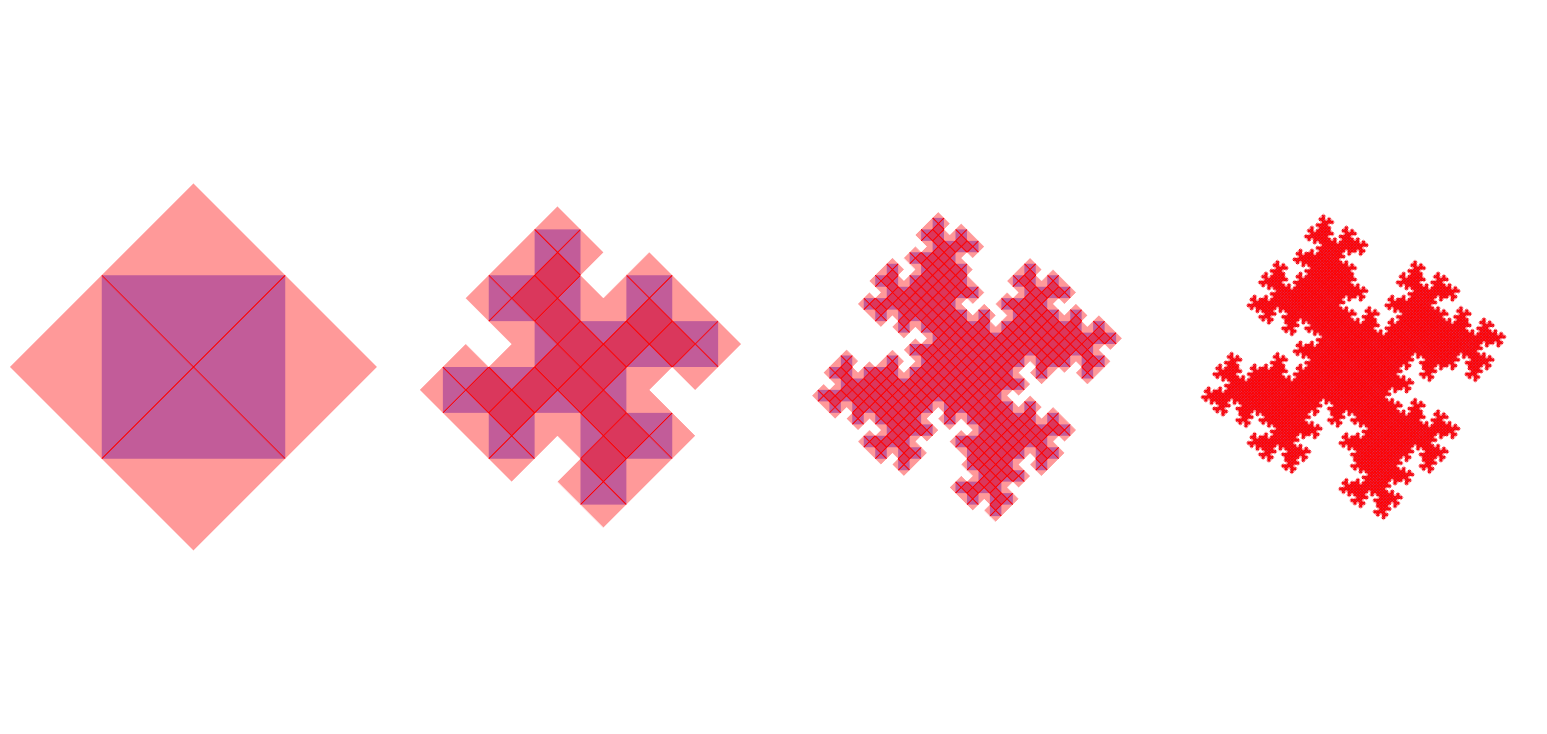}
\caption{The non-monotonic square snowflake prefractals $\SqS_j$ (union of the blue and red parts), 
the boundary approximations $\Delta_j$ (union of the blue and pink parts), 
the inner prefractals $\SqS_j^-$ (red part only), 
and the outer prefractals $\SqS_j^+$ (union of the red, blue and pink parts),
for $j=0,1,2,3$.}
\label{fig:SquareSnowflakeWithUpsilonShapes}
\end{figure}

We now define two sequences of open and closed polygons, respectively:
$$\SqS_j^-:=\SqS_j\setminus\Delta_j,\qquad\SqS_j^+:=\SqS_j\cup \Delta_j,\qquad j=0,1,2\ldots.$$
The open inner prefractals $\SqS_j^-$ are nested and increasing, the closed outer prefractals $\SqS_j^+$ are nested and decreasing, and they approximate from inside and outside the non-monotonic prefractals $\SqS_j$: 
$$\SqS_j^-\subset\SqS_{j+1}^-,\qquad
\SqS_j^+\supset\SqS_{j+1}^+,\qquad
\SqS_j^-\subset\SqS_j\subset\SqS_j^+,\qquad
\partial\SqS_j\subset\Delta_j=\SqS_j^+\setminus\SqS_j^-.$$
This monotonicity implies that we can define two limits 
$\SqS^-:=\bigcup_{j\in\N_0} \SqS_j^-\subset\SqS^+:=\bigcap_{j\in\N_0} \SqS_j^+$
and that they are open and closed, respectively. 
The inner prefractals $\SqS_j^-$, the outer prefractals $\SqS_j^+$ and the boundary approximations $\Delta_j$ are shown in Figure \ref{fig:SquareSnowflakeWithUpsilonShapes}.
(Note that $\SqS_0^-=\emptyset$.) 

Having defined $\Gamma^-$, we are now in a position to prove that it is thick using Proposition \ref{prop:thick}.

\begin{figure}[htb]\centering
\begin{tikzpicture}[scale=3.5]
\def\Incr{1/8};
\draw[fill=black!20!white](0,0)--(\Incr,-\Incr)--(3*\Incr,\Incr)--(.5,0)--(3*\Incr,-\Incr)--(5*\Incr,-3*\Incr)--(9*\Incr,\Incr)
--(7*\Incr,3*\Incr)--(1,0.5)--(9*\Incr,3*\Incr)--(11*\Incr,5*\Incr)--(7*\Incr,9*\Incr)--(5*\Incr,7*\Incr)--(.5,1)
--(5*\Incr,9*\Incr)--(3*\Incr,11*\Incr)--(-\Incr,7*\Incr)--(\Incr,5*\Incr)--(0,.5)--(-\Incr,5*\Incr)--(-3*\Incr,3*\Incr)--(0,0);
\draw[thick,fill=black!40!white](0,0)--(1/4,0)--(1/4,1/4)--(1/2,1/4)--(1/2,0)--(1/2,-1/4)--(3/4,-1/4)--(3/4,0)--(1,0)
--(1,1/4)--(3/4,1/4)--(3/4,1/2)--(1,1/2)--(5/4,1/2)--(5/4,3/4)--(1,3/4)--(1,1)
--(3/4,1)--(3/4,3/4)--(1/2,3/4)--(1/2,1)--(1/2,5/4)--(1/4,5/4)--(1/4,1)--(0,1)
--(0,3/4)--(1/4,3/4)--(1/4,1/2)--(0,1/2)--(-1/4,1/2)--(-1/4,1/4)--(0,1/4)--(0,0);
\draw[fill=black!60!white](1/8,1/8)--(3/8,3/8)--(5/8,1/8)--(1/2,0)--(5/8,-1/8)--(7/8,1/8)--(5/8,3/8)--(7/8,5/8)--(1,1/2)
--(9/8,5/8)--(7/8,7/8)--(5/8,5/8)--(3/8,7/8)--(1/2,1)--(3/8,9/8)--(1/8,7/8)--(3/8,5/8)--(1/8,3/8)--(0,1/2)--(-1/8,3/8)--(1/8,1/8);
\fill(.5+\Incr+0.03,3/4+\Incr+0.03)circle(.02);\draw(.68,.85)node{$x^+$};
\fill(5/8-1/48,1-1/48)--(5/8-1/48,1+1/48)--(5/8+1/48,1+1/48)--(5/8+1/48,1-1/48);
\draw[->](.73,1.13)--(.64,1.02);\draw(.8,1.15)node{$Q^e$};
\draw[ultra thick](6*\Incr,1)--(5*\Incr,7*\Incr)--(.5,1)--(5*\Incr,9*\Incr)--(6*\Incr,1);
\draw[ultra thick](5/8,1/8)--(6/8,2/8)--(7/8,1/8)--(6/8,0/8)--(5/8,1/8);
\fill(6/8-1/48,1/8-1/48)--(6/8-1/48,1/8+1/48)--(6/8+1/48,1/8+1/48)--(6/8+1/48,1/8-1/48);
\draw[->](1,-.1)--(6/8+1/48,1/8-1/48);\draw(1.05,-.15)node{$Q^i$};
\draw[decorate,decoration={brace,amplitude=5},thin](5/4,1)--(5/4,3/4);\draw(1.35,.85)node{$\ell_j$};
\fill(7/8-0.03,1/8+0.03)circle(.02);\draw(.94,.18)node{$x^-$};
\draw(0.02,3/8-0.02)node{$\SqS_j^-$};
\draw(-2.9*0.0625,3/8-0.045)node{$\SqS_j$};
\draw(-1.1*0.0625,0.165)node{$\SqS_j^+$};
\end{tikzpicture}
\caption{A schematic representation of the proof of conditions \eqref{Cond2} and \eqref{Cond3} for the square snowflake $\SqS^-$ in the proof of Proposition~\ref{prop:SqSthick3propr}.
Given $x^-\in\partial\SqS_j^-$ we construct $Q^i\in \mathcal Q(\SqS_j^-)$; given $x^+\in\partial\SqS_j^+$ we construct $Q^e\in \mathcal Q((\SqS_j^+)^c)$.
This illustration shows the case $j=1$.}
\label{fig:SqScond23}
\end{figure}
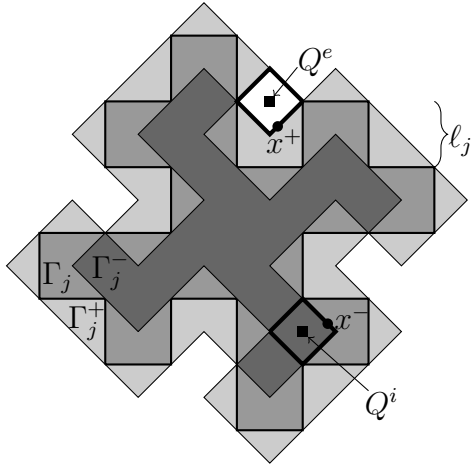
\begin{prop}\label{prop:SqSthick3propr} 
The square snowflake domain $\SqS^-$ is thick, with $\overline{\SqS^-}=\SqS^+$ 
and $\widetilde{A}^s_{p,q}(\Gamma^-)=A^s_{p,q,\overline{\Gamma^-}}=A^s_{p,q,\Gamma^+}$ for all $s\in\R\setminus\{0\}$ and $1<p,q<\infty$, and $\widetilde{H}_{p}^{s}(\Gamma^-)=H_{p,\overline{\Gamma^-}}^{s}=H_{p,{\Gamma^+}}^{s}$ for all $s\in\R$ and $1<p<\infty$.
\end{prop}
\begin{proof}
Again, we show that the sequences $\SqS^\pm_j$ satisfy the assumptions of Proposition~\ref{prop:thick}. First we note that since $\partial\SqS^-\subset\Delta_j$ and $|\Delta_j|=2^{1-j}$ for all $j\in\N$, we have $|\partial\SqS^-|=0$. 
Next we verify the three conditions \eqref{Cond1}, \eqref{Cond2} and \eqref{Cond3}. 
To that end we set $\xi=1/4$ and $j_0=1$.

To prove \eqref{Cond1}, take any $x\in\Delta_j$.
By the definition of $\Delta_j$, $x\in S_{j,k}$, for some $k\in\{1,\ldots,N_j\}$.
At least one vertex of the tilted square $S_{j,k}$ belongs to $\partial \SqS_j^-$ and at least one to $\partial \SqS_j^+$, so $\max\{\dist(x,\partial\SqS_j^-)$, $\dist(x,\partial\SqS_j^+)\}$ $\le \xi^j$, which is the length of the diagonal of $S_{j,k}$.
Thus condition \eqref{Cond1} holds with $c=1$.

Now take $x^-\in\partial\SqS_j^-$.
Since $\SqS_j^-$ is the union of tilted squares of side $\frac1{\sqrt2}\xi^j$, $x^-$ lies on the boundary of one of these squares.
Denote by $Q^i\in\mathcal Q(\SqS_j^-)$ the square with side length $\frac14\xi^j$ (aligned to the Cartesian axes) centred at the centre of this tilted square. For an illustration see Figure \ref{fig:SqScond23}. 
Then it is elementary to check that condition \eqref{Cond2} holds with $c_1^-=c_2^-=\frac14$, $c_3^-=\frac1{4\sqrt2}$ and $c_4^-=\frac38$.

Similarly, any $x^+\in\partial\SqS_j^+$ lies on the boundary of a tilted square contained in $(\SqS_j^+)^c$, so the same reasoning permits the construction of $Q^e\in \mathcal Q((\SqS_j^+)^c)$ such that \eqref{Cond3} holds with $c_m^+=c_m^-$ for $m=1,\ldots,4$. Again, for an illustration see Figure \ref{fig:SqScond23}. 
\end{proof}

As with Proposition \ref{prop:class-snow-d-set}, the following result will be useful in section \ref{sec:Density}.

\begin{prop}\label{prop:square-snow-d-set}
The boundary of the square snowflake introduced above is a $d$-set for $d=\frac{3}{2}$.
\end{prop}

\begin{proof}
The structure of the proof is similar to that of Proposition \ref{prop:class-snow-d-set}, so we shall be briefer here.

Let $d=\frac32$ and $\xi=\frac{1}{4}$.

\emph{Step 1.} By similar reduction arguments as in Step 1 of the proof of Proposition \ref{prop:class-snow-d-set}, it is enough to prove that the part of the boundary {built around} the segment $[(0,0),(1,0)]$ in $\R^2$ is a $d$-set. To that effect, we intersect $\partial\Gamma^-$ and all sets involved in the definition of $\Gamma^-$ with the quarter plane $\Pi:=\{y\le x,\; y<1-x\} \cup \{(1,0)\}$.

\emph{Step 2.} We shall use the notation that has been used already in this subsection, except that we prepend the fraction $\frac{1}{4}$ to it. In particular, $\frac{1}{4}\Delta_j:=\Delta_j \cap \Pi$ and the part of the boundary to be considered is $\frac{1}{4}\partial \Gamma^- := \partial\Gamma^- \cap \Pi$, so that clearly $\frac{1}{4}\partial \Gamma^- = \bigcap_{j=0}^\infty \frac{1}{4}\Delta_j$.
From the way \eqref{Cond1} is proved in Proposition \ref{prop:SqSthick3propr}, we see that also the following holds:
$$\forall x\in {\textstyle \frac{1}{4}}\Delta_j, \;\; \dist(x,{\textstyle \frac{1}{4}}\partial\Gamma_j) < \xi^{j}\, \text{ and }\, \dist(x,{\textstyle \frac{1}{4}}\partial\Gamma^-)\leq \xi^{j}.$$
The first inequality comes from the fact that, in the square $S_{j,k}$ mentioned in the proof, the whole 
diagonal $E_{j,k}$ belongs to ${\textstyle \frac{1}{4}}\partial\Gamma_j$. The second inequality is due to the fact that the endpoints of $ E_{j,k} $ belong to ${\textstyle \frac{1}{4}}\Delta_{j'}$ for all $j'$, and therefore to ${\textstyle \frac{1}{4}}\partial\Gamma^-$. Then $d_H({\textstyle \frac{1}{4}}\partial\Gamma_j,{\textstyle \frac{1}{4}}\partial\Gamma^-)\leq \xi^j$, so that
\begin{equation}\label{eq:convsquaresnow}
\textstyle \frac{1}{4}\partial\Gamma_j \xrightarrow[j \to \infty]{} {\textstyle \frac{1}{4}}\partial\Gamma^- 
\end{equation}
in the Hausdorff metric in the space of non-empty compact subsets of $\R^2$.

\emph{Step 3.} We are going to show now that ${\textstyle \frac{1}{4}}\partial\Gamma^-$ is also the fractal determined by eight contractions $\psi_i$, $i=1,\ldots,8$, in $\R^2$ and that these contractions are indeed similarities with contraction ratio equal to $\xi$ and satisfy the open set condition. Then we can conclude, as in Step 3 of the proof of Proposition \ref{prop:class-snow-d-set}, that ${\textstyle \frac{1}{4}}\partial\Gamma^-$ is a compact $D$-set with $D\geq 0$ determined by $\sum_{i=1}^8 \xi^D =1$, from which it follows that $D=d$, finishing the proof.

The mentioned contractions are defined in the following way, using the notation from the proof of Proposition \ref{prop:class-snow-d-set}:
\begin{align*}
	\psi_1 &= H(\xi), &
	\psi_2 &= T(\xi,0)\circ R({\textstyle\frac\pi2}) \circ H(\xi),\\
	\psi_3 &= T(\xi,\xi) \circ H(\xi), &
	\psi_4 &= T(2\xi,\xi)\circ R({\textstyle-\frac\pi2}) \circ H(\xi), \\
	\psi_5 &= T(2\xi,0)\circ R({\textstyle-\frac\pi2}) \circ H(\xi), &
	\psi_6 &= T(2\xi,-\xi) \circ H(\xi),\\
	\psi_7 &= T(3\xi,-\xi)\circ R({\textstyle\frac\pi2}) \circ H(\xi),&
	\psi_8 &= T(3\xi,0) \circ H(\xi).
\end{align*}
\noindent These are, clearly, similarities of ratio $\xi$. They determine the unique non-empty compact set $K$ in $\R^2$ such that $K=\psi(K):=\bigcup_{i=1}^8 \psi_i(K)$ and which can be obtained as 
\eqref{eq:Kaslimit} (with $\psi$ as just defined)
for any non-empty compact $\Lambda\subset \R^2$, the limit being taken in the sense of the Hausdorff metric.

Choosing $\Lambda=[(0,0),(1,0)]$, it is easy to see from \eqref{eq:SqSrule} that $\psi^j(\Lambda)= \frac{1}{4}\partial\Gamma_j$, $j\in \N$. 
Combining this with \eqref{eq:convsquaresnow} and \eqref{eq:Kaslimit} (with $\psi$ as defined above), we get that $K=\frac{1}{4}\partial\Gamma^-$.
	
The proof finishes by observing that the similarities above satisfy the open set condition, for which we can take the interior of $\frac{1}{4}\Delta_0$  as the required open set.
\end{proof}

\subsection{Interior regular domains}
\label{subsec:IntRegDom}

Since in Corollary \ref{cor:Hs0} below we give a result concerning interior regular domains (see Definition \ref{def:intreg}) the boundary of which are $d$-sets with $0<d<n$ (see Definition \ref{def:dSet}), we would like to show here that all the snowflakes considered in this section \ref{sec:ExamplesThick} are examples of such domains. And since we have already proved in Propositions \ref{prop:class-snow-d-set} and \ref{prop:square-snow-d-set} that their boundaries are $d$-sets with $0<d<n$, it only remains to show that such domains are interior regular.

Actually, we are going to prove something more general, namely that any $I$-thick domain whose boundary satisfies the ball condition (see Definition \ref{def:ballcond}) is interior regular. This applies to our snowflakes because, on the one hand, we have already proved in Propositions \ref{prop:SNthick3propr} and \ref{prop:SqSthick3propr} that they are $I$-thick (and even thick) domains and, on the other hand, all $d$-sets with $0\leq d<n$ satisfy the ball condition (cf. \cite[Prop.~4.3]{Cae02} for the case $d\not=0$; the case $d=0$ is trivial).

\begin{defn}\label{def:intreg}
A domain $\Omega \subset \R^n$ is said to be \emph{interior regular}\footnote{This notion is taken from \cite[eq.~(4.95)]{Tri08}, except that we don't require $\Omega$ to be bounded nor to satisfy $(\overline{\Omega})^{\circ}=\Omega$. Note that if we replace the requirement of closedness in the Definition \ref{def:dSet} of $d$-set by openness, then being interior regular is equivalent to being an $n$-set.} if there is a constant $c>0$ such that $|\Omega \cap Q| \geq c\, |Q|$ for any cube $Q$ with side length at most $1$ centered at any point in $\partial\Omega$.
\end{defn}

\begin{defn}\label{def:ballcond}
(\!\!\cite[Def.~18.10]{Tri97})
A non-empty closed set $S$ in $\R^n$ is said to satisfy the \emph{ball condition} if there is $\eta\in(0,1)$ such that, for any ball $B(x,r)$ centered at $x\in S$ and with radius $r\in(0,1)$, there exists a ball $B(y,\eta r)$ such that
$$B(y,\eta r) \subset B(x,r) \quad \text{and} \quad B(y,\eta r) \cap S = \emptyset.$$
\end{defn}

\begin{rem}\label{rem:boundarydistance}
If necessary replacing $\eta$ by $\eta/2$, we can always assume that $S$ satisfying the ball condition also satisfies
$$\dist(B(y,\eta r), S) \geq \eta r.$$
\end{rem}
A set that satisfies the ball condition is also called \emph{porous}, e.g.\ in \cite[Def.~3.4, Rmk.~3.5]{Tri08}.

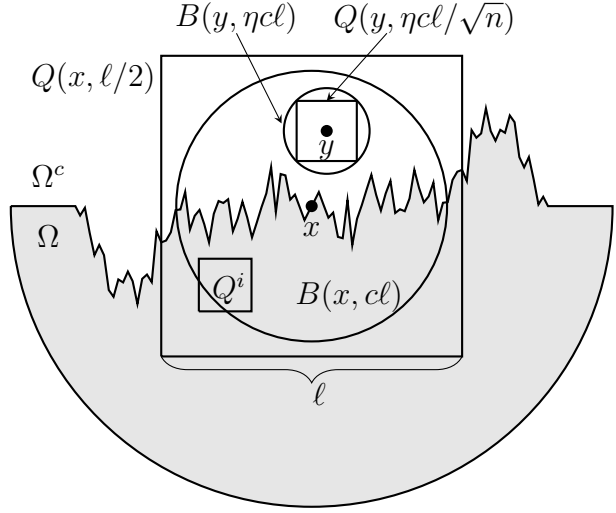
\begin{figure}[htb]\centering
\begin{tikzpicture}[scale=1]
\draw[thick,fill=black!10!white](-4,0)--
(-3.1416, 0.0000)--(-3.0781,-0.1243)--(-3.0147,-0.0221)--(-2.9512,-0.3560)--(-2.8877,-0.6527)--
(-2.8243,-0.5163)--(-2.7608,-0.6521)--(-2.6973,-1.1166)--(-2.6339,-1.2063)--(-2.5704,-0.9755)--
(-2.5069,-1.0152)--(-2.4435,-0.8705)--(-2.3800,-0.9987)--(-2.3165,-1.3010)--(-2.2531,-0.9673)--
(-2.1896,-0.9103)--(-2.1261,-1.0113)--(-2.0627,-0.3339)--(-1.9992,-0.4821)--(-1.9357,-0.4579)--
(-1.8723,-0.3479)--(-1.8088,-0.0532)--(-1.7453,-0.0447)--(-1.6819, 0.1892)--(-1.6184,-0.3129)--
(-1.5549,-0.1284)--(-1.4915,-0.0508)--(-1.4280,-0.2406)--(-1.3645,-0.4833)--(-1.3011,-0.1022)--
(-1.2376, 0.0677)--(-1.1741, 0.0417)--(-1.1107, 0.0915)--(-1.0472,-0.3350)--(-0.9837,-0.1681)--
(-0.9203, 0.0587)--(-0.8568,-0.1509)--(-0.7933,-0.4165)--(-0.7299,-0.2194)--(-0.6664,-0.0132)--
(-0.6029, 0.1691)--(-0.5395, 0.5182)--(-0.4760,-0.0796)--(-0.4125, 0.4270)--(-0.3491, 0.4092)--
(-0.2856, 0.0864)--(-0.2221, 0.1537)--(-0.1587,-0.0145)--(-0.0952,-0.1864)--(-0.0317,-0.0594)--
( 0.0317, 0.0594)--( 0.0952, 0.1864)--( 0.1587, 0.0145)--( 0.2221,-0.1537)--( 0.2856,-0.0864)--
( 0.3491,-0.4092)--( 0.4125,-0.4270)--( 0.4760, 0.0796)--( 0.5395,-0.5182)--( 0.6029,-0.1691)--
( 0.6664, 0.0132)--( 0.7299, 0.2194)--( 0.7933, 0.4165)--( 0.8568, 0.1509)--( 0.9203,-0.0587)--
( 0.9837, 0.1681)--( 1.0472, 0.3350)--( 1.1107,-0.0915)--( 1.1741,-0.0417)--( 1.2376,-0.0677)--
( 1.3011, 0.1022)--( 1.3645, 0.4833)--( 1.4280, 0.2406)--( 1.4915, 0.0508)--( 1.5549, 0.1284)--
( 1.6184, 0.3129)--( 1.6819,-0.1892)--( 1.7453, 0.0447)--( 1.8088, 0.0532)--( 1.8723, 0.3479)--
( 1.9357, 0.4579)--( 1.9992, 0.4821)--( 2.0627, 0.3339)--( 2.1261, 1.0113)--( 2.1896, 0.9103)--
( 2.2531, 0.9673)--( 2.3165, 1.3010)--( 2.3800, 0.9987)--( 2.4435, 0.8705)--( 2.5069, 1.0152)--
( 2.5704, 0.9755)--( 2.6339, 1.2063)--( 2.6973, 1.1166)--( 2.7608, 0.6521)--( 2.8243, 0.5163)--
( 2.8877, 0.6527)--( 2.9512, 0.3560)--( 3.0147, 0.0221)--( 3.0781, 0.1243)--( 3.1416,-0.0000)
--(4,0) arc (360:180:4);
\draw[thick] (-2,-2)--(2,-2)--(2,2)--(-2,2)--(-2,-2);\draw (-2.9,1.7)node{$Q(x,\ell/2)$};
\fill(0,0)circle(.08);   \draw (0,-.3)node{$x$};
\draw[thick](0,0)circle(1.8); \draw (.5,-1.2)node{$B(x,c\ell)$};
\fill(0.2,1)circle(.08);   \draw (0.2,.75)node{$y$};
\draw[thick](0.2,1)circle(.57); 
\draw (-1,2.5)node{$B(y,\eta c\ell)$};          \draw[-stealth](-1,2.3)--(-.4,1.1);
\draw (1.5,2.5)node{$Q(y,\eta c\ell/\sqrt n)$}; \draw[-stealth](1,2.3)--(.2,1.4);
\draw[thick] (-.2,.6)--(.6,.6)--(.6,1.4)--(-.2,1.4)--(-.2,.6);
\draw[thick] (-1.5,-1.4)--(-.8,-1.4)--(-.8,-.7)--(-1.5,-.7)--(-1.5,-1.4);\draw (-1.1,-1.1)node{$Q^i$};
\draw[decorate,decoration={brace,amplitude=10}](2,-2)--(-2,-2); \draw (0.1,-2.5)node{$\ell$};
\draw (-3.5,-.4)node{$\Omega$};
\draw (-3.5,.4)node{$\Omega^c$};
\end{tikzpicture}
\caption{A schematic representation of the proof of Proposition~\ref{prop:InteriorReg}. 
Here the case $B(y,\eta c \ell) \subset \overline{\Omega}^c$ is depicted.}
\label{fig:InteriorReg}
\end{figure}

\begin{prop}\label{prop:InteriorReg}
If $\Omega$ is an $I$-thick domain in $\R^n$ and $\partial \Omega$ satisfies the ball condition, then $\Omega$ is interior regular.
\end{prop}
\begin{proof}
From the $\eta \in (0,1)$ coming from the ball condition satisfied by $\partial\Omega$, fix $j_0=1$, $c_1=2\eta/\sqrt{n}$, $c_2=4\eta/\sqrt{n}$, $c_3=\eta$ and $c_4=2$ in the definition of $I$-thickness applied to $\Omega$. Consider then the constants $c_5$, $c_6$, $c_7$ and $c_8$ that come out from that definition and set $c = 1/(2(\sqrt{n}\,c_6+c_8+1))$. 

Let $x$ be any point of $\partial\Omega$ and let $Q(x,\ell/2)$ be a cube centered at $x$ with side length $\ell\leq 1$. Consider the ball $B(x,c \ell)$ and observe that $B(x,c\ell) \subset Q(x,\ell/2)$. By the ball condition satisfied by $\partial\Omega$ and Remark \ref{rem:boundarydistance} there exists $B(y,\eta c \ell)$ such that
$$
B(y,\eta c \ell) \subset B(x,c \ell) \quad \text{and} \quad \dist\big(B(y,\eta c \ell),\partial\Omega\big) \geq \eta c \ell.
$$

One of the following two situations must happen:
$$B(y,\eta c \ell)\subset \Omega \quad \text{or} \quad B(y,\eta c \ell) \subset \overline{\Omega}^c.$$
In the first case we have that
$$|Q(x,\ell/2)\cap \Omega| \geq |B(y,\eta c \ell)| \geq |Q(y,\eta c \ell /\sqrt{n})| = (2\eta c/\sqrt{n})^n |Q(x,\ell/2)|,$$
as required. In the second case, start by considering a cube $Q(y,\eta c \ell/\sqrt{n}) \in {\cal Q}(B(y,\eta c \ell)) \subset {\cal Q}(\Omega^c)$ and $j\in\N$ such that $2^{-j}<c \ell \leq 2^{-(j-1)}$ and observe that
$$c_{1}2^{-j}\leq l\big(Q(y,\eta c \ell/\sqrt{n})\big)\leq c_{2}2^{-j}
\quad\mbox{ and }\quad 
c_{3}2^{-j}\leq{\rm dist}\big(Q(y,\eta c \ell/\sqrt{n}),\partial\Omega\big)\leq c_{4}2^{-j}$$
for the constants $c_1$, $c_2$, $c_3$ and $c_4$ given above.
Then such $Q(y,\eta c \ell/\sqrt{n})$ is an exterior cube with respect to the $I$-thick domain $\Omega$, therefore there exists an interior cube $Q^i \in {\cal Q}(\Omega)$ such that
$$c_5 c \ell/2 \leq l(Q^{i})\leq c_6 c \ell\quad\mbox{ and }\quad 
c_7 c \ell/2 \leq{\rm dist}(Q^{i},\partial\Omega)\leq{\rm dist}\big(Q(y,\eta c \ell/\sqrt{n}),Q^{i}\big)
\leq c_8 c \ell.$$
For $z\in Q^i$ it holds by the choice of $c$ that
\begin{align*}
 |x-z| < \diam(Q^i) + \dist\big(Q(y,\eta c \ell / \sqrt{n}),Q^i\big) + \sup_{w\in Q(y,\eta c \ell / \sqrt{n})}|w-x|\leq (\sqrt{n}\,c_6 + c_8 + 1)c\ell=\frac\ell2,
\end{align*}
so that $Q^i \subset B(x,\ell/2) \subset Q(x,\ell/2)$. Hence
$$|Q(x,\ell/2)\cap \Omega| \geq |Q^i| \geq (c_5 c \ell/2)^n = (c_5 c/2)^n|Q(x,\ell/2)|,$$
which finishes the proof.
A sketch of the construction of the cubes and the balls involved in the proof is shown in Figure~\ref{fig:InteriorReg}.
\end{proof}
To summarise: all snowflakes introduced in this section, either classical or square, are interior regular, thick domains whose boundaries are compact $d$-sets satisfying the ball condition.

\section{Density of \texorpdfstring{$A_{p_1,q_1,\Gamma}^{s_1}$ in $A_{p_2,q_2,\Gamma}^{s_2}$ for $s_1>s_2$}{subscript spaces}}
\label{sec:Density}

In this section we give sufficient conditions under which $A_{p_1,q_1,\Gamma}^{s_1}$ is dense in $A_{p_2,q_2,\Gamma}^{s_2}$ for $\Gamma$ a closed set with empty interior. 
Recall that the spaces $A_{p,q,\Gamma}^{s}$ were defined in \eqref{eqn:SubscriptDef}, and that Proposition \ref{prop:Nullity} provides necessary and sufficient conditions on $s,p,q$ and $\Gamma$ for $A_{p,q,\Gamma}^{s}$ to be non-trivial. 
Our main focus is on the case where $\Gamma$ is a $d$-set for some $0<d<n$, which allows us to connect the spaces $A_{p,q,\Gamma}^{s}$ to certain trace spaces on $\Gamma$. 
We remind the reader that $d$-sets were defined in Definition \ref{def:dSet}.

Before tackling the density question for the $A^s_{p,q,\Gamma}$ spaces on $d$-sets with $0<d<n$, it is instructive to consider the limit case $d=0$, for which the spaces $A^s_{p,q,\Gamma}$ have a simple and explicit characterization. This allows us to give a rather complete answer to the density question, which provides a foretaste of the results obtained for the case $0<d<n$ later in the section. In particular, we note that, for $d=0$, $A^{s+1}_{p,q,\Gamma}$ is never dense in $A^{s}_{p,q,\Gamma}$ provided the latter is non-trivial.

\begin{prop}\label{prop:d0}
Let $\Gamma\subset\R^n$ be a non-empty compact $d$-set for $d=0$. Then $\Gamma$ is a finite set and, for all $0<p_1,p_2,q_1,q_2<\infty$ and $s_1,s_2\in \R$, with $\lfloor\cdot\rfloor$ denoting the integer part,
\begin{enumerate}[(i)]
\item if $s_1+n\left(1-\frac 1{p_1}\right)\ge0\text{ and }s_2+n\left(1-\frac 1{p_2}\right)\ge0$ then  $A^{s_1}_{p_1,q_1,\Gamma}=A^{s_2}_{p_2,q_2,\Gamma}=\{0\}$; 
\item if $\left\lfloor s_1+n\left(1-\frac 1{p_1}\right)\right\rfloor=\left\lfloor s_2+n\left(1-\frac 1{p_2}\right)\right\rfloor<0$ then $A^{s_1}_{p_1,q_1,\Gamma}=A^{s_2}_{p_2,q_2,\Gamma}\neq \{0\}$, with equivalent quasi-norms; 
\item if 
$\left\lfloor s_2+n\left(1-\frac 1{p_2}\right)\right\rfloor<0$ and $\left\lfloor s_1+n\left(1-\frac 1{p_1}\right)\right\rfloor>\left\lfloor s_2+n\left(1-\frac 1{p_2}\right)\right\rfloor$ then $A^{s_1}_{p_1,q_1,\Gamma}\subsetneqq A^{s_2}_{p_2,q_2,\Gamma}$ and the inclusion is \emph{not} dense. 
\end{enumerate}
\end{prop}
\begin{proof}
The fact that any compact $0$-set is finite follows trivially from the fact that $\cH^0$ is the counting measure. Without loss of generality it suffices to consider a set containing a single point, e.g.\ $\Gamma=\{0\}$, for which $A^s_{p,q,\Gamma}=A^s_{p,q,\{0\}}$, the subspace of $\Aspq$ of the elements supported at the origin. 
It is a standard result in distribution theory (see e.g.\ \cite[Thm.~3.9]{McLean}) that the only elements of $\cS'(\R^n)$ supported in $\{0\}$ are finite linear combinations of the delta function $\delta$ and its derivatives.
Reasoning for higher order derivatives of $\delta$ similarly as in \cite[Rmk~2.2.4.3]{RuS96} and afterwards applying  \cite[Prop.~2.3.2.2]{Tri83}, 
one can see that for any multi-index $\beta\in\N_0^n$, $D^\beta\delta\in A^{-n\left(1-\frac 1p\right)-\epsilon-|\beta|}_{p,q}(\R^n)\setminus A^{-n\left(1-\frac 1p\right)-|\beta|}_{p,q}(\R^n)$ for all $0<p,q<\infty$ and $\epsilon>0$.
So $A^s_{p,q,\{0\}}=\mathrm{span}\{D^\beta\delta, |\beta|<-s-n(1-\frac 1p)\}$, from which the basic statements of parts (i) and (ii) of the proposition follow immediately. The statement about equivalent quasi-norms in (ii) follows because $A^{s_1}_{p_1,q_1,\Gamma}=A^{s_2}_{p_2,q_2,\Gamma}$ is finite-dimensional. 
For part (iii), density fails because the two spaces $A^{s_1}_{p_1,q_1,\Gamma}$ and $A^{s_2}_{p_2,q_2,\Gamma}$ are finite dimensional with different dimension. 
(From an analytical perspective, this corresponds to the fact that it is not possible to approximate $D^\beta\delta$ with lower-order derivatives of $\delta$ centred at the same point.) 
\end{proof}

To study the case $0<d<n$ we need to consider traces on $d$-sets. The following proposition is a consequence of
\cite[\S18.5 and Cor.~18.12(i)]{Tri97} and the fact that any $d$-set
with $0<d<n$ satisfies the ball condition (recall Definition \ref{def:ballcond}) --- see, e.g., \cite[Prop.~4.3]{Cae02}. 
We mention also the important monograph \cite{JW84}, which contains many further results about traces on $d$-sets. Here, and henceforth, given $p\in(0,\infty)$, we denote by $L_{p}(\Gamma):=L_{p}(\Rn,\mathcal{H}^{d}|_{\Gamma})$ the complex $L_{p}$ space with respect to the restriction measure $\mathcal{H}^{d}|_{\Gamma}$ defined by $\mathcal{H}^{d}|_{\Gamma}(M):=\mathcal{H}^{d}(M\cap\Gamma)$ for all $\mathcal{H}^{d}$-measurable subsets $M$ of $\Rn$, 
equipped with the quasi-norm (norm if $p\geq 1$)
\[
\|f|L_{p}(\Gamma)\|:=\left(\int_{\Rn}|f(x)|^{p}\,\mathcal{H}^{d}|_{\Gamma}(dx)\right)^{1/p}=\left(\int_{\Gamma}|f(\gamma)|^{p}\,\mathcal{H}^{d}(d\gamma)\right)^{1/p},
\]
where the last identity holds because the support of
$\mathcal{H}^{d}|_{\Gamma}$ is exactly $\Gamma$. 
It is standard that $L_p(\Gamma)$ is a quasi-Banach space (Banach space if $p\geq1$). 

\begin{prop}
\label{prop:traceLp}Let $\Gamma$ be a $d$-set in $\Rn$ with $0<d<n$.
Let $0<p<\infty$. Then there exists $c>0$ such that
\[
\|\varphi|_{\Gamma}|L_{p}(\Gamma)\|\leq c\,\|\varphi|B_{p,\min\{1,p\}}^{\frac{n-d}{p}}(\Rn)\|\qquad\mbox{for all }\varphi\in\mathcal{S}(\Rn),
\]
and hence by completion there exists a unique continuous linear
operator $tr_{\Gamma}:\, B_{p,\min\{1,p\}}^{\frac{n-d}{p}}(\Rn)\to L_{p}(\Gamma)$
such that $tr_{\Gamma}f=f|_{\Gamma}$ whenever $f\in\mathcal{S}(\Rn)$.
Moreover, $tr_{\Gamma}$ is surjective
and there exist $c_{1},c_{2}>0$ such that
\[
c_{1}\|f|L_{p}(\Gamma)\|\leq\inf\,\|g|B_{p,\min\{1,p\}}^{\frac{n-d}{p}}(\Rn)\|\leq c_{2}\|f|L_{p}(\Gamma)\|\qquad\mbox{for all }f\in L_{p}(\Gamma),
\]
where the infimum runs, for each fixed $f\in L_{p}(\Gamma)$, over
all $g\in B_{p,\min\{1,p\}}^{\frac{n-d}{p}}(\Rn)$ such that $tr_{\Gamma}g=f$.
\end{prop}
Having defined $tr_\Gamma$ on $B_{p,\min\{1,p\}}^{\frac{n-d}{p}}(\Rn)$, for each $m\in\N_0$ we can define on $B_{p,\min\{1,p\}}^{\frac{n-d}{p}+m}(\Rn)$ the vector-valued trace operator
\begin{align*}
\Tr_{\Gamma,m} f:=(tr_{\Gamma} D^{\beta}f)_{0\leq |\beta|\leq m},
\end{align*}
consisting of the traces of all distinct partial distributional derivatives of order $\leq m$. This defines a continuous linear operator (not surjective in general) 
\[\Tr_{\Gamma,m}: B_{p,\min\{1,p\}}^{\frac{n-d}{p}+m}(\Rn) \to \big(L_p(\Gamma)\big)^{N_m},\]
where $N_m=\left|\{|\beta|\leq m\}\right|=\binom{n+m}{n}$ is the number of distinct partial distributional derivatives of order $\leq m$. When $m=0$, $N_0=1$ and $\Tr_{\Gamma,0}$ coincides with $tr_\Gamma$. 

For $s>\frac{n-d}{p}+m$ and $0<p,q<\infty$ we have the embedding
\begin{align}
\label{eqn:AspqBEmbed}
\Aspq \hookrightarrow B_{p,\min\{1,p\}}^{\frac{n-d}{p}+m}(\Rn),
\end{align}
so the restriction of $\Tr_{\Gamma,m}$ to such $\Aspq$ gives a continuous linear operator $\Tr_{\Gamma,m}|_{A_{p,q}^{s}}$ from $\Aspq$ into $(L_p(\Gamma))^{N_m}$. 
The range space $\Tr_{\Gamma,m}(\Aspq)$ is then linearly isomorphic to the quotient space $\Aspq/\ker(\Tr_{\Gamma,m}|_{A_{p,q}^{s}})$, which motivates the following:
\begin{defn}
\label{def:tracespace}
Let
$\Gamma$ be a $d$-set in $\Rn$ with $0<d<n$. Let $0<p,q<\infty$, $m\in\N_0$ 
and $s>\frac{n-d}{p}+m$. 
Define $\mathbb{A}_{p,q,m}^{s-\frac{n-d}{p}}(\Gamma)$ to be the vector space $\Tr_{\Gamma,m}(\Aspq)\subset (L_p(\Gamma))^{N_m}$ endowed with the
quasi-norm inherited from the quotient quasi-norm in $\Aspq/\ker(\Tr_{\Gamma,m}|_{A_{p,q}^{s}})$, i.e. for $(f^\beta)_{0\leq |\beta|\leq m} \in\mathbb{A}_{p,q,m}^{s-\frac{n-d}{p}}(\Gamma)$
\begin{align}
\|(f^\beta)_{0\leq |\beta|\leq m}|\mathbb{A}_{p,q,m}^{s-\frac{n-d}{p}}(\Gamma)\|:=\|[g]|\Aspq/\ker(\Tr_{\Gamma,m}|_{A_{p,q}^{s}})\|
=\inf_{g\in[g]}\,\|g|\Aspq\|,\label{eq:quotientqnorm}
\end{align}
where $[g]$ stands for the equivalence class containing all $g\in\Aspq$ such that $\Tr_{\Gamma,m}g=(f^\beta)_{0\leq |\beta|\leq m}$.\!
\end{defn}

Naturally, when $A=B$ or $A=F$, $\mathbb{A}$ should be replaced by $\mathbb{B}$ or $\mathbb{F}$ respectively, and $\IH^s_{p,m}(\Gamma)$ stands for $\mathbb{F}^s_{p,2,m}\GG$.

Since $\Aspq$ is complete and $\ker(\Tr_{\Gamma,m}|_{A_{p,q}^{s}})$ is
closed, $\mathbb{A}_{p,q,m}^{s-\frac{n-d}{p}}(\Gamma)$
is also complete (see \cite[Thm.~II.5.1]{TL80} for the case of normed spaces, but the proof can be adapted
to quasi-normed spaces).
Furthermore, by standard arguments it follows that the restricted operator 
\begin{equation}
\Tr_{\Gamma,m}|_{A_{p,q}^{s}}:\,\Aspq\longrightarrow\mathbb{A}_{p,q,m}^{s-\frac{n-d}{p}}(\Gamma)\label{eq:traceA}
\end{equation}
is continuous and surjective, and by the density of the embedding \rf{eqn:AspqBEmbed} (which follows because $\mathcal{S}(\Rn)$ is dense in both spaces),  
\begin{align}
\label{eqn:TraceDensity}
\mathbb{A}_{p,q,m}^{s-\frac{n-d}{p}}(\Gamma)\hookrightarrow \Tr_{\Gamma,m}\big(B_{p,\min\{1,p\}}^{\frac{n-d}{p}+m}(\Rn)\big)
\quad \text{with dense image in the $\big(L_p(\Gamma)\big)^{N_m}$ quasi-norm}.
\end{align}

In the following remark, and henceforth, the notation $\| f \| \lesssim \| g \|'$ indicates that there exists a constant $c>0$ (independent of $f,g$) such that $\| f \| \leq c\| g \|'$.

\begin{rem}\label{JW}
Here we show the connection with the trace spaces defined by Jonsson and Wallin in \cite{JW84}. Let $\Gamma$ be a $d$-set	in $\R^n$ with $0<d<n$. Let $1\leq p,q<\infty$ in the case of $B$ spaces and $1<p<\infty$, $q=2$ in the case of $F$ spaces. Let $m\in\N_0$ and $\frac{n-d}{p}+m<s\leq \frac{n-d}{p}+m+1$. For each $f\in L_1^{\rm loc}(\Rn)$ let $\bar{f}$ denote the strictly defined function given (a.e. in $\Rn$) by
$$\bar{f}(x):=\lim_{r\to 0} \frac{1}{|B(x,r) 
|}\int_{B(x,r)} f(y)\, dy.$$
It was proved by Jonsson and Wallin --- see the statements in \cite[Thms. VI.1 and VII.1, pp. 141 and 182]{JW84} --- that
$$T_{\Gamma,m}: f \longmapsto (\overline{D^\beta f}|_\Gamma)_{0\leq |\beta|\leq m}$$
establishes a continuous linear operator from $B^s_{p,q}(\Rn)$ onto a so-called Besov space $B^{s-\frac{n-d}{p}}_{p,q}(\Gamma)$, a subspace of $\big(L_p(\Gamma)\big)^{N_m}$, and from $H^s_p(\Rn)=F^s_{p,2}(\Rn)$ onto $B^{s-\frac{n-d}{p}}_{p,p}(\Gamma)$. And, moreover, that these operators admit bounded right inverses which are linear and acting in the same way as long as $s-\frac{n-d}{p}$ stays strictly between $m$ and $m+1$. Comparing with the way we have defined $\Tr_{\Gamma,m}$ and $\mathbb{A}^{s-\frac{n-d}{p}}_{p,q,m}(\Gamma)$, a density argument as in \cite[VIII.1.3, p. 211]{JW84} shows that $\Tr_{\Gamma,m}$ and $T_{\Gamma,m}$ act in the same way in the mentioned domains and that $\mathbb{B}^{s-\frac{n-d}{p}}_{p,q,m}(\Gamma)$ and $B^{s-\frac{n-d}{p}}_{p,q}(\Gamma)$, as well as $\mathbb{H}^{s-\frac{n-d}{p}}_{p,m}(\Gamma)$ and $B^{s-\frac{n-d}{p}}_{p,p}(\Gamma)$, coincide as sets. We claim that, besides coinciding as sets, in each of these pairs the norms are equivalent. 
Given any $g\in B^{s-\frac{n-d}{p}}_{p,q}(\Gamma)=\mathbb{B}^{s-\frac{n-d}{p}}_{p,q,m}(\Gamma)$ and any $f\in B^s_{p,q}(\Rn)$ such that $g=\Tr_{\Gamma,m} f = T_{\Gamma,m} f$, the continuity of $T_{\Gamma,m}|_{\Bspq}$ gives $\| g | B^{s-\frac{n-d}{p}}_{p,q}(\Gamma)\| \lesssim \| f | B^s_{p,q}(\Rn) \|$, so that, by \eqref{eq:quotientqnorm},
$$\| g | B^{s-\frac{n-d}{p}}_{p,q}(\Gamma)\| \lesssim \| g | \mathbb{B}^{s-\frac{n-d}{p}}_{p,q,m}(\Gamma)\|.$$
On the other hand, the existence of a bounded right inverse of $T_{\Gamma,m}|_{\Bspq}$ implies that there exists $h\in B^s_{p,q}(\Rn)$ such that $g=T_{\Gamma,m}h=\Tr_{\Gamma,m}h$ and $\|h|B^s_{p,q}(\Rn)\|\lesssim \| g |B^{s-\frac{n-d}{p}}_{p,q}(\Gamma)\|$, so that, using the continuity of \eqref{eq:traceA},
$$\| g | \mathbb{B}^{s-\frac{n-d}{p}}_{p,q,m}(\Gamma)\| \lesssim \| g | B^{s-\frac{n-d}{p}}_{p,q}(\Gamma)\|.$$
Similar arguments prove the equivalence of the $\mathbb{H}^{s-\frac{n-d}{p}}_{p,m}(\Gamma)$ and $B^{s-\frac{n-d}{p}}_{p,p}(\Gamma)$ norms.
\end{rem}

\begin{rem}\label{McLeantrace}
Under suitable regularity assumptions, the trace spaces of Definition~\ref{def:tracespace} coincide with classical trace spaces arising in PDE theory. 
For example, when $\Gamma$ is either the graph of a $C^{k-1,1}$ function $\zeta:\R^{n-1}\to\R$, $k\in\N$, or the boundary of a $C^{k-1,1}$ domain (both special cases of a $d$-set with $d=n-1$), 
McLean \cite[pp.~98--99]{McLean} defines the Hilbert space $H^s(\Gamma)\subset L^2(\Gamma)$ for $0\le s\le k$ as the push-forward of $H^s(\R^{n-1})$ under suitable coordinate charts. 
By \cite[Thm.~3.37]{McLean}, for $0<s\le k-\frac12$ the space $H^{s}\GG$ is the range of the classical trace operator $\gamma:H^{s+1/2}\OO\to L^2(\Gamma)$, 
which is defined by $\gamma u:=U|_\Gamma$ for $u=U|_\Omega\in \cD(\R^n)|_\Omega$, and by density for general $u\in H^{s+1/2}(\Omega)$, and $\gamma:H^{s+1/2}\OO\to H^s(\Gamma)$ has a bounded right inverse. 
Hence, for such $\Gamma$ and $0<s\le k-\frac12$, one can prove, using a similar argument to that employed in Remark \ref{JW}, that our space $\IH^s_{2,0}\GG$ and the space $H^s\GG$ defined in \cite{McLean} are linearly and topologically isomorphic.
\end{rem}

If we now restrict ourselves to the case $1< p,q<\infty$ then we can deduce by standard Banach space results the following density result, which will be important later. 

\begin{prop}
\label{prop:densedual}
Let $\Gamma$ be a $d$-set in $\Rn$ with $0<d<n$. 
Let $1<p,q<\infty$, $m\in\N_0$ and $s>\frac{n-d}{p}+m$. 
Then $\big(\Tr_{\Gamma,m}(B_{p,1}^{\frac{n-d}{p}+m}(\Rn))\big)'\hookrightarrow \big(\mathbb{A}_{p,q,m}^{s-\frac{n-d}{p}}(\Gamma)\big)'$ with dense image. 
\end{prop}
\begin{proof}
This is a consequence of the following general fact: if $X$ and $Y$ are Banach spaces such that $X$ is continuously embedded in $Y$ with dense image and $X$ is reflexive, then $Y'$ is densely embedded in $X'$. To prove this, let $E:X\to Y$ be the embedding of $X$ into $Y$. Since $E(X)$ is dense in $Y$, the ``if'' part of \cite[Thm.~3.1.17(b)]{megginson1998introduction} implies that the adjoint $E':Y'\to X'$ is injective. To show that $E'(Y')$ is dense in $X'$, by the ``only if'' part of \cite[Thm.~3.1.17(b)]{megginson1998introduction} it suffices to show that $E'':X''\to Y''$ (the adjoint of the adjoint) is injective. But, since $X$ is reflexive, $E''=\phi_Y\circ E \circ \phi_X^{-1}$ where $\phi_X:X\to X''$ and $\phi_Y:Y\to Y''$ are the canonical embeddings, so the required injectivity of $E''$ follows from that of $E$, $\phi_Y$ and $\phi_X^{-1}$.

To prove Proposition \ref{prop:densedual} we apply this result with $X=\mathbb{A}_{p,q,m}^{s-\frac{n-d}{p}}(\Gamma)$ and $Y=\Tr_{\Gamma,m}\big(B_{p,1}^{\frac{n-d}{p}+m}(\Rn)\big)$. The required dense embedding of $X$ in $Y$ is provided by \rf{eqn:TraceDensity} (specialised to the case $1<p<\infty$), and the reflexivity of $X$ follows from the reflexivity of $A_{p,q}^s(\Rn)$ (see \cite[\S2.6.1--2]{Tri78}), the surjectivity of \rf{eq:traceA}, and \cite[Cor.~1.11.22]{megginson1998introduction}.
\end{proof}

We now aim to establish a connection between the dual space $\big(\mathbb{A}_{p,q,m}^{s-\frac{n-d}{p}}(\Gamma)\big)'$ and the space of distributions $A^{-s}_{p',q',\Gamma}$. 
For this, we turn again to \rf{eq:traceA}, and note that since $\Tr_{\Gamma,m}|_{A_{p,q}^{s}}$ is surjective onto $\mathbb{A}_{p,q,m}^{s-\frac{n-d}{p}}(\Gamma)$ by the definition of this space, by \cite[Thm.~2.19, Rmk.~20]{Brezis} the adjoint operator 
\[ 
\big(\Tr_{\Gamma,m}|_{A_{p,q}^{s}}\big)': \big(\mathbb{A}_{p,q,m}^{s-\frac{n-d}{p}}(\Gamma)\big)' \to \big(\Aspq\big)' = I_{p',q'}^{-s,A}\big(A_{p',q'}^{-s}(\Rn)\big)
\]
is injective and a linear and topological isomorphism onto its range, which satisfies
\begin{equation}
\mathcal{R}\big((\Tr_{\Gamma,m}|_{A_{p,q}^{s}})'\big)=\big(\!\ker(\Tr_{\Gamma,m}|_{A_{p,q}^{s}})\big)^{a},\qquad s>\frac{n-d}{p}+m,\quad1 < p,q<\infty.\label{eq:rangeandkernel}
\end{equation}
The following proposition, which identifies $\ker(\Tr_{\Gamma,m}|_{A_{p,q}^{s}})$, is a generalisation of Triebel's \cite[Prop.~19.5]{Tri01}, which considered only the case $m=0$ and $\Gamma$ compact. 
Our arguments here, even in the case $m=0$, differ in some parts from Triebel's, since we consider that Triebel's proof does not provide enough evidence for the statement of \cite[Prop.~19.5]{Tri01}.
Specifically, it appears to presume that $\overline{D^\beta f}(x)=0$ $\mu$-a.e.\ on $\Gamma$ implies that $\overline{D^\beta f}(x)=0$ $(s-|\beta|,p)$-q.e.\ on $\Gamma$ for $0\leq |\beta| \leq m$ (cf.\ Eqns \eqref{capqe} and \eqref{0muae} below), which a priori is not obvious to us (though it comes as a consequence for the functions $f$ in the spaces $H^s_p(\Rn)$ or $B^s_{p,q}(\Rn)$ under the conditions of the proposition below after this has been proved). 

\begin{prop}
\label{prop:netrusov}Let $\Gamma$ be a $d$-set in $\Rn$
with $0<d<n$. Let $1<p<\infty$, $1\leq q<\infty$, $m\in\N_0$ and $\frac{n-d}{p}+m<s<\frac{n-d}{p}+m+1$.
Then $\mathcal{D}(\Gamma^{c})$ is dense in $\ker(\Tr_{\Gamma,m}|_{F_{p,2}^{s}})$
and in $\ker(\Tr_{\Gamma,m}|_{B_{p,q}^{s}})$.
That is, 
$$
\ker(\Tr_{\Gamma,m}|_{F_{p,2}^{s}})=\widetilde{F}_{p,2}^{s}(\Gamma^{c})=\widetilde{H}_{p}^{s}(\Gamma^{c}) 
\qquad\text{and}\qquad\ker(\Tr_{\Gamma,m}|_{B_{p,q}^{s}})=\widetilde{B}_{p,q}^{s}(\Gamma^{c}).
$$
\end{prop}
\begin{proof}
Since the inclusion $\supset$ is clear, we concentrate on proving the reverse one.

\emph{Step 1.} $\ker(\Tr_{\Gamma,m}|_{H^s_p})\cap C^m(\Rn)$ approximates $\ker(\Tr_{\Gamma,m}|_{H^s_p})$ in $H^s_p(\Rn)$.

Let $f\in \ker(\Tr_{\Gamma,m}|_{H^s_p})$. Then there exists $(f_j)_{j\in\N} \subset {\cal D}(\R^n)$ such that
\begin{equation}\label{DconvH}
f_j \xrightarrow[j \to \infty]{} f \quad \mbox{in }\, H^s_p(\Rn),
\end{equation}
and, by the continuity of $\Tr_{\Gamma,m}|_{H^s_p}$ (Eqn.~\eqref{eq:traceA}),
\begin{equation}\label{TDconvH}
\Tr_{\Gamma,m}|_{H^s_p} f_j \xrightarrow[j \to \infty]{} \Tr_{\Gamma,m}|_{H^s_p} f = 0 \quad \mbox{in }\, \mathbb{H}^{s-\frac{n-d}{p}}_{p,m}(\Gamma).
\end{equation}
By Remark \ref{JW}, and denoting by ${\cal E}_{\Gamma,m}$ the appropriate bounded right inverse mentioned there, 
\begin{equation}\label{g_j}
g_j := {\cal E}_{\Gamma,m}(\Tr_{\Gamma,m}|_{H^s_p} f_j)\; \in\,H^s_p(\Rn)
\end{equation}
is such that $\|g_j |H^s_p(\Rn) \| \lesssim \|\Tr_{\Gamma,m}|_{H^s_p} f_j | \mathbb{H}^{s-\frac{n-d}{p}}_{p,m}(\Gamma)\|$ and
\begin{equation}\label{Trg=Trf}
\Tr_{\Gamma,m} g_j = \Tr_{\Gamma,m} f_j,
\end{equation}
and by \eqref{TDconvH} it follows that
\begin{equation}\label{gconv0}
g_j \xrightarrow[j \to \infty]{} 0 \quad \mbox{in }\, H^s_p(\Rn).
\end{equation}
Define now, for any $j\in\N$, $h_j:=f_j-g_j \in H^s_p(\Rn)$. We have, by \eqref{DconvH}, \eqref{gconv0} and \eqref{Trg=Trf}, that
$$h_j \xrightarrow[j \to \infty]{} f \quad \mbox{in }\, H^s_p(\Rn) \qquad  \mbox{and} \qquad h_j \in \ker(\Tr_{\Gamma,m}|_{H^s_p}).$$
So, our claim above will be proved if we show that $h_j\in C^m(\Rn)$.
Clearly, it is enough to prove that this is the case for the functions $g_j$. And we shall prove this by combining the definition \eqref{g_j} of $g_j$ with the properties of ${\cal E}_{\Gamma,m}$.

We recall, from the discussion in Remark \ref{JW}, that ${\cal E}_{\Gamma,m}$ comes from \cite{JW84}; more precisely it is the operator named $\cal E$ in \cite[Thm.~VII.3, p.~197]{JW84}, taking $k$ there equal to our $m$ here. As is explicitly mentioned in \cite[p.~197]{JW84}, it is the same operator as that considered in \cite[Thm.~VI.3, p.~155]{JW84}, and it is defined in \cite[pp.~156--157]{JW84}.
For $s-\frac{n-d}{p} \in (m,m+1)$, it acts in exactly the same way.
Since $f_j \in {\cal D}(\Rn) \subset B^{s-\frac{n-d}{p}}_{\infty,\infty}(\Rn)$, with $s-\frac{n-d}{p}-\frac{n-d}{\infty} \in (m,m+1)$, then, by \cite[Thms.~VI.1 and VI.3, pp.~141 and 155]{JW84}, $\Tr_{\Gamma,m}|_{H^s_p} f_j = ( D^\beta f_j|_\Gamma )_{0\leq |\beta|\leq m} = T_{\Gamma,m}|_{B^{s-\frac{n-d}{p}}_{\infty,\infty}} f_j \in B^{s-\frac{n-d}{p}}_{\infty,\infty}(\Gamma)$ and $g_j \in B^{s-\frac{n-d}{p}}_{\infty,\infty}(\Rn) = {\cal C}^{s-\frac{n-d}{p}}(\Rn)$, the last identity coming from \cite[p.~8]{JW84} or \cite[(2.3.5.1), p.~51]{Tri83}, the space ${\cal C}^{s-\frac{n-d}{p}}(\Rn)$ being called a Lipschitz type space in \cite[p.~2]{JW84} or a Zygmund space in \cite[p.~36]{Tri83}. In any case, what matters for us is that the elements of ${\cal C}^{s-\frac{n-d}{p}}(\Rn)$ belong to $C^m(\Rn)$.

\emph{Step 2.}  $\cD(\Gamma^c)$ approximates $\ker(\Tr_{\Gamma,m}|_{H^s_p})\cap C^m(\Rn)$ in $H^s_p(\Rn)$.

Let $f \in \ker(\Tr_{\Gamma,m}|_{H^s_p})\cap C^m(\Rn)$. We shall use Netrusov's theorem \cite[Thm.~10.1.1, p.~281]{AdHe}, which, in particular, states that the desired approximability by elements of ${\cal D}(\Gamma^c)$ holds provided
\begin{equation}\label{capqe}
\overline{D^\beta f}(x)=0 \quad (s-|\beta|,p)\mbox{-q.e. on } \Gamma,\quad 0\leq |\beta|<s,
\end{equation}
where $(\sigma,p)$-q.e.\ means up to a set of zero capacity $C_{\sigma,p}$ (cf.~\cite[Def.~2.2.6, p.~20]{AdHe} for the definition). Recall that the bar over the function stands for the corresponding strictly defined function as in Remark \ref{JW} above. 
From \cite[Thm.~5.1.9]{AdHe} it follows that 
\begin{equation}\label{0cap}
C_{\sigma,p}(\Gamma)=0 \quad \mbox{when }\; 1<p<\infty \; \mbox{ and }\; 0< \sigma < \frac{n-d}{p}.
\end{equation}
Applying this to $\sigma=s-|\beta|$ for $m+1\leq |\beta|<s$, we have that $C_{s-|\beta|,p}(\Gamma)=0$ when $m+1\leq |\beta|<s$, so that \eqref{capqe} holds trivially for these values of $\beta$. For the remaining values of $\beta$, the assumption that $f \in \ker(\Tr_{\Gamma,m}|_{H^s_p})$ implies that 
\begin{equation}\label{0muae}
D^\beta f(x)=\overline{D^\beta f}(x)=0 \quad \cH^d\mbox{-a.e. on } \Gamma, \quad 0\leq |\beta| \leq m
\end{equation}
(cf. also the discussion in Remark \ref{JW} above). We claim that from this and the continuity of $D^\beta f$ it follows that
$$D^\beta f(x)=0 \quad \mbox{everywhere on } \Gamma, \quad 0\leq |\beta| \leq m,$$
which proves \eqref{capqe} for the remaining values of $\beta$. To prove this claim, observe that if there were a point $x\in\Gamma$ with $D^\beta f(x) \not= 0$, then there would exist $r>0$ such that $D^\beta f(y)=0$ for all $y \in B(x,r)$. 
But since $0 < r^d \lesssim \mu(\Gamma \cap B(x,r))$, we would then have a contradiction with \eqref{0muae}.

This finishes the proof of the proposition for the $H^s_p$ spaces.

\emph{Step 3.} Extension to $B^s_{p,q}$ spaces.

The result can be extended to the $B_{p,q}^{s}$ spaces by interpolation from $H^s_p$ as in Triebel's proof in \cite[pp.~262--263]{Tri01}, using here, for each $m$, the $s$-independent (for $s \in (\frac{n-d}{p}+m,\frac{n-d}{p}+m+1)$) bounded linear right inverses mentioned in Remark \ref{JW} above. 
\end{proof}

\begin{rem}
The technique used in the above proof, of reduction to functions with enough regularity, is taken from \cite[Step 1 of proof of Prop.~3.5]{FaJa99}, with a reference to \cite[(i) of proof of Thm.~1]{marschall1987trace}. 
In both \cite{FaJa99} and \cite{marschall1987trace}, for the crucial part corresponding to verifying that $g_j$ has the appropriate regularity, the reader is invited to check a related proof in \cite{Ste70}. By contrast, in our proof above we give a complete (and short) justification of this step using a result readily available in \cite{JW84}. In the case $m=0$ a more direct proof can also be seen in \cite[proof of Prop.~2.26]{CC08}.
\end{rem}

\begin{rem}
\label{rem:LimitCase1}
It is an open problem whether the result of Proposition \ref{prop:netrusov} holds also when $s = \frac{n-d}{p}+m+1$. If this were the case, we could update several of the results which follow by the inclusion of a corresponding limiting situation. See Remark \ref{rem:LimitCase2}. 
\end{rem}

From Proposition \ref{prop:netrusov} it is straightforward to prove the following result --- which may be of independent interest --- regarding a characterization of $H^s_{p,0}(\Omega)$, the closure of $C_0^\infty(\Omega)$ in $H^s_p(\Omega)$, in terms of the kernel of a trace operator. 

\begin{cor}[Traces from $\Omega\ne\Rn$ to $\partial\Omega$]
\label{cor:Hs0}
Let $\Gamma$, $d$, $p$, $m$ and $s$ be as in Proposition \ref{prop:netrusov}, and suppose further that $\Gamma=\partial\Omega$ with $\Omega$ an \emph{interior regular} domain (recall Definition \ref{def:intreg}). Then
$$H^s_{p,0}(\Omega)=\ker(\Tr^\Omega_{\Gamma,m}),$$
where, given any $f\in H^s_p(\Omega)$, 
\begin{equation}\label{traceOmega}
\Tr^\Omega_{\Gamma,m}f := \Tr_{\Gamma,m}u \; \mbox{ for any } u\in H^s_p(\Rn) \mbox{ such that } u|_\Omega=f,
\end{equation}
this definition being independent of the choice of $u$ due to the assumption of interior regularity \cite[Prop.~3.3]{Cae00}.
\end{cor}

\begin{rem} Corollary \ref{cor:Hs0} extends, e.g., \cite[Thm.~3.40]{McLean} (which considers only $p=2$, and requires that $\Omega$ is of class $C^{\lceil s\rceil-1,1}$, and in the Lipschitz case that $s$ is restricted to 
$\frac{n-(n-1)}{2}+0 = \frac{1}{2} < s \leq 1 < \frac{n-(n-1)}{2}+0+1$) and \cite[Thm.~3.5]{FaJa99} (which assumes $\Omega$ is a bounded $(\varepsilon,\delta)$ domain, and considers only $m=0$, with $s$ restricted to $\frac{n-d}{p}+0 = \frac{n-d}{p} < s \leq 1 < \frac{n-d}{p}+0+1$). This is because the traces in these two references also satisfy \eqref{traceOmega} (cf.\ Remarks \ref{McLeantrace} and \ref{JW} and \cite[Thm.~3.2]{FaJa99}), McLean's Lipschitz domains are assumed to have compact boundaries, so clearly are interior regular domains, and $(\varepsilon,\delta)$ domains are special cases of interior regular domains (cf.\ \cite[Prop.~1, p.~119]{Wall91}).
Similar results hold for $B^s_{p,q}$ spaces fitting the hypotheses of Proposition \ref{prop:netrusov}, thus also extending corresponding results known in more restricted settings (cf., e.g., \cite[Thm.~1, p.~49]{marschall1987trace}). See also a related result in \cite[Thm. 3]{Wall91}. We recall that all the snowflakes considered in \S\ref{sec:ExamplesThick} are examples of (bounded) interior regular domains whose boundaries are $d$-sets with $0<d<n$ --- cf. \S\ref{subsec:IntRegDom}.
\end{rem}

\begin{rem}
\label{rem:notation}We shall use the shortcut $A_{p,(q)}^{s}(\Rn)$ to deal simultaneously with both $H_{p}^{s}(\Rn)$ and $\Bspq$ in the above context, and adapt similarly the other notation.
\end{rem}

We can now make the connection with the spaces we want to consider. Combining Proposition~\ref{prop:netrusov} with (\ref{eq:rangeandkernel}) and Proposition \ref{prop:annihilators} reveals that
\[\mathcal{R}\big((\Tr_{\Gamma,m}|_{A_{p,(q)}^{s}})'\big)=I_{p',(q')}^{-s,A}(A_{p',(q'),\Gamma}^{-s}),\]
and this completes the proof of one of the major results of this section. 
\begin{thm}
\label{thm:gammasubgamma}Let $\Gamma$ be a $d$-set in $\Rn$
with $0<d<n$. Let $1<p,q<\infty$, $m\in\N_0$ and $\frac{n-d}{p}+m<s<\frac{n-d}{p}+m+1$.
Then, with the notation set above, the operator
\begin{equation}\label{eq:AdjTrMapping}
(I_{p',(q')}^{-s,A})^{-1}\circ(\Tr_{\Gamma,m}|_{A_{p,(q)}^{s}})'
\quad : \quad \big(\mathbb{A}_{p,(q),m}^{s-\frac{n-d}{p}}(\Gamma)\big)'
\quad \to \quad A_{p',(q'),\Gamma}^{-s}
\end{equation}
is a linear and topological isomorphism.
\end{thm}

The significance of Theorem~\ref{thm:gammasubgamma} is that it identifies, via the adjoint of the restricted trace operator $\Tr_{\Gamma,m}|_{A_{p,(q)}^{s}}$, a space of distributions $A_{p',(q'),\Gamma}^{-s}$ defined on $\R^n$ and supported in $\Gamma$, with the dual $\big(\mathbb{A}_{p,(q),m}^{s-\frac{n-d}{p}}(\Gamma)\big)'$ of a trace space on $\Gamma$. (The nature of this identification is discussed further in Remark \ref{rem:idGamma} below.) As a result, we can 
deduce density results for the $A_{p',(q'),\Gamma}^{-s}$ spaces from the corresponding density results for the $\big(\mathbb{A}_{p,(q),m}^{s-\frac{n-d}{p}}(\Gamma)\big)'$ spaces. Indeed, by combining Theorem \ref{thm:gammasubgamma} with Proposition \ref{prop:densedual}, we obtain another of our main results. Note that in the following theorem we have switched $(s,p,q)\leftrightarrow(-s,p',q')$ compared to Theorem \ref{thm:gammasubgamma}, to put the focus on the $A_{p,(q),\Gamma}^{s}$ spaces.
We highlight that we have given the positive answer to Q2 promised at the beginning of the section.

\mypropbox{
\begin{thm}\label{thm:maindensity2}
Let $\Gamma$ be a $d$-set in $\Rn$ with $0<d<n$. 
Let 
$$1<p,q_1,q_2<\infty,\quad m\in\N_0\qquad\text{and}\qquad-\frac{n-d}{p'}-m-1<s_{2}  \leq  s_{1}<-\frac{n-d}{p'}-m.$$
Then
$$B_{p,q_1,\Gamma}^{s_{1}}\;\text{is dense in}\; B_{p,q_2,\Gamma}^{s_{2}}
\qquad\text{and}\qquad H_{p,\Gamma}^{s_{1}}\;\text{is dense in}\;
H_{p,\Gamma}^{s_{2}}.$$
\end{thm}}
\begin{proof}
First recall that, by \cite[Prop.~2.3.2.2]{Tri83} and \eqref{eqn:SubscriptDef}, 
$B_{p,q_1,\Gamma}^{s_{1}}\subset B_{p,q_2,\Gamma}^{s_{2}}$ and $H_{p,\Gamma}^{s_{1}}\subset H_{p,\Gamma}^{s_{2}}$.
The density assertion follows from Theorem~\ref{thm:gammasubgamma}, combined with the fact that $\big(\mathbb{B}_{p',q_1',m}^{-s_1-\frac{n-d}{p'}}(\Gamma)\big)'\hookrightarrow \big(\mathbb{B}_{p',q_2',m}^{-s_2-\frac{n-d}{p'}}(\Gamma)\big)'$ 
and
$\big(\mathbb{H}_{p',m}^{-s_1-\frac{n-d}{p'}}(\Gamma)\big)'\hookrightarrow \big(\mathbb{H}_{p',m}^{-s_2-\frac{n-d}{p'}}(\Gamma)\big)'$ 
with dense image, which holds because by Proposition~\ref{prop:densedual} all these spaces contain $\big(\Tr_{\Gamma,m}\big(B_{p',1}^{\frac{n-d}{p'}+m}(\Rn)\big)\big)'$ as a dense subspace. 
We note that all the embeddings and identifications involved are compatible, e.g.\ 
$I_{p,q_1}^{s_1,B}=I_{p,q_2}^{s_2,B}|_{B^{s_1}_{p,q_1}}$, $\Tr_{\Gamma,m}|_{B_{p',q_2'}^{-s_2}}=(\Tr_{\Gamma,m}|_{B_{p',q_1'}^{-s_1}})|_{B_{p',q_2'}^{-s_2}}$ and similarly for the $H^s_p$ spaces.
The basic structure underlying the proof is summarised in Figure~\ref{fig:Jungle}.
\begin{figure}[htb]
\begin{center}
\hrule
\begin{tikzpicture}
\matrix[matrix of math nodes, row sep={60pt,between origins},text height=1.5ex, text depth=0.25ex] (s)
{|[name=X2]| A^s_{p,q}(\R^n) & \subdense &
|[name=X1]| B^{\frac{n-d}p}_{p,1}(\R^n)\\
|[name=Y2]| \mathbb A^{s-\frac{n-d}p}_{p,q}\GG & \subdense &
|[name=Y1]| L_p(\Gamma) \\};
\draw[->>] 
 (X2) edge node[auto]{$tr_\Gamma|_{A^s_{p,q}}$} (Y2)
 (X1) edge node[auto]{$tr_\Gamma$} (Y1) ;
\end{tikzpicture}\vrule
\begin{tikzpicture}
\matrix[matrix of math nodes, row sep={60pt,between origins},text height=1.5ex, text depth=0.25ex] (s)
{|[name=V1]| B^{-\frac{n-d}p,\Gamma}_{p',\infty}& \subcl &
|[name=H1]| B^{-\frac{n-d}p}_{p',\infty}(\R^n) & \subdense &
|[name=H2]| A^{-s}_{p',q'}(\R^n) & \supcl &
|[name=V2]| A^{-s}_{p',q',\Gamma} \\
&&|[name=X1s]| \big(B^{\frac{n-d}p}_{p,1}(\R^n)\big)' & \subdense &
  |[name=X2s]| \big(A^s_{p,q}(\R^n)\big)'\\
&&|[name=Y1s]| L_{p'}\GG & \subdense &
  |[name=Y2s]| \big(\mathbb A^{s-\frac{n-d}p}_{p,q}\GG\big)'\\};
\draw[right hook->]
 (Y2s) edge[bend left=0] node[auto,pos=0.5]{$tr_\Gamma|_{A^s_{p,q}}'$} (X2s)
 (Y1s) edge[bend left=0] node[auto,pos=0.5]{$tr_\Gamma'$} (X1s) 
 (V1) edge[bend left=30] node[auto]{$\subdense$} (V2) ;
\draw[right hook->>]
 (H1) edge node[auto]{$I^{-\frac{n-d}p,B}_{p',\infty}$} (X1s)
 (H2) edge node[auto]{$I^{-s,A}_{p',q'}$} (X2s)
 (Y2s) edge[bend right=30] (V2)
 (Y1s) edge[bend left=30] node[auto,pos=.0,rotate=-45]{\small$(I^{-\frac{n-d}p,B}_{p',\infty})^{-1}\circ tr_\Gamma'$} (V1) ;
\end{tikzpicture}\hrule
\end{center}
\caption{Diagram summarising the main relations between the function spaces described in \S\ref{sec:Density}, in the case $m=0$, $\frac{n-d}{p}<s<\frac{n-d}{p}+1$ and $1<p,q<\infty$.
Here $\subdense$ and $\subcl$ denote embeddings with dense and closed ranges, respectively, and $\hookrightarrow$ and $\twoheadrightarrow$ denote injective and surjective bounded mappings, respectively. The space $B^{-\frac{n-d}p,\Gamma}_{p',\infty}$ is defined in Remark~\ref{rem:TraceKerR}.
The density of the embeddings for the spaces on $\R^n$ (i.e.\ $A^s_{p,q}(\R^n) \subdense B^{\frac{n-d}p}_{p,1}(\R^n)$) implies the same for the trace spaces ($\mathbb A^{s-\frac{n-d}p}_{p,q}\GG \subdense L_p(\Gamma)$), their duals
($L_{p'}\GG \subdense (\mathbb A^{s-\frac{n-d}p}_{p,q}\GG)'$) and their identifications 
($B^{-\frac{n-d}p,\Gamma}_{p',\infty}\subdense A^{-s}_{p',q',\Gamma}$) via $tr_\Gamma'$.
Then Theorem~\ref{thm:maindensity2} follows from the fact the latter density holds for different values of $s$ (and $q$).}
\label{fig:Jungle}
\end{figure}
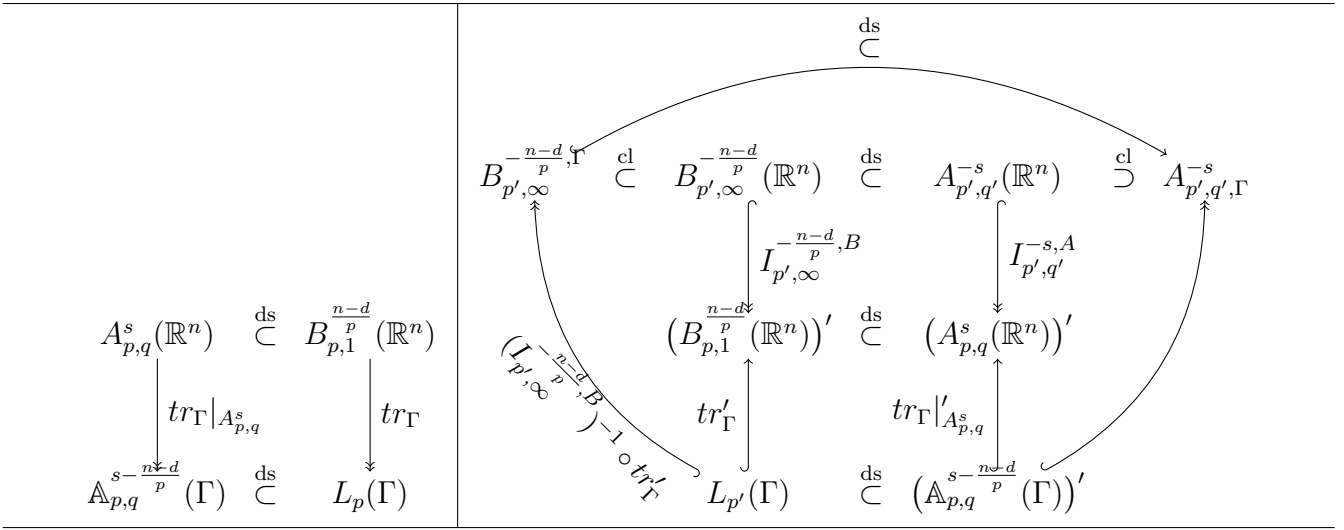
\end{proof}

\begin{rem}[Limiting case]
\label{rem:LimitCase2}
If the result of Proposition \ref{prop:netrusov} could be extended to include the limiting case $s=(n-d)/p + m+1$ (which seems to be an open problem --- see Remark \ref{rem:LimitCase1}), then the density result of Theorem \ref{thm:maindensity2} would extend to the limiting case $s_2=-(n-d)/p'-m-1$, as holds for the case of an $(n-1)$-dimensional hyperplane and the $H^s$ spaces, as discussed in \S\ref{sec:Intro}.
\end{rem}

\begin{rem}[The adjoint of the trace operator is the identification operator]
\label{rem:idGamma} \hspace{.1mm}
To give a more concrete description of the identification of $A_{p',(q'),\Gamma}^{-s}$ and $(\mathbb{A}_{p,(q),m}^{s-\frac{n-d}{p}}(\Gamma))'$
provided by Theorem \ref{thm:gammasubgamma}, we point out that, as discussed by Triebel in \cite[\S9.2]{Tri01} (in the case $A=F$), for $m=0$ the adjoint operator $(\Tr_{\Gamma,0}|_{A_{p,q}^{s}})'=(tr_{\Gamma}|_{A_{p,q}^{s}})':\big(\mathbb{A}^{s-(n-d)/p}_{p,q,0}(\Gamma)\big)'\to (\Aspq)' = I^{-s}_{p',q'}(\Aspqp)$ appearing in Theorem \ref{thm:gammasubgamma} can be viewed as an extension (by density) of the standard \emph{identification operator} $id_{\Gamma}$ identifying $L_{p'}$ functions on $\Gamma$ with tempered distributions on $\Rn$ (see e.g.\ \cite[\S9.2]{Tri01} and \cite[Eqn.~(18.6)]{Tri97}).

In more detail, it is well-known that the dual space of $L_p(\Gamma)$ can be realised as $L_{p'}(\Gamma)$ using the identification $I_{p'}:L_{p'}(\Gamma)\to (L_{p}(\Gamma))'$ defined for $f\in L_{p'}(\Gamma)$ and $g\in L_{p}(\Gamma)$ by $(I_{p'}f)(g):=\int_{\Gamma}f(\gamma)g(\gamma)\,\mathcal{H}^{d}(d\gamma)$. 
Also, by \cite[Thm.~2.11.2]{Tri83} the duality result in Proposition \ref{prop:dualA} extends to the case $A=B$, $1<p<\infty$ and $q=1$ ($q'=\infty$), giving an isomorphism
\[I^{-s,B}_{p',\infty}:B^{-s}_{p',\infty}(\Rn) \to \left(B^s_{p,1}(\Rn)\right)', \]
which extends by density the action of tempered distributions on elements of $\mathcal{S}(\Rn)$. 
Recalling that $tr_{\Gamma}:\, B_{p,1}^{\frac{n-d}{p}}(\Rn)\to L_{p}(\Gamma)$ is surjective, the adjoint 
\[tr_\Gamma'\;:\;\big(L_p(\Gamma)\big)'
=I_{p'}\big(L_{p'}(\Gamma)\big)
\;\to \;\big(B_{p,1}^{\frac{n-d}{p}}(\Rn)\big)' 
= I^{-\frac{n-d}{p},B}_{p',\infty}\big(B^{-\frac{n-d}{p}}_{p',\infty}{(\R^n)}\big)\]
is an isomorphism onto its image, and acts by
\[ tr_\Gamma'(l)(g) = l(tr_\Gamma g) = \int_{\Gamma}f_l(\gamma)(tr_\Gamma g)(\gamma)\,\mathcal{H}^{d}(d\gamma), \qquad l\in \big(L_p(\Gamma)\big)', g\in B_{p,1}^{\frac{n-d}{p}}(\Rn),\]
where $f_l=(I_{p'})^{-1}(l)\in L_{p'}(\Gamma)$. 
In particular, taking $g\in\mathcal{S}(\Rn)$ and replacing $tr_\Gamma g$ by $g|_\Gamma$ (as per the definition of $tr_\Gamma$) we recover Triebel's identification operator (cf.\ \cite[\S9.2]{Tri01} and \cite[Eqn.~(18.6)]{Tri97})
\[ (id_\Gamma f)(g)= tr_\Gamma'\big(I_{p'}(f)\big)(g) = \int_{\Gamma}f(\gamma)(g|_\Gamma)(\gamma)\,\mathcal{H}^{d}(d\gamma), \qquad f\in L_{p'}(\Gamma), g\in \mathcal{S}(\Rn).\]
\end{rem}
\begin{rem}[The kernel of the trace operator and the range of its adjoint]\label{rem:TraceKerR}
In \cite{Tri97,Tri01}, Triebel describes in some detail the mapping properties of $tr_\Gamma$ and $tr_\Gamma'$. 
In particular, in \cite[Thm.~18.2]{Tri97} it is proved that the range of $tr_\Gamma'$ satisfies 
\[\mathcal{R}(tr_\Gamma')=I^{-\frac{n-d}{p},B}_{p',\infty}(B_{p',\infty}^{-\frac{n-d}{p},\Gamma}),\] 
where
\[
B_{p',\infty}^{-\frac{n-d}{p},\Gamma}:=\big\{g\in B_{p',\infty}^{-\frac{n-d}{p}}(\Rn):\,\left\langle g,\varphi\right\rangle =0\;\mbox{ for all }\varphi\in\mathcal{S}(\Rn)\mbox{ such that }\varphi|_{\Gamma}=0\big\}.
\]
Obviously we have the inclusion
\begin{align}\label{eqn:SupSub}
 B_{p',\infty}^{-\frac{n-d}{p},\Gamma} \subset B_{p',\infty,\Gamma}^{-\frac{n-d}{p}}.
\end{align}
and, since 
\[ \big(\widetilde{B}^{\frac{n-d}p}_{p,1}(\Gamma^c)\big)^a 
= I^{-\frac{n-d}p,B}_{p',\infty}(B^{-{\frac{n-d}p}}_{p',\infty,\Gamma})\]
(as is easily proved by the same argument used to prove Proposition \ref{prop:annihilators}), 
we recover another obvious inclusion:
\begin{align} \label{eqn:kertr}
\ker{tr_\Gamma} = 
{}^a\big(\mathcal R (tr_\Gamma')\big)=
{}^a\Big(I^{-\frac{n-d}{p},B}_{p',\infty}\big(B_{p',\infty}^{-\frac{n-d}{p},\Gamma}\big) \Big)
\supset{}^a\Big(I^{-\frac{n-d}{p},B}_{p',\infty}\big(B_{p',\infty,\Gamma}^{-\frac{n-d}{p}}\big) \Big)
= \widetilde{B}^{\frac{n-d}{p}}_{p,1}(\Gamma^c).  
\end{align}
If $d>n-1$ then we have equality in \rf{eqn:SupSub} (and hence in \rf{eqn:kertr}) --- see \cite[\S9.34(viii)]{Tri01}. But for $0<d\leq n-1$ the inclusions in \rf{eqn:SupSub} and \rf{eqn:kertr} may be strict (cf.\ the discussion in \cite[\S17.3, p.~126]{Tri97}). 
\end{rem}

The following is a simple corollary of Proposition \ref{prop:netrusov}, relating to properties of the $\widetilde{A}^{s_2}_{p,(q)}$ spaces. 
\begin{cor}\label{cor:TildeCap}
Let $\Gamma$, $p$, $q$, and $m$ be as in Proposition \ref{prop:netrusov}. Then for 
$\frac{n-d}{p}+m<s_2 \leq  s_1<\frac{n-d}{p}+m+1$ it holds that
\begin{align*}
\widetilde{A}^{s_2}_{p,(q)}(\Gamma^c)\cap A^{s_1}_{p,(q)}(\Rn) = \widetilde{A}^{s_1}_{p,(q)}(\Gamma^c).
\end{align*}
In particular, 
\begin{align*}
\widetilde{H}^{s_2}_{p}(\Gamma^c)\cap H^{s_1}_{p}(\Rn) = \widetilde{H}^{s_1}_{p}(\Gamma^c).
\end{align*}
\end{cor}
\begin{proof}
It suffices to prove the inclusion $\widetilde{A}^{s_2}_{p,(q)}(\Gamma^c)\cap A^{s_1}_{p,(q)}(\Rn) \subset \widetilde{A}^{s_1}_{p,(q)}(\Gamma^c)$, since the reverse inclusion is obvious. So let $u\in \widetilde{A}^{s_2}_{p,(q)}(\Gamma^c)\cap A^{s_1}_{p,(q)}(\Rn)$, which implies that $\Tr_{\Gamma,m}|_{A_{p,(q)}^{s_2}}u=0$. Then since $A_{p,(q)}^{s_1}$ is continuously embedded in $A_{p,(q)}^{s_2}$ and $\Tr_{\Gamma,m}|_{A_{p,(q)}^{s_1}} = (\Tr_{\Gamma,m}|_{A_{p,(q)}^{s_2}})|_{A_{p,(q)}^{s_1}}$ we have that $\Tr_{\Gamma,m}|_{A_{p,(q)}^{s_1}}u = \Tr_{\Gamma,m}|_{A_{p,(q)}^{s_2}}u=0$, which by Proposition \ref{prop:netrusov} implies that $u\in \widetilde{A}^{s_1}_{p,(q)}(\Gamma^c)$.
\end{proof}

\begin{rem}\label{rem:CounterexDensity}
The condition on the regularity exponents $s_1$ and $s_2$ in Theorem~\ref{thm:maindensity2} cannot in general be dispensed with. For brevity we focus on the $H^s_{p,\Gamma}$ spaces, but the $B^s_{p,q,\Gamma}$ case is analogous.

We first recall that for any $d$-set $\Gamma$ with $0<d<n$, we have by Proposition~\ref{prop:Nullity} that if $s_1 >-\frac{n-d}{p'}>s_2$ then $H^{s_1}_{p,\Gamma}=\{0\}$ and $H^{s_2}_{p,\Gamma}\ne\{0\}$, so $H^{s_1}_{p,\Gamma}$ is not dense in $H^{s_2}_{p,\Gamma}$ in this case.
If moreover $\Gamma$ is either compact or a $d$-dimensional hyperplane, by Remark~\ref{rem:nullity} this holds also for $s_1\ge-\frac{n-d}{p'}>s_2$.

For a counterexample to density when both $H^{s_1}_{p,\Gamma}$ and $H^{s_2}_{p,\Gamma}$ are non-trivial, let $0<d<n$ be an integer and let $\Gamma=\{(x_1,\ldots,x_n):\; x_1^2+\cdots+x_d^2\le1,\;x_{d+1}=\cdots=x_n=0\}$ be the unit $d$-dimensional closed disc embedded in $\R^n$.
We shall show that, with $s_M:=\frac{n-d}{p'}+M$, the inclusion 
\begin{align}
\label{eqn:nondensity}
H^{-s_M+\epsilon}_{p,\Gamma}\subset H^{-s_M-\epsilon'}_{p,\Gamma}
\;\text{ is \emph{not dense}, for all }
\epsilon,\epsilon'>0, \;M\in\N,\; 1<p<\infty.
\end{align}
To prove \rf{eqn:nondensity} for a given $M\in\N$, $1<p<\infty$ and $\epsilon,\epsilon'>0$, it suffices to exhibit $\xi_M\in H^{-s_M-\epsilon'}_{p,\Gamma}$ and $u_M\in \cD(\R^n)$ such that $\langle\xi_M,u_M \rangle\neq 0$ but $\langle v,u_M \rangle=0$ for all $v\in H^{-s_M+\epsilon}_{p,\Gamma}$, since then 
it cannot be possible to find a sequence of elements of $v\in H^{-s_M+\epsilon}_{p,\Gamma}$ approximating $\xi_M$.
Explicitly, we define $\xi_M\in\cS'(\R^n)$ to be the ``$M$th derivative in the $n$th Cartesian coordinate of a 
$d$-dimensional delta'': 
$\xi_M(\varphi)=\int_\Gamma\frac{\partial^M\varphi}{\partial x_n^M}(\gamma)\cH^d(d\gamma)$ for $\varphi\in\cS(\R^n)$. 
Then $\xi_M\in H^{-s_M-\epsilon'}_{p,\Gamma}(\R^{n})$, e.g.\ by \cite[Prop.~A.1]{HewMoi:15}, noting that $\xi_M$ is --- up to a constant factor (depending on the Hausdorff measure normalisation) --- an $M$th distributional derivative of the tensor product ($\chi_{\{|x|\le 1\}}\otimes\delta_0$) between the characteristic function of the unit disc in $\R^d$ ($\chi_{\{|x|\le 1\}}\in L_p(\R^d)$) and a delta function in $\R^{n-d}$ ($\delta_0\in H^{-\frac{n-d}{p'}-\epsilon'}_{p}(\R^{n-d})$ by \cite[Rmk.~2.2.4.3]{RuS96}). 
Next we define $u_M\in\cD(\Rn)$ to be the cut-off polynomial $u_M(x):=\chi(x)x_n^M$, for some $\chi\in\cD(\R^n)$ taking the constant value 1 in a neighbourhood of $\Gamma$. 
Clearly $\langle\xi_M,u_M\rangle=M!\cH^d(\Gamma)>0$. 
Furthermore, $u_M\in\Aspq$ for all $0<p,q<\infty$ and all $s\in\R$, and $u_M\in\ker\Tr_{\Gamma,m}$ if and only if $m<M$, so that in particular, by Proposition~\ref{prop:netrusov}, $u_M\in \tH^{s_M-\epsilon}_{p'}(\Gamma^c)$.
But by Proposition \ref{prop:annihilators} this implies that $\langle v,u_M \rangle=0$ for all $v\in H^{-s_M+\epsilon}_{p,\Gamma}$, and hence \rf{eqn:nondensity} is proved.

We note that the function $u_M$ constructed above satisfies $u_M\in \tH^{s_M-\epsilon}_{p'}(\Gamma^c)\cap H^{s_M+\epsilon'}_{p'}(\R^n) \setminus \tH^{s_M+\epsilon'}_{p'}(\Gamma^c)$,
showing that the conditions on $s_1,s_2$ in Corollary~\ref{cor:TildeCap} also cannot in general be dispensed with.

The above analysis holds, with appropriate modifications, 
with $\Gamma$ replaced by a $d$-dimensional hyperplane or any sufficiently smooth $d$-dimensional manifold. 
However, to our knowledge, whether the conditions on the regularity exponents are as close to optimal in the case of $d$-sets with $d\notin\N_0$ is an open problem.
\end{rem}

\begin{rem}\label{rem:CounterexDensityRough}
Without the assumption that $\Gamma$ is a $d$-set, the density of the subscript spaces into one another can fail completely.
Let $d_1,d_2,\ldots$ be a sequence spanning all rational numbers in the interval $(0,n)$, and 
let $K_1,K_2,\ldots$ be a sequence of compact subsets of $[0,1]^n$ such that each $K_j$ is a $d_j$-set; an explicit construction of these sets in terms of ``Cantor dusts'' is given in \cite[Thm.~4.5]{HewMoi:15}, where $K_j$ is denoted $F^{(n)}_{d_j,\infty}$.
Define $\Gamma:=\{ 0\}\cup\bigcup_{j\in\N} 2^{-j-1} (K_j + 3e_1)$, i.e.\ the union of translated and scaled copies of all the $K_j$ in such a way that they lie at positive distance from one another. (Here the point $0$ is added to make $\Gamma$ compact, and $e_1$ is the unit vector along the first Cartesian axis.)
Fix $1<p<\infty$ and $-\frac n{p'}\le s_2<s_1\le 0$.
Then by Proposition \ref{prop:Nullity} there is $0\ne u\in H^{s_2}_{p,K_j}$ for $j$ such that $s_2<-\frac{n-d_j}{p'}<s_1$, while $H^{s_1}_{p,K_j}=\{0\}$.
Thus $\hat u$, obtained from scaling and translating $u$ to be supported in $ 2^{-j-1} (K_j + 3 e_1)$, belongs to $H^{s_2}_{p,\Gamma}$ but cannot be approximated by elements of $H^{s_1}_{p,\Gamma}$.
Hence $\Gamma$ is a compact set such that 
$H^s_{p,\Gamma}=\{0\}$ for all $s\ge0$ and $H^s_{p,\Gamma}\neq \{0\}$ for all $s<0$, and
$$
H^{s_1}_{p,\Gamma}\; \text{ is {\em not dense} in }\; H^{s_2}_{p,\Gamma}
\qquad
\text{for all } 1<p<\infty \text{ and }-\frac n{p'}\le s_2<s_1\le 0.
$$
\end{rem}

\section*{Acknowledgements}
We would like to thank Prof.\ H.~Triebel for putting the second and third named authors in contact with the first, thus inducing the fruitful collaboration which gave rise to this paper, and for some helpful preliminary discussions about the issues dealt with herein. 
We would also like to thank Prof.\ S.~Chandler-Wilde, with whom the second and third named authors have been collaborating regularly, for providing a stimulating atmosphere for the commencement of this project during a visit by all three authors to the University of Reading in April 2018.

\addcontentsline{toc}{section}{References}

\end{document}